\documentclass[a4paper,11pt,bibtotoc,idxtotoc]{scrreprt}

\usepackage{scrpage2}

\usepackage[british]{babel}

\usepackage{amsmath,amsthm,amssymb}
\swapnumbers

\usepackage[utf8x]{inputenc}

\usepackage{graphicx}
\usepackage[matrix,arrow]{xy}

\usepackage{dezert}

\titlehead{\centering\Large\sffamily
Charles University in Prague\par
\smallskip
Faculty of Mathematics and Physics\par}

\subject{Doctoral Thesis}
\title{Homomorphisms\\and Structural Properties\\of Relational Systems}
\author{Jan Foniok}
\date{}
\publishers{Department of Applied Mathematics\\[3pt]
Supervisor: Prof. RNDr. Jaroslav Nešetřil, DrSc.\\[3pt]
Branch I4: Discrete Models and Algorithms}

\begin{document}

\maketitle

\addchap{Preface}

My interest in graph homomorphisms dates back to the Spring School of
Combinatorics in 2000.  The School is traditionally organised by the
Department of Applied Mathematics of the Charles University in a hilly
part of the Czech Republic; in 2000 it was one of the few times when the
Spring School was international not only in terms of its participants,
but also in terms of its venue. Participants will never forget carrying
a blackboard across the border between Germany and the Czech Republic
and the exciting boat trip on the Vltava.

The study text on homomorphisms~\cite{Nes:SS}, specially prepared for
the Spring School, aroused my curiosity that has eventually resulted
in both my master's thesis and this doctoral dissertation.

The study of graph homomorphisms was pioneered by G.~Sabidussi,
Z.~He\-dr\-lín and A.~Pultr in the 1960's. It was part of an attempt to
develop a theory of general mathematical structures in the framework
of algebra and category theory. Many nice and important results have
emerged from their work and the work of their followers. Even so, until
recently many graph theorists would not include homomorphisms among the
topics of central interest in graph theory.

Nevertheless, graph homomorphisms and structural properties of graphs
have recently attracted much attention of the mathematical community. The
reason may be in part that the homomorphism point of view has proved
useful in various areas ranging from colouring and graph reconstruction
to applications in artificial intelligence, telecommunication, and even
statistical physics. A book~\cite{HelNes:GrH} now exists that introduces
the topic and brings together the most important parts of the theory
and its applications. This thesis surveys my small contribution to the
ongoing research in this area.

\bigskip

Some results contained in this thesis have been published or accepted
for publication.

\begin{itemize}
\frenchspacing

\item[\cite{FNT:WG06}]
J.~Foniok, J.~Ne{\v s}et{\v r}il, and C.~Tardif.
\newblock Generalised dualities and finite maximal antichains.
\newblock In F.~V. Fomin, editor, {\em Graph-Theoretic Concepts in Computer
  Science (Proceedings of WG 2006)}, volume 4271 of {\em Lecture Notes in
  Comput. Sci.}, pages 27--36. Springer-Verlag, 2006.

\item[\cite{FNT:Eurocomb7}]
J.~Foniok, J.~Ne{\v s}et{\v r}il, and C.~Tardif.
\newblock On finite maximal antichains in the homomorphism order.
\newblock {\em Electron. Notes Discrete Math.}, 29:389--396, 2007.

\item[\cite{FNT:GenDu}]
J.~Foniok, J.~Ne{\v s}et{\v r}il, and C.~Tardif.
\newblock Generalised dualities and maximal finite antichains in the
  homomorphism order of relational structures.
\newblock {\em European J. Combin.}, to appear.

\end{itemize}

\subsection*{Acknowledgements}

\vbox{
I thank my supervisor, Jaroslav Nešetřil, for introducing me to
the world of homomorphisms and for persistent support of my research.
Discussions with him have always been inspiring, but at the same time
compelling and enjoyable.

Many thanks to Claude Tardif for wonderful cooperation and for his
surrealistic attitudes. Claude's contribution to my being able to write
this thesis is indeed valuable.

For numerous suggestions I am grateful to
Manuel Bodirsky,
Julia Böttcher
and Zdeněk Hedrlín.
I also acknowledge the generous support of the Institute for Theoretical Computer
Science (ITI) in Prague and the EU Research Training Network COMBSTRU.
Last but not least I thank Ida Švejdarová for lending me a pen when
I needed one.
}

\subsection*{Disclaimer}

This is not the official version of the thesis. The printed version
uses typesetting and fonts that are not in the public domain. A better
version can be downloaded from
http://www.ifor.math.ethz.ch/staff/foniok.

\tableofcontents
\addchap{Notation}

\begin{tabular}{l@{\quad\dots\quad}p{.76\textwidth}}
$\bs A$&	base set of a \Ds~$A$ [see~\ref{def:122}]\\
$A\to B$&	$A$ is homomorphic to~$B$\\
$A\heq B$&	$A$ is homomorphically equivalent to~$B$; $A\to B$ and $B\to A$\\
$A\inc B$&	$A$ and $B$ are incomparable; $A\notto B$ and $B\notto A$\\
$A/{\approx}$&	factor structure [see~\ref{dfn:factor}]\\
$\Block(A)$&	see~\ref{dfn:incgraph}\\
$C^B$&		exponential structure [see~\ref{dfn:expon}]\\
$\CD$&		category of \Ds s; the homomorphism order of \Ds s\\
$D(F)$&		dual of the \Dt~$F$\\
$D(\F)$&	dual of a finite set~$\F$ of \Dt s\\
$\D(\F)$&	finite dual set of a finite set~$\F$ of \Df s\\
$\DSh(A)$&	directed shadow of~$A$ [see~\ref{dfn:ori-shadow}]\\
$\toD$&		$\{ A : A\to D \}$\\
$\toDd $&	$\{ A : A\to D \text{ for some } D\in\D\}$\\
$\ntoDd $&	$\{ A : A\notto D \text{ for all } D\in\D\}$\\
$\Delta$&	a type of relational structures; $\Delta=(\delta_i:i\in I)$\\
$f:A\to B$&	$f$ is a homomorphism from~$A$ to~$B$\\
$f[S]$&		image of the set~$S$ under the mapping~$f$, if $S$~is a subset of the
		domain of~$f$; $f[S]=\bigl\{f(s):s\in S\bigr\}$\\
$f(e)$&		if $e=(u_1,u_2,\dotsc,u_k)$, \vfil\goodbreak then
		$f(e)=\bigl(f(u_1),f(u_2),\dotsc,f(u_k)\bigr)$\\
$f\restriction T$&
		restriction of the function~$f$ to a subset~$T$ of the domain of~$f$\\
$\Fnto$&	$\{ A : F\notto A\}$\\
$\Ffto$&	$\{ A : F\to A \text{ for some } F\in\F\}$\\
$\Ffnto$&	$\{ A : F\notto A \text{ for all } F\in\F\}$\\
$\gimel(x)$&	height label of the vertex~$x$ [see~\ref{dfn:htlab}]\\
$I$&		set of indices; $\Delta=(\delta_i:i\in I)$\\
$\Inc(A)$&	incidence graph of~$A$ [see~\ref{dfn:incgraph}]\\
$K_k$&		complete graph on $k$ vertices\\
$p\vee q$&	supremum of $p$ and $q$; join in a lattice\\
$p\wedge q$&	infimum of $p$ and $q$; meet in a lattice\\
$p\heyt q$&	Heyting operation [see~\ref{dfn:heytal}]\\
$\vec P_k$&	directed path with $k$ edges\\
$R_i(A)$&	the $i$th edge set of a \Ds~$A$ [see~\ref{def:122}]\\
$\downs S$&	downset generated by~$S$\\
$\ups S$&	upset generated by~$S$\\
\end{tabular}

\noindent
\begin{tabular}{l@{\quad\dots\quad}p{.68\textwidth}}
$\Sh(A)$&	shadow of~$A$ [see~\ref{dfn:shadow}]\\
$\vec T_k$&	transitive tournament with $k$ vertices\\
$[a]_{\approx}$&class of the equivalence~$\approx$ that contains~$a$\\
$\top$&		\Ds\ with one vertex and all loops; $\bs\top=\{1\}$;
		$R_i(\top)={\bs\top}^{\delta_i}$ for all $i\in I$\\
$\bot$&		\Ds\ with one vertex and no edges; $\bs\bot=\{1\}$;
		$R_i(\bot)=\emptyset$ for all $i\in I$\\
$\coprod_{j\in J} A_j$; $A+B$&
		sum of relational structures [see~\ref{dfn:sum}]\\
$\prod_{j\in J} A_j$; $A\times B$&
		product of relational structures [see~\ref{dfn:prod}]
\end{tabular}

\setchapterpreamble[u]{%
	\dictum[Johann Wolfgang von~Goethe]{%
	Whatever you can do or dream you can, begin it.
	Boldness has genius, magic and power in it.
	Begin~it~now.}}

\chapter{Introduction}

\section{Motivation and overview}

In this thesis, we study homomorphisms of finite relational
structures. Finite relational structures can be viewed in several
ways. The view we adopt consists in seeing them as a generalisation
of graphs. Relational structures may actually be described as oriented
uniform hypergraphs with coloured edges. There are three main differences
from ordinary graphs: edges are ordered, they are tuples of possibly
more than two vertices, and there are various kinds of edges.

Homomorphisms are mappings between vertex sets of relational
structures. Homomorphisms preserve edges; so the image of an edge is an
edge. Moreover, it is an edge of the same kind.

Thus homomorphisms endow graphs and relational structures with an
algebraic structure that will be familiar to an algebraist or category
theorist.

The unifying concept in the thesis is the question of existence of
homomorphisms. It interconnects the two main topics presented here.

The first topic is homomorphism dualities.  There the existence of a
homomorphism between structures is equivalent to the non-existence of a
homomorphism between other structures.  In particular, we study situations
where a class of relational structures is characterised both by the
non-existence of a homomorphism from some \emph{finite} set of structures,
and by the existence of a homomorphism to some other \emph{finite} set
of structures. Such situations are called finite homomorphism dualities.
We provide a full characterisation of finite homomorphism dualities.

The other topic is the homomorphism order, where the existence of a
homomorphism defines a relation that turns out to induce a partial order on
the class of relational structures.  We examine especially finite maximal
antichains in the homomorphism order. We find a surprising correspondence
between maximal antichains and finite dualities. Many finite maximal
antichains have the splitting property; we derive a structural condition
on those antichains that do not have this property.

The main results of the thesis are the characterisation of all finite
homomorphism dualities (Theorem~\ref{thm:finite-character}) and the
splitting property of finite maximal antichains in the homomorphism
order with described exceptions (Theorem~\ref{thm:aasplit} and
Theorem~\ref{thm:mac=dual}).

Other results include a new construction of dual structures,
which generalises two previous constructions (Section~\ref{sec:3constr}).
Furthermore we extend our results on homomorphism dualities for relational
structures into the context of lattices (Section~\ref{sec:heyt}).
And finally we state several consequences of these results in the area
of computational complexity (Chapter~\ref{chap:complex}).

\section{Relational structures}

First things first. We study homomorphisms of relational structures,
so let us first define relational structures.

\begin{dfn}
\label{dfn:structure}
A \deph{type}\index{type $\Delta$} $\Delta$ is a sequence $(\delta_i:
i\in I)$ of positive integers; $I$~is a finite set of indices. A
(finite) \deph{relational structure}\index{relational structure} $A$
of type $\Delta$ is a pair $\bigl(X, (R_i: i\in I)\bigr)$, where $X$
is a finite nonempty set and $R_i\subseteq X^{\delta_i}$; that is,
$R_i$ is a $\delta_i$-ary relation on~$X$.  Relational structures of
type~$\Delta$ are denoted by capital letters $A$, $B$, $C$,~\dots
\end{dfn}

\begin{pgf}
There are many natural examples of relational structures. Perhaps the simplest
are digraphs (with loops allowed), which are simply \Ds s of type $\Delta=(2)$.
This example is also the motivation for our terminology.
The class of
all partially ordered sets is a subclass of the class of all \Ds s
for $\Delta=(2)$, requiring that the relation be reflexive, transitive
and antisymmetric.

The class of all $(2)$-structures whose relation is symmetric is the
class of all undirected graphs. This is an important example and we
shall keep it in mind as we define properties and operations on \Ds s;
all of them carry over immediately to undirected graphs.
For an undirected graph~$G$, the base set is usually called the \emph{vertex set} and denoted
by~$V(G)$.
\end{pgf}

\begin{pgf}
In a logician's words, relational structures are models of theories with no function
symbols; and finite relational structures are such models in the theory of finite sets.
\end{pgf}

\begin{dfn}
\label{def:122}
If $A=\bigl(X,(R_i:i\in I)\bigr)$, the \deph{base set}~$X$ is denoted by~$\bs{A}$
and the relation~$R_i$ by~$R_i(A)$. We often refer to a relational
structure of type~$\Delta$ as \deph{\Ds}.\index{$\Delta$-structure} The type~$\Delta$
is almost always fixed in the following text. The
elements of the base set are called \deph{vertices} and the elements of
the relations~$R_i$ are called \deph{edges}\index{edge!of relational structure}.
For the set of all edges of a \Ds~$A$ we use the notation~$R(A)$, that is
\[R(A) := \bigcup_{i\in I} R_i(A).\]
To distinguish between various relations of a \Ds\ we speak about
\deph{kinds of edges} (so the elements of~$R_i(A)$ are referred to as
the \emph{edges of the $i$th kind}).
\end{dfn}

\begin{unary}
Some of the relations of \Ds s may be unary (this is the same as saying
that some of the numbers~$\delta_i$ may be equal to one). In this thesis,
however, we consider only relational structures with \emph{no} unary
relations. This is a rather technical assumption.  All the results and
proofs remain valid even for structures with unary relations, but some
adjustments would have to be made.

For example, in Theorem~\ref{thm:mac=dual} we suppose that there are
at most two relations. In fact, we should suppose that there are at
most two relations of arity greater than one and an arbitrary number of
unary relations.

Elsewhere the statements would have to be slightly altered, like in
Theorem~\ref{thm:dualsize}.  If we allow unary relations, the upper
bound has to be replaced by $n^{n+|I|}$.

Such adjustments would, in our opinion, make the text less clear and
more difficult to read, so we find it useful to assume that structures
have no unary relations.
\end{unary}

Substructures of relational structures are what we are familiar with
from graph theory as \emph{induced subgraphs}. It is possible to define
``non-induced'' substructures as well, but in our context it is more
convenient not to do so.

\begin{dfn}
A \deph{substructure} of a \Ds~$A$ is any structure $\bigl(S,(R'_i:
i\in I)\bigr)$ such that $S$~is a subset of~$\bs A$ and $R'_i=R_i
\cap S^{\delta_i}$. It is also called the \deph{substructure of~$A$
induced by~$S$}; and it is a \deph{proper substructure} if $S\subsetneqq\bs A$.
\end{dfn}

Next we introduce a factorisation construction.  This construction can
be viewed as gluing equivalent vertices together.  An example of the
construction (for digraphs) is in Figure~\ref{fig:factor}.

\begin{dfn}\label{dfn:factor}
Let $A$ be a \Ds\ and let~$\approx$~be an equivalence relation
on~$\bs A$. We define the \deph{factor structure}~$A/{\approx}$ to
be the \Ds\ whose base set is the set of all equivalence classes of
the relation~$\approx$, so $\bs{A/{\approx}} = \bs A/{\approx}$, and
a $\delta_i$-tuple of equivalence classes is in the relation~$R_i$
of~$A/{\approx}$ if we can find an element in each of the classes
such that these elements form an edge of~$A$,
\[
R_i(A/{\approx}) =
	\Bigl\{\bigl([a_1]_\approx, [a_2]_\approx,\dotsc,[a_{\delta_i}]_\approx\bigr) :
	(a_1,a_2,\dotsc,a_{\delta_i})\in R_i(A)\Bigr\},\ i\in I.
\]
\end{dfn}

\picture{intro.3}{Factor structure}{factor}

Finally, let us once more stress the two restrictions posed on relational structures
in this thesis: All structures we consider are \emph{finite} and have
\emph{no unary relations}.

\section{Homomorphisms}

Familiar with the notion of relational structures, we continue by defining
homomorphisms. As one might expect, they are mappings of base sets that
preserve all the relations.

\begin{dfn}
Let $A$ and $A'$ be two relational structures of the same type~$\Delta$. A
mapping $f: \bs A\to\bs{A'}$ is a \deph{homomorphism}\index{homomorphism}
from~$A$ to~$A'$ if for every $i\in I$ and for every $u_1, u_2,
\dotsc, u_{\delta_i}\in \bs A$ the following implication holds:
\[(u_1, u_2, \dotsc, u_{\delta_i})\in R_i(A) \quad\Rightarrow\quad
\bigl(f(u_1), f(u_2), \dotsc, f(u_{\delta_i})\bigr)\in R_i(A').\]
\end{dfn}

\begin{dfn}
The fact that
$f$ is a homomorphism from~$A$ to~$A'$ is denoted by~\[f:A\to A'.\]
If there exists a homomorphism from~$A$ to~$A'$, we say that $A$~is
\deph{homomorphic}\index{homomorphic} to~$A'$ and write
$A\to A'$; otherwise we write $A\notto A'$. If $A$~is homomorphic to~$A'$
and at the same time $A'$~is homomorphic to~$A$, we say that $A$ and $A'$
are \deph{homomorphically equivalent}\index{homomorphically equivalent}
and write $A\heq A'$.  If on the other hand there
exists no homomorphism from~$A$ to~$A'$ and no homomorphism from~$A'$
to~$A$, we say that $A$ and~$A'$ are \deph{incomparable} and write $A\inc A'$.
\end{dfn}

\begin{dfn}
A homomorphism from~$A$ to itself
is called an \deph{endomorphism}\index{endomorphism} of~$A$.
\end{dfn}

Next we give two simple examples of homomorphisms.

\begin{exa}
This is a trivial example: Let $\Delta=(2)$, so we consider digraphs. Let
$\vec P_k$ be the directed path on vertices $\{0,1,\dotsc,k-1\}$ with edge
set $\bigl\{(j,j+1):j=0,1,\dotsc,k-1\bigr\}$. Then $f:V(\vec P_k)\to V(G)$
is a homomorphism if and only if $f(0),f(1),\dotsc,f(k)$ is a directed
walk in~$G$.
\end{exa}

\begin{exa}
\label{exa:k-colouring}
For an undirected graph~$G$, a homomorphism to the complete graph~$K_k$ is essentially a
$k$-colouring of~$G$: imagine the vertices of~$K_k$ as colours, and since edges of~$G$ have to be
preserved, distinct vertices of~$G$ have to be mapped to distinct vertices of~$K_k$, in other
words they have to be assigned distinct colours.
\end{exa}

\begin{compohom}
It is a very important aspect of homomorphisms that they compose~--
the composition of two homomorphisms is a homomorphism as well. This
composition operation endows a set of \Ds s with a structure of an
algebraic flavour. Homomorphisms of relational structures share this
property with morphisms of other structures, like topological spaces,
semigroups, monoids, partial orders and many others. This flavour is discussed in
more detail in Section~\ref{sec:category} but it is omnipresent
throughout the thesis.
\end{compohom}

\begin{dfn}
As usual, a homomorphism from~$A$ to~$A'$ is an
\deph{isomorphism}\index{isomorphism} if it is a bijection and $f^{-1}$ is
a homomorphism of~$A'$ to~$A$. If there exists an isomorphism between~$A$
and~$A'$, we say that $A$ and $A'$ are \deph{isomorphic}\index{isomorphic}
and write $A\isom A'$\index{$A\isom A'$}. An isomorphism of $A$ with
itself is called an \deph{automorphism}\index{automorphism} of~$A$.
\end{dfn}

\section{Retracts and cores}

In this section, we introduce the notion of \emph{cores}. Cores are \Ds
s that are minimal in the following sense: A core is not homomorphically
equivalent to any smaller structure.

For the formal definition, we use retractions. A retraction is an
endomorphism that does not move any vertex in its image.

\begin{dfn}
Let $A$ be a \Ds. An endomorphism $f:A\to A$ is a
\deph{retraction}\index{retraction} if it leaves its image fixed, in
other words if $f(x)=x$ for all $x\in f[\bs A]$.  A substructure $B$
of~$A$ is called a \deph{retract}\index{retract} of~$A$ if there exists
a retraction of~$A$ onto~$B$; a retract is \deph{proper} if it is a proper
substructure.
\end{dfn}

Later we will need the fact that a retract is homomorphically equivalent
to the original structure.

\begin{lemma}
\label{lem:reteq}
If $B$ is a retract of~$A$, then $A$ and $B$ are homomorphically equivalent.
\end{lemma}

\begin{proof}
If $B$ is a retract of $A$, then $B$~is a substructure of~$A$ and so the
identity mapping is a homomorphism from~$B$ to~$A$. On the other hand,
the retraction is a homomorphism from~$A$ to~$B$.
\end{proof}

\begin{dfn}
A \Ds\ $C$ is called a \deph{core}\index{core} if it has no
proper retracts. A retract $C$ of~$A$ is called a \deph{core of~$A$}
if it is a core.
\end{dfn}

Several other conditions are equivalent to the one we chose for the
definition of a core.

\begin{lemma}[Characterisation of cores]
For a \Ds~$C$ the following conditions are equivalent.
\begin{itemize}
\item[(1)] $C$ is a core (that is, $C$ has no proper retracts).
\item[(2)] $C$ is not homomorphic to any proper substructure of~$C$.
\item[(3)] Every endomorphism of~$C$ is an automorphism.
\end{itemize}
\end{lemma}

\begin{proof}
(1)${}\Rightarrow{}$(2):
Suppose that $f$~is a homomorphism of~$C$ to a proper substructure
of~$C$. Then $f$~is a permutation of its image~$f[\bs C]$ and so there
exists a positive integer~$k$ such that $f^k$~restricted to~$f[\bs
C]$ is the identity mapping. Then $f^k$~is a retraction of~$C$ onto a
proper retract.

(2)${}\Rightarrow{}$(3): Let $f: C\to C$ be
an endomorphism. The mapping~$f$ is surjective, so it is a bijection
($\bs C$~is finite). There exists a positive integer~$k$ such that $f^k$
is the identity mapping. Then $f^{-1}=f^{k-1}$ is a homomorphism, so
$f$~is an automorphism.

(3)${}\Rightarrow{}$(1):
Every retraction is an endomorphism, so every retraction of~$C$ is an
automorphism. Thus the image of every retraction of~$C$ is~$C$, hence
$C$~has no proper retracts.
\end{proof}

Next we prove that every relational structure has exactly one core
(up to isomorphism). For the proof we use two lemmas.

\begin{lemma}\label{lemma:suriso}
Let $A$ and $B$ be two \Ds s. If there exist surjective
homomorphisms $f:A\to B$ and $g:B\to A$, then $A$ and $B$ are isomorphic.
\end{lemma}

\begin{proof}
The existence of $f$ shows that $|\bs A|\ge|\bs B|$ and the existence
of~$g$ shows $|\bs A|\le|\bs B|$, so $|\bs A|=|\bs B|$. This means
that $f$ is a bijection and that $g\circ f$ is a bijection, so there
exists a positive integer $k$ such that $(g\circ f)^k$ is the identity
mapping. Then $(g\circ f)^{k-1}\circ g = f^{-1}$, therefore $f^{-1}$
is a homomorphism and $f$ an isomorphism.
\end{proof}

\begin{lemma}\label{lemma:coreveryuniq}
Let $C$ and $C'$ be two cores. If $C$ and $C'$ are homomorphically equivalent,
they are isomorphic.
\end{lemma}

\begin{proof}
Let $f: C\to C'$ and $g: C'\to C$. The mapping $f\circ g$ is an
endomorphism of~$C'$. Since $C'$ is a core, it is an automorphism,
therefore both $f$ and $g$ are surjective. The existence of an
isomorphism follows from Lemma~\ref{lemma:suriso}.
\end{proof}

This lemma implies in particular that every \Ds\ is homomorphically equivalent
to at most one core.  The next proposition asserts that it is actually
\emph{exactly} one.

\begin{prop}\label{prop:allhavecore}
Every \Ds\ $A$ has a unique core~$C$ (up to isomorphism).
Moreover, $C$~is the unique core to which $A$~is homomorphically
equivalent.
\end{prop}

\begin{proof}
Proof of existence: select the retract~$C$ of~$A$ of the smallest
size. Any potential proper retract~$C'$ of~$C$ would be a smaller retract
of~$A$, and so $C$ is a core.

Let $C$ and $C'$ be two distinct cores of~$G$. By the definition of a
core, there exist homomorphisms $f: G\to C$ and $f': G\to C'$. The
restrictions $f\restriction C'$ and $f'\restriction C$ show that $C$
and $C'$ are homomorphically equivalent; by Lemma~\ref{lemma:coreveryuniq}
they are isomorphic.
\end{proof}

\begin{cor}
A \Ds~$C$ is a core if and only if it is not homomorphically equivalent
to a \Ds\ with fewer vertices.
\end{cor}

\begin{proof}
Suppose that $C$ and $D$ are homomorphically equivalent and $D$~has
fewer vertices than~$C$.  Then the core of~$D$ has at most~$|\bs D|$
vertices, but $|\bs D|<|\bs C|$, so $C$~is not the unique core to which
$D$~is homomorphically equivalent. Hence $C$~is not a core.

Conversely, if $C$~is not a core, then it is homomorphically equivalent
to its core~$C'$. The core~$C'$ is a proper retract of~$C$, hence it
has fewer vertices.
\end{proof}

Thus we can summarise four equivalent definitions of a core.

\begin{cor}[Characterisation of cores revisited]
For a \Ds~$C$ the following conditions are equivalent.
\begin{itemize}
\item[(1)] $C$ is a core (that is, $C$ has no proper retracts).
\item[(2)] $C$ is not homomorphic to any proper substructure of~$C$.
\item[(3)] Every endomorphism of~$C$ is an automorphism.
\item[(4)] $C$~is not homomorphically equivalent
to a \Ds\ with fewer vertices.
\qed
\end{itemize}
\end{cor}

\section{The category of relational structures and homomorphisms}
\label{sec:category}

It is little surprising that structures of an algebraic nature,
like \Ds s, with suitably selected mappings among them, should form
a category.  Here, we observe some basic properties of the category of
\Ds s and homomorphisms. The reader may consult \cite{AHS:AbsConCat}
or~\cite{BarWel:Cat} for an introduction to category theory.

\begin{pgf}
\label{pgf:cat}
For a fixed type~$\Delta=(\delta_i : i\in I)$, let $\CD$ be the category
of all \Ds s (objects) and their homomorphisms (morphisms).
\end{pgf}

\begin{dfn}
\label{dfn:sum}
Let $J$ be a nonempty finite index set and let $A_j$, $j\in J$ be \Ds s.
We define the \deph{sum} $\coprod_{j\in J}A_j$\index{sum of \Ds s} to be
the \emph{disjoint union} of the structures~$A_j$; formally the base set of
the sum is defined by
\[
\bs{\coprod_{j\in J}A_j} = \bigcup_{j\in J} \left(\{j\}\times \bs{A_j}\right),\\
\]
and the relations by
\begin{multline*}
R_i\left(\coprod_{j\in J}A_j\right) = \bigcup_{j\in J}
	\Bigl\{\bigl( (j,x_1), (j,x_2), \dotsc, (j,x_{\delta_i}) \bigr) :\\
	(x_1, x_2, \dotsc, x_{\delta_i}) \in R_i(A_j) \Bigr\},
	\quad i\in I.
\end{multline*}
\end{dfn}

\begin{prop}
The \Ds~$A=\coprod_{j\in J}A_j$ with embeddings $\iota_j: A_j\to A$ such that $\iota_j:
x\mapsto (j,x)$ for $j\in J$ is the coproduct of the \Ds s $A_j$ in the category~$\CD$.
\end{prop}

\begin{proof}
Clearly the embeddings $\iota_j$ are homomorphisms. Suppose we have another structure $A'$ and
homomorphisms $g_j: A_j\to A'$. Then $f: (j,x)\mapsto g_j(x)$ is the unique homomorphism
from~$A$ to~$A'$ such that $f\iota_j=g_j$ for all $j$ in~$J$.
\end{proof}

\begin{cor}
\label{cor:sum}
The sum~$\coprod_{j\in J}A_j$ is homomorphic to a \Ds~$B$ if and only if $A_j\to B$ for every
$j\in J$.
\end{cor}

For the sum of two \Ds s, we use the notation $A+B$, or more generally,
we can write $A_1+A_2+\dotsb+A_n$ for the sum of~$n$~structures.

We remark here that the sum of graphs, as defined above, is the usual operation of disjoint
union of graphs.

We go on to describe product in the category of \Ds s.

\begin{dfn}
\label{dfn:prod}
Let $J=\{1,2,\dotsc,n\}$ be an index set and let $A_j$, $j\in J$ be \Ds s.
We define the \deph{product}~$\prod_{j\in J}A_j=A$\index{product of \Ds
s}\index{$\prod_{j\in J}A_j$} to be the \Ds\ whose base set is
the Cartesian product of the vertex sets of the factors, and there is an
edge of a kind if and only if there is an edge of the same kind in each
projection. Formally,
\begin{align*}
\bs{A} &= \prod_{j\in J} \bs{A_j},\\
R_i\left(A\right) &=
	\Bigl\{\bigl( (x_{1,1}, x_{1,2}, \dotsc, x_{1,n}),
	(x_{2,1}, x_{2,2}, \dotsc, x_{2,n}),\dotsc,
	(x_{\delta_i,1}, x_{\delta_i,2}, \dotsc, x_{\delta_i,n}) \bigr) :\\
	&\qquad(x_{1,j}, x_{2,j}, \dotsc, x_{\delta_i,j})
	\in R_i(A_j) \text{ for all $j\in J$}\Bigr\},
	\qquad i\in I.
\end{align*}
\end{dfn}

\begin{exa}
Let $\Delta=(2,2)$. Let $A$ and $B$ be the \Ds s depicted
in Figure~\ref{fig:2ds}.  Figure~\ref{fig:prod} shows the
product~$A\times B$.
\picture{intro.5}{The \Ds s $A$ and $B$}{2ds}
\picture{intro.4}{The product $A\times B$}{prod}
\end{exa}

\begin{prop}\label{prop:products}
The \Ds~$A=\prod_{j\in J}A_j$ with projections $\pi_j: A\to A_j$ such that $\pi_j:
(x_1,x_2,\dotsc,x_n)\mapsto x_j$ for $j\in J$ is the product of the \Ds s $A_j$
in the category~$\CD$.
\end{prop}

\begin{proof}
The projections are indeed homomorphisms; and whenever $A'$~is a \Ds\ such that $g_j:A'\to A_j$
are homomorphisms, then $f: x\mapsto\bigl(g_1(x),g_2(x),\dotsc,g_n(x)\bigr)$ is the unique
homomorphism from~$A'$ to~$A$ such that $\pi_j f=g_j$ for all $j$ in~$J$.
\end{proof}

\begin{cor}
\label{cor:product}
A \Ds~$B$ is homomorphic to the product~$\prod_{j\in J}A_j$ if and only if $B\to A_j$ for every
$j\in J$.
\end{cor}

Analogously to sums, we use the convenient notation $A\times B$ and $A_1\times
A_2\times \dotsb \times A_n$ for products.

Finally, we introduce exponentiation in~$\CD$. The definition is somewhat
technical, but exponentiation is important in the context of homomorphism
dualities.

\begin{dfn}
\label{dfn:expon}
Let $B$ and $C$ be two \Ds
s. We define the \deph{exponential structure}~$C^B$\index{exponential
structure}\index{$C^B$} to be the \Ds\ whose base set
is \[\bs{C^B}=\{f: \text{$f$ is a mapping from $\bs B$ to $\bs C$}\},\] and the
$i$th relation  is the set of all
$\delta_i$-tuples $(f_1,f_2,\dotsc,f_{\delta_i})$ such that whenever
the $\delta_i$-tuple
$(b_1,b_2,\dotsc,b_{\delta_i})$ is an element of~$R_i(B)$, then
\[\bigl(f_1(b_1),f_2(b_2),\dotsc,f_{\delta_i}(b_{\delta_i})\bigr)
\in R_i(C).\]
\end{dfn}

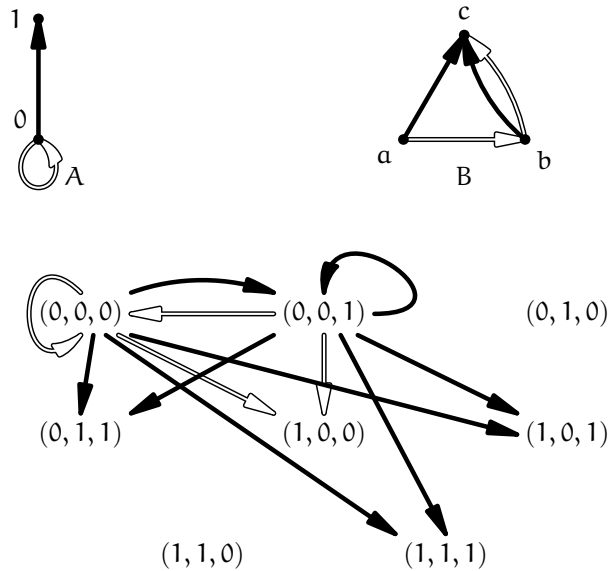
\picture{intro.1}{An example of an exponential structure}{exponent}
\begin{exa}
An example of two \Ds s $A$ and $B$ and the exponential structure~$A^B$ for
$\Delta=(2,2)$ is shown in Fig.~\ref{fig:exponent}. A mapping $f:\bs B\to
\bs A$ is represented by the triple $\bigl(f(a),f(b),f(c)\bigr)$. The
existence of a unique black edge in~$A$ means that for two functions
$f_1,f_2:\bs B\to \bs A$ to be connected with an edge in~$A^B$, it must
hold that $f_1(a)=f_1(b)=0$ (the initial vertices of all black edges
in~$B$ must be mapped to the initial vertex of the only black edge
in~$A$) and $f_2(c)=1$ (a similar condition for the terminal vertices
of the black edges).  Similarly, for outlined edges, $(f_1,f_2)$ is an
outlined edge in~$A^B$ if and only if $f_1(a)=f_1(b)=f_2(b)=f_2(c)=0$.
\end{exa}

\begin{prop}\label{prop:exponential}
For \Ds s $B$ and $C$, the exponential structure $C^B$ together with the homomorphism
$\eval:C^B\times B\to C$ defined by $\eval(f,b):=f(b)$ is an exponential object in the
category~$\CD$.
\end{prop}

\begin{proof}

If
$\bigl((f_1,b_1),(f_2,b_2),\dotsc,(f_{\delta_i},b_{\delta_i})\bigr)\in R_i(C^B\times B)$,
then
\[(f_1,f_2,\dotsc,f_{\delta_i})\in R_i(C^B)\quad\text{and}\quad
	(b_1,b_2,\dotsc,b_{\delta_i})\in R_i(B)\]
by the definition of product. It follows from the definition of edges
in the exponential structure that
\[\bigl(f_1(b_1),f_2(b_2),\dotsc,f_{\delta_i}(b_{\delta_i})\bigr)\in R_i(C).\]
So the mapping $\eval$ is indeed a homomorphism.

Let $A$ be a structure and $g:A\times B\to C$ a homomorphism. Define $\lambda g:
\bs{A}\to\bs{C^B}$ by setting $\lambda g(a)(b):=g(a,b)$. Suppose $(a_1,a_2,\dots,a_{\delta_i})\in
R_i(A)$. Now, if $(b_1,b_2,\dotsc,b_{\delta_i})\in R_i(B)$, then
\begin{multline*}
\bigl(\lambda
g(a_1)(b_1), \lambda g(a_2)(b_2), \dotsc, \lambda g(a_{\delta_i})(b_{\delta_i})\bigr) =\\
\bigl(g(a_1,b_1),g(a_2,b_2),\dotsc,g(a_{\delta_i},b_{\delta_i})\bigr)\in R_i(C)
\end{multline*}
because $g$~is a homomorphism. Therefore $\bigl(\lambda g(a_1),\lambda
g(a_2),\dotsc,\lambda g(a_{\delta_i})\bigr)\in R_i(C^B)$ and $\lambda g$
is a homomorphism from~$A$ to~$C^B$; it is easy to check that it is
the only such homomorphism that satisfies $\eval\circ(\lambda g\times
\id_B)=g$, that is, that the following diagram commutes.
\[\xy
\xymatrix{
A \ar@{{}{--}>}[dd]_{\lambda g}
	&&A\times B \ar@{{}{--}>}[dd]_{\lambda g \times \id_B} \ar[ddrr]^g\\ \\
C^B
	&&C^B \times B \ar[rr]_\eval&& C
}
\endxy\]
\end{proof}

\begin{cor}
\label{cor:expon}
For any \Ds s $A$, $B$ and $C$,
\[A\to C^B \quad\text{if and only if}\quad A\times B \to C.\]
\end{cor}

\begin{proof}
If $g:A\times B\to C$, then $\lambda g$ defined in the proof
of Proposition~\ref{prop:exponential} is a homomorphism from~$A$
to~$C^B$. Conversely, if $h:A\to C^B$, then $\eval\circ(h\times\id_B)$
is a homomorphism from~$A\times B$ to~$C$.
\end{proof}

\begin{pgf}
The base set of $C^B$ consists of all mappings from~$\bs B$ to $\bs C$, but not all of them are
homomorphisms. It follows immediately from the definition of relations of the exponential
structure, that homomorphisms from~$B$ to~$C$ are exactly those elements $f$ of~$\bs{C^B}$, for
which the $\delta_i$-tuple $(f,f,\dotsc,f)$ is in~$R_i(C^B)$ for all $i\in I$.
\end{pgf}

\begin{prop}
The category $\CD$ of \Ds s and homomorphisms is Cartesian closed.
\end{prop}

\begin{proof}
Let $V=\{1\}$ and let $\top$ be the \Ds\ such that
\begin{align*}
\bs\top &= V \quad\text{and}\\
R_i(\top) &= V^{\delta_i} \quad\text{for $i\in I$}.
\end{align*}
There exists exactly one homomorphism from any \Ds\ to~$\top$, namely the
constant mapping to~$1$. Hence $\top$ is the terminal object of~$\CD$.
Using Proposition~\ref{prop:products}, we see that there are finite
products in~$\CD$; and by Proposition~\ref{prop:exponential}, there are
exponential objects.
\end{proof}

\section{Connectedness and irreducibility}

In this section, we define connected relational structures. Our notion of
connectedness generalises weak connectedness of digraphs. We define it in
two principally different ways~-- using auxiliary undirected multigraphs
(\emph{shadows} and \emph{incidence graphs}), and by specifying structural
conditions. We show that these two ways are equivalent.

Irreducibility is a dual notion to connectedness in the category~$\CD$. We
mention a famous problem connected to irreducibility: Hedetniemi's
product conjecture.

\begin{dfn}\label{dfn:ori-shadow}
The \deph{directed shadow}\index{directed shadow} of a \Ds~$A$ is the directed
multigraph $\DSh(A)$ whose vertices are the elements
of~$\bs A$ and there is one edge from~$a$ to~$b$ for each occurrence of
the vertices~$a,b$ in an edge in some~$R_i(A)$ of arity~$\delta_i\ge2$
such that $(a_1, \dots, a_{\delta_i})\in R_i(A)$ with $a_j=a$, $a_{j+1}=b$
for some~$1\le j<\delta_i$.
\end{dfn}

\begin{dfn}\label{dfn:shadow}
The \deph{shadow}\index{shadow} of a \Ds~$A$ is the undirected
multigraph~$\Sh(A)$\index{$\Sh(A)$} that is created from~$\DSh(A)$ by
replacing every directed edge with an undirected edge (the symmetrisation
of~$\DSh$).
\end{dfn}

\begin{exa}
\label{exa:shadow}
Let $\Delta=(2,3)$, let $A$ be a \Ds,
\[A=\Bigl(\{1,2,\dotsc,6\},\bigl(\{(3,2),(6,3),(6,5)\},
\{(1,5,6),(4,4,1),(4,5,2)\}\bigr)\Bigr).\]
The the directed shadow~$\DSh(A)$ and the shadow~$\Sh(A)$ of the \Ds~$A$
are shown in Fig.~\ref{fig:shadow}. The loop at the vertex~$4$ is caused
by the triple $(4,4,1)\in R_3(A)$.
\picture{intro.2}{The directed shadow and the shadow of a \Ds}{shadow}
\end{exa}

Shadows ``preserve homomorphisms''~-- a homomorphism of \Ds s is also
a homomorphism of their shadows. The converse, however, is not true
in general.

\begin{lemma}\label{lem:sh-preserves}
If $f:A\to B$ is a homomorphism of \Ds s, then $f$~is a graph homomorphism
from~$\Sh(A)$ to~$\Sh(B)$.
\end{lemma}

\begin{proof}
If $\{u,v\}$ is an edge of~$\Sh(A)$, then by definition there is an edge $e\in R_i(A)$ for
some~$i\in I$ such that $u$ and~$v$ appear as consecutive vertices in~$e$. Therefore $f(u)$
and~$f(v)$ appear as consecutive vertices in the edge~$f(e)$ of~$B$, hence $\{u,v\}$ is an edge
of~$\Sh(B)$.
\end{proof}

\begin{dfn}\label{dfn:incgraph}
The \deph{incidence graph}\index{incidence
graph}~$\Inc(A)$\index{$\Inc(A)$} of a \Ds~A is the bipartite multigraph
$(V_1\cup V_2, E)$ with parts $V_1=\bs A$ and
\[V_2=\Block(A):=\Bigl\{\bigl(i,(a_1,\dots,a_{\delta_i})\bigr) : i\in I,\
(a_1,\dots,a_{\delta_i})\in R_i(A)\Bigr\},\]
and one edge between $a$ and $\bigl(i,(a_1,\dots,a_{\delta_i})\bigr)$
for each occurrence of~$a$ as some $a_j$ in an edge
$(a_1,\dots,a_{\delta_i})\in R_i(A)$.
\end{dfn}

\begin{exa}
The incidence graph~$\Inc(A)$ of the \Ds~$A$ from Example~\ref{exa:shadow} is shown in
Figure~\ref{fig:inci}.
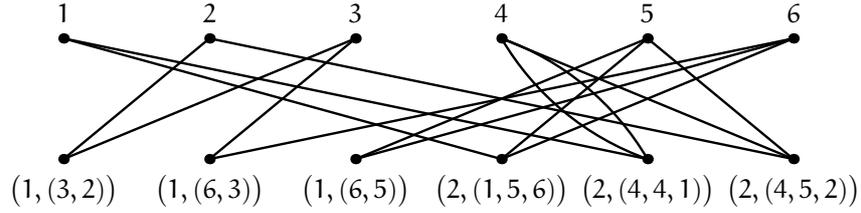
\picture{intro.6}{The incidence graph}{inci}
\end{exa}

Next we formulate three structural conditions and show that they are
equivalent with each other as well as with the connectedness of the
shadow and the incidence graph.

\begin{lemma}\label{lem:connected}
For a core~$G$, the following are equivalent:

(1) If $G\to A+B$ for some structures $A$, $B$, then $G\to A$ or $G\to B$.

(2) If $G\sim A+B$ for some structures $A$, $B$, then $B\to A$ or $A\to B$.

(3) If $G\sim A+B$ for some structures $A$, $B$, then $G\sim A$ or $G\sim B$.

(4) The shadow $\Sh(G)$ is connected.

(5) The incidence graph $\Inc(G)$ is connected.
\end{lemma}

\begin{proof}\mbox{}

\mbox{(1)${}\Rightarrow{}$(2):}
If $G\sim A+B$, then $G\to A+B$, and using (1) we have $G\to A$ or
$G\to B$. In the first case $B\to A+B\sim G\to A$, hence $B\to A$. In
the latter case $A\to A+B\sim G\to B$, and so $A\to B$.

\mbox{(2)${}\Rightarrow{}$(3):}
Suppose $G\sim A+B$. By (2) we have $B\to A$, and therefore $A\sim
A+B\sim G$; or we have $A\to B$, and then $B\sim A+B\sim G$.

\mbox{(3)${}\Rightarrow{}$(1):}
Let $G\to A+B$. Using distributivity, we have $(G\times A)+(G\times B)\sim G \times (G+A)
\times (G+B) \times (A+B)\sim G$, since $G$ is homomorphic to all other factors.

\mbox{(3)${}\Rightarrow{}$(5):}
Suppose that (3) holds but $\Inc(G)$~is disconnected. Let $A'$ be a component of~$\Inc(G)$ and
$B'=\Inc(G)-A'$; let~$A$~be the substructure of~$G$ induced by~$V(A')\cap\bs G$ and let~$B$~be
the substructure of~$G$ induced by~$V(B')\cap\bs G$. Then $G=A+B$ but both $A$ and~$B$ are
proper substructures of~$G$, so if $G\heq A$ or $G\heq B$, then $G$~is not a core, a
contradiction.

\mbox{(5)${}\Rightarrow{}$(4):}
In $\Inc(G)$, if any two vertices $u,v\in\bs G$ have a common neighbour~$e$ in $\Block(G)$,
they belong to the same edge~$e$ of~$G$, and so there is a path from~$u$ to~$v$ in~$\Sh(G)$.
Therefore the existence of a path from~$u$ to~$v$ in~$\Inc(G)$ implies the existence of a path
from~$u$ to~$v$ in~$\Sh(G)$.

\mbox{(4)${}\Rightarrow{}$(1):}
First observe that since every edge of $A+B$ is either an edge of~$A$
or an edge of~$B$, we have that $\Sh(A+B)=\Sh(A)+\Sh(B)$. Let
$f:G\to A+B$. Then $f:\Sh(G)\to\Sh(A+B)=\Sh(A)+\Sh(B)$ by
Lemma~\ref{lem:sh-preserves}. Because $\Sh(G)$ is connected, $f[\bs
G]\subseteq\bs A$ or $f[\bs G]\subseteq\bs B$: otherwise there is an
edge~$\{u,v\}$ of~$\Sh(G)$ such that $f(u)$ is a vertex of~$\Sh(A)$ and
$f(v)$~is a vertex of~$\Sh(B)$, but then $\{f(u),f(v)\}$ is not an edge of
$\Sh(A)+\Sh(B)$, a contradiction with $f$ being a homomorphism. Therefore
$f:G\to A$ or $f:G\to B$.
\end{proof}

\begin{pgf}
The conditions (1)--(3) above are equivalent even for structures that are
not cores; if a \Ds~$G$~satisfies (1)--(3), its core satisfies all the
conditions (1)--(5).  Similarly, the conditions (4)--(5) are equivalent
for all structures~$G$.
\end{pgf}

\begin{dfn}
A \Ds\ is called \deph{connected} if it satisfies the equivalent
conditions (4)--(5) of Lemma~\ref{lem:connected}.
Maximal connected substructures of a \Ds~$A$ are called the \deph{components} of~$A$.
\end{dfn}

\begin{pgf}
It is easy to see that every \Ds\ is the sum of its components; and that the
decomposition into components is unique.
\end{pgf}

A part of the previous lemma holds in the dual category too.

\begin{lemma}\label{lem:irreducible}
For a structure $G$, the following are equivalent:

(1) If $A\times B\to G$ for some structures $A$, $B$, then $A\to G$ or $B\to G$.

(2) If $G\sim A\times B$ for some structures $A$, $B$, then $B\to A$ or $A\to B$.

(3) If $G\sim A\times B$ for some structures $A$, $B$, then $G\sim A$ or $G\sim B$.
\end{lemma}

\begin{proof}
Repeat the proof of Lemma~\ref{lem:connected}; reverse all arrows and replace $+$ with~$\times$.
\end{proof}

\begin{dfn}
A \Ds\ is called \deph{irreducible} if it satisfies the equivalent
conditions of Lemma~\ref{lem:irreducible}.
\end{dfn}

\begin{pgf}
One might expect to find conditions involving $\Inc(G)$ and $\Sh(G)$, similar
to (4) and (5) of Lemma~\ref{lem:connected} to be equivalent with
irreducibility of a \Ds. But in spite of being dual to connectedness,
irreducibility is much more tricky. For example, it is an easy exercise
to devise an efficient algorithm for testing connectedness of a \Ds. On
the other hand,  irreducibility is not even known to be decidable.
\end{pgf}

\begin{pgf}
\label{pgf:nondecomposable}
Another difference from connectedness is that it is not true that every
\Ds\ can be decomposed as a finite product of irreducible structures
(dually to the decomposition into connected components).
\end{pgf}

Here it is worthwhile to mention that our terminology is
inspired by algebra, and lattice theory in particular (see
Section~\ref{sec:introrder}). In the past, other names for this property
have been used; namely \emph{productive} in~\cite{NesPul:SubFac} and
\emph{multiplicative} in~\cite{HHMN:Mult} and by many other authors. We
believe that the word \emph{irreducible} fits the meaning of the property
better.

The term \emph{multiplicative} is motivated by a property of irreducible
structures described in the next paragraphs.

\subsection*{Hedetniemi's product conjecture}

\begin{dfn}
Let $\A$ be a class of \Ds s that is closed under homomorphic
equivalence. We say that $\A$~is \deph{multiplicative} if $A_1\in\A$
and $A_2\in\A$ implies that ${A_1\times A_2\in\A}$.
\end{dfn}

A famous conjecture of Hedetniemi~\cite{Hed:Con} states that the chromatic
number of the product of two (undirected) graphs is equal to the minimum
of the chromatic numbers of the two graphs. This conjecture is equivalent
to the following.

\begin{conj}
The class $\bar\K_k$ of all graphs that are not $k$-colourable
is multiplicative for every positive integer~$k$.
\end{conj}

For some graph classes, such as the class of all graphs that are not homomorphic to
a fixed graph, multiplicativity has an equivalent description. This is in fact
not restricted to graphs, but holds for relational structures as well.

\begin{prop}
Let $H$ be a \Ds. Then the class
\[\HH := \{X : X\notto H\}\]
of all \Ds s that are not homomorphic to~$H$ is multiplicative if and
only if $H$~is irreducible.
\end{prop}

\begin{proof}
Let $H$ be irreducible. By condition~(1) of Lemma~\ref{lem:irreducible},
if the product~$A_1\times A_2$ is homomorphic to~$H$, then $A_1\to
H$ or $A_2\to H$. Thus if $A_1\times A_2\notin\HH$, then $A_1\notin\HH$
or $A_2\notin\HH$.

Conversely, if $\HH$~is multiplicative, then condition~(1) of
Lemma~\ref{lem:irreducible} is satisfied. Hence the \Ds~$H$ is
irreducible.
\end{proof}

Recall from Example~\ref{exa:k-colouring} that a graph is $k$-colourable
if and only if it is homomorphic to the complete graph~$K_k$. Thus
Hedetniemi's conjecture has another equivalent formulation.

\begin{conj}
The complete graph~$K_k$ is irreducible for every positive integer~$k$.
\end{conj}

The conjecture is evidently true for $k=1$. It merely says that if
the product $A_1\times A_2$ has no edges, then $A_1$ or $A_2$ has no
edges either.

It is not difficult to show that if both $A_1$ and $A_2$ are
non-bipartite, then $A_1\times A_2$ is non-bipartite too; this is the
case $k=2$. El-Zahar and Sauer~\cite{EZaSau:Chro} have proved that $K_3$
is irreducible. Almost no other examples of irreducible (undirected)
graphs are known, though. Some more irreducible graphs have been found
by Tardif~\cite{Tar:Mul}. Nevertheless, Hedetniemi's product conjecture
remains wide open for general~$k$.

\section{Paths, trees and forests}

Trees are very important in the context of homomorphism dualities (see
especially Theorem~\ref{thm:dual-pairs}). We define them as structures
whose shadow is a tree.  Forests are then structures consisting of tree
components. In addition, we define paths in \Ds s as ``linear trees'':
every edge has at most two neighbouring edges.

\begin{dfn}
A \Ds~$A$ is called a \deph{\Dt} or simply a \deph{tree} if $\Sh(A)$ is a tree; it is called a
\deph{\Df} or just a \deph{forest} if $\Sh(A)$ is a forest.
\end{dfn}

\begin{pgf}
Notice that $A$ is a \Dt\ if and only if $\Inc(A)$ is a tree; and that
$A$ is a \Df\ if and only if each component of~$A$ is a \Dt.
\end{pgf}

\begin{dfn}
A \Dt~$P$ is called a \deph{\Dp} if every edge of~$P$ intersects at most
two other edges and every vertex of~$P$ belongs to at most two edges.
\end{dfn}

\begin{pgf}
In every \Dp\ with at least two edges there are two edges (\deph{end
edges}) such that each of them shares a vertex with exactly one other
edge. Any other edge (\deph{middle edge}) shares two of its vertices,
each with one other edge.
\end{pgf}

As a generalisation of acyclic directed graphs, that is directed graphs
without directed cycles, we introduce acyclic relational structures.

\begin{dfn}
\label{dfn:acyclic}
A \Ds~$A$ is called \deph{acyclic} if there is no directed cycle in its directed
shadow~$\DSh(A)$.
\end{dfn}

Evidently, every \Dt\ is acyclic. We need this fact especially in
Section~\ref{sec:3constr}.

\section{Height labelling and balanced structures}

Balanced digraphs are digraphs that are homomorphic to a directed
path. The name originates from the fact that the number of forward
edges and the number of backward edges is the same along every cycle in
a balanced digraph (and such cycles are called balanced). For digraphs
being homomorphic to a directed path is the same as being homomorphic
to an oriented forest.  For relational structures it is not, and the
forest definition is the suitable one.

\begin{dfn}
\label{dfn:balanced}
We say that a \Ds~$A$ is \deph{balanced} if $A$~is homomorphic to a \Df.
\end{dfn}

For a balanced digraph, the height of a vertex is defined as the length
of the longest directed path ending in that vertex. As a generalisation,
we define height labelling of relational structures.

\begin{dfn}
\label{dfn:htlab}
Let $A$ be a \Ds\ and let $\gimel$ be a labelling of its vertices
with $(\sum_{i\in I}{\delta_i}-|I|)$-tuples of integers, indexed by
$(i,1),(i,2),\dotsc,(i,\delta_i-1)$, $i\in I$.

We say that $\gimel$~is a \deph{height labelling} of~$A$ if
\begin{itemize}
\item there exists a vertex~$x$ of~$A$ such that $\gimel(x)=(0,0,\dotsc,0)$, and
\item whenever $(x_1,x_2,\dotsc,x_{\delta_i})\in R_i(A)$ and $1\le j<\delta_i$, then
\begin{equation}\label{eq:1.1}
\begin{aligned}
\bigl(\gimel(x_{j+1})\bigr)_{(i,j)} &= \bigl(\gimel(x_{j})\bigr)_{(i,j)}+1\text{,\quad and}\\
\bigl(\gimel(x_{j+1})\bigr)_{(i',j')} &= \bigl(\gimel(x_{j})\bigr)_{(i',j')}\text{\quad
	for $(i',j')\ne (i,j)$.}
\end{aligned}
\end{equation}
\end{itemize}
\end{dfn}

The first condition in the above definition is purely technical; it
facilitates the exposition of arguments. It is possible to omit it
without altering the essence of the definition.

\begin{prop}
\label{prop:htlab}
If $A$~is a balanced \Ds, then $A$~has a height labelling. If a height
labelling of a connected structure exists, it is unique up to an additive
constant vector.
\end{prop}

\begin{proof}
First we prove that a height labelling exists for any \Dt~$T$: pick an
arbitrary vertex~$x$ and set $\gimel(x)=(0,0,\dotsc,0)$. On all other
vertices the labelling is defined recursively. If a neighbour~$y$
of~$x$ in the shadow~$\Sh(T)$ has already been assigned a label, the
label~$\gimel(x)$ is determined by the conditions~\eqref{eq:1.1}.  In this
way, a label is assigned to every vertex, because $\Sh(T)$~is a tree.

By separately labelling each component, we get a height labelling also for any \Df.

Now if $f:A\to T$ is a surjective homomorphism from a balanced \Ds~$A$
to a \Df~$T$, let us fix a height labelling of~$T$ and for a vertex~$x$
of~$A$ define
\[\gimel(x):=\gimel\bigl(f(x)\bigr).\]
The labelling defined in this way is a height labelling; it is not difficult
to check the conditions~\eqref{eq:1.1}.

Uniqueness follows from the fact that the difference between the labels
of two vertices of~$A$ depends only on what edges of~$A$ generated the
path between the two vertices in the shadow~$\Sh(A)$.
\end{proof}

\begin{pgf}
\label{pgf:htlabcnt}
We may picture height labelling as if we replace each edge of the $i$th
kind by ${\delta_i-1}$ binary edges of various kinds and then count forward
edges and backward edges. A height labelling for a structure exists if
and only if each path (in the shadow) between any fixed pair of vertices counts
the same difference of the numbers of forward edges and backward edges
for all kinds. This is related to balanced structures and the height of
a path, as defined for digraphs in~\cite{NesZhu:Path}.
\end{pgf}

\section{Partial orders}
\label{sec:introrder}

Here we state the definitions and elementary facts about
partial orders and lattices that we need later (chiefly in
Chapter~\ref{chap:order}).  Several introductory books on this topic
are available, such as the recommended one by Davey and Priestley
~\cite{DavPri:Introduction}.

\begin{dfn}
A \deph{partial order} is a binary relation~$\preceq$ over a set~$P$
which is reflexive, antisymmetric, and transitive, that is for all $a$,
$b$, and $c$ in~$P$, we have that:
\begin{enumerate}
\item $a\preceq a$ (reflexivity);
\item if $a\preceq b$ and $b\preceq a$ then $a = b$ (antisymmetry); and
\item if $a\preceq b$ and $b\preceq c$ then $a\preceq c$ (transitivity).
\end{enumerate}

A set with a partial order is called a \deph{partially ordered set} or a \deph{poset}.
So formally, a partially ordered set is an ordered pair $(P,{\preceq})$,
where $P$~is called the \deph{base set} and $\preceq$~is a partial order over~$P$.
\end{dfn}

\begin{dfn}
\label{dfn:downset}
A subset $Q\subseteq P$ of the poset~$(P,{\preceq})$ is a \deph{downset}
if, for all elements $p$ and~$q$, if $p$~is less than or equal to~$q$
and $q$~is an element of~$Q$, then $p$~is also in~$Q$. For an arbitrary
subset $S\subseteq P$ the \deph{downset generated by~$S$} is the smallest
downset that contains~$S$; it is denoted by~$\downs S$.
\end{dfn}

\begin{dfn}
\label{dfn:upset}
Similarly, a subset $Q\subseteq P$ of the poset~$(P,{\preceq})$
is an \deph{upset} if, for all elements $p$ and~$q$, if $q\preceq p$
and $q$~is an element of~$Q$, then $p$~is also in~$Q$. For an arbitrary
subset $S\subseteq P$ the \deph{upset generated by~$S$} is the smallest
upset that contains~$S$; it is denoted by~$\ups S$.
\end{dfn}

\begin{dfn}
Suppose $\preceq$ is a partial order on a nonempty set~$P$. Then the
elements $p,q\in P$ are said to be \deph{comparable} provided $p\preceq
q$ or $q\preceq p$.  Otherwise they are called \deph{incomparable}.
A subset~$Q$ of~$P$ is an \deph{antichain} if all elements of~$Q$ are
pairwise incomparable. An antichain is \deph{maximal} if there is no
other antichain strictly containing it.
\end{dfn}

\begin{dfn}
Given a subset~$Q$ of a poset~$P$, the \deph{supremum} of~$Q$ is the
least element of~$P$ that is greater than or equal to each element of~$Q$;
so the supremum of~$Q$ is an element~$u$ in~$P$ such that
\begin{enumerate}
\item $x \preceq u$ for all $x$ in~$Q$, and
\item for any $v$ in~$P$ such that $x \preceq v$ for all $x$ in~$Q$ it holds
that $u \preceq v$.
\end{enumerate}
Dually, the \deph{infimum} of~$Q$ is the greatest element of~$P$ that is
less than or equal to each element of~$Q$.
\end{dfn}

\begin{dfn}
A \deph{lattice} is a poset whose subsets of size two all have a supremum
(called \deph{join}; the join of $p$ and~$q$ is denoted by $p\vee q$)
and an infimum (called \deph{meet}; the meet of $p$ and~$q$ is denoted
by $p\wedge q$).
\end{dfn}

There are many natural examples of lattices.

\begin{exa}\label{exa:latt-natur}
The natural numbers in their usual order form a lattice, under the
operations of minimum and maximum.

The positive integers also form a lattice under the operations of taking
the greatest common divisor and least common multiple, with divisibility
as the order relation: $a \le b$ if $a$~divides~$b$.
\end{exa}

\begin{exa}\label{exa:latt-set}
For any set~$A$, the collection of all subsets of~$A$ can be ordered via
subset inclusion to obtain a lattice bounded by~$A$ itself and the empty
set. Set intersection and union interpret meet and join, respectively.
\end{exa}

\begin{exa}\label{exa:latt-vect}
All subspaces of a vector space~$V$ form a lattice. Here the meet is the
intersection of the subspaces and the join of two subspaces $W$ and $W'$
is the minimal subspace that contains the union~$W\cup W'$.
\end{exa}

An important property of lattices is distributivity. We meet it
again in Lemma~\ref{lem:heyt-distr} and the following paragraphs.

\begin{dfn}
A lattice~$L$ is \deph{distributive} if the following identity holds
for all $x$, $y$, and~$z$ in~$L$:
\[x\wedge (y\vee z) = (x\wedge y) \vee (x\wedge z).\]
\end{dfn}

This says that the meet operation preserves non-empty finite joins. It
is a basic fact of lattice theory that the above condition is equivalent
to its dual:
\[x\vee (y\wedge z) = (x\vee y) \wedge (x\vee z).\]

\begin{exa}
Not every lattice is distributive. From our examples above, the
three lattices in \ref{exa:latt-natur} and \ref{exa:latt-set} are
all distributive. However, the lattice of subspaces of a vector space
from~\ref{exa:latt-vect} is not distributive: Consider the space~$\R^2$
and let $A$, $B$ and~$C$ be three distinct subspaces of dimension~$1$
(straight lines). Then $A\vee B=\R^2$ and $A\cap C=B\cap C=0$, and thus
\[C=(A\vee B) \cap C \ne (A\cap C)\vee (B\cap C) = 0.\]
\end{exa}

Next we introduce  the notion of a gap, which is an island of non-density
in a partially ordered set.

\begin{dfn}
\label{dfn:gap}
A pair $(p,q)$ of elements of a poset~$P$ is a \deph{gap} if $p\prec q$,
and for every $r\in P$, if $p\preceq r\preceq q$ then $r=p$ or $r=q$.
\end{dfn}

In a dense poset there exists $r$ such that $p\prec r\prec q$ for any
pair~$(p,q)$ with $p\prec q$. Thus a poset is dense if and only if it
has no gaps. Gaps are further discussed in Sections~\ref{sec:gaps}
and~\ref{sec:heyt}.

\setchapterpreamble[u]{%
	\dictum[Kingsley Amis, Lucky Jim]{%
	It was one~more argument to~support his~theory that
	nice~things are~nicer than~nasty~ones.}}

\chapter{Homomorphism dualities}

\bigskip\bigskip

Homomorphism dualities are situations where a class of relational
structures is characterised in two ways: by the non-existence of a
homomorphism \emph{from} some fixed set of structures, and by the
existence of a homomorphism \emph{to} some other fixed set of structures.

For example, an (undirected) graph is bipartite if and only if it
is homomorphic to the complete graph~$K_2$; so the existence of a
homomorphism to~$K_2$ is determined by the non-existence of a homomorphism
from odd cycles.

An important aspect of dualities is that in some cases they make the respective class
more accessible.  Duality guarantees the existence of a certificate for
positive as well as negative answers to the membership problem. In both
cases it is a homomorphism, either a homomorphism from the ``forbidden''
set or a homomorphism to the other (``dual'') set of structures.

In special cases this provides an example of a \emph{good
characterisation} in the sense of Edmonds~\cite{Edm:PTF}. This may
mean that an effective algorithm for testing membership is available.

Homomorphism dualities have been studied for many years now.
Their roots go back to the early 1970's, appearing already in
Nešetřil's textbook on graph theory~\cite{N:TG}. The pioneering
work was done by Nešetřil and Pultr~\cite{NesPul:SubFac}. They
found the first instances of homomorphism dualities~-- the
duality pairs of directed paths and transitive tournaments (see
Example~\ref{exa:dipath-dual}). More duality pairs were later
discovered by Komárek~\cite{Kom:Somenew}, and his work led to a
characterisation theorem for digraphs~\cite{Kom:Phd}. Nešetřil and
Tardif~\cite{NesTar:Dual} found a connection to the homomorphism order
(Chapter~\ref{chap:order}) and generalised the notion for relational structures
and together with the author of this thesis they have recently found a
full characterisation~\cite{FNT:GenDu} (Theorem~\ref{thm:finite-character}
here). Meanwhile, dualities have again been generalised in two different
contexts~\cite{NesOss:ResDual,NesPulTar:HeytDual}.

In this chapter, we first investigate dualities with homomorphisms to
a single structure characterised by the non-existence of a homomorphism
from a single structure (duality pairs) and show some properties of such
structures. Then we fully characterise finite dualities. Finally, we address
some extremal problems related to the size of the involved structures.

\section{Duality pairs}
\label{sec:pairs}

Duality pairs are the simplest cases of homomorphism dualities. Here
a class of relational structures is characterised by the existence of
a homomorphism to a single structure~$D$ and at the same time by the
non-existence of a homomorphism from a single structure~$F$.

\begin{dfn}
\label{dfn:dpair}
Let $F$, $D$ be relational structures. We say that the pair $(F,D)$ is a \deph{duality pair}
if for every structure~$A$ we have $F\to A$ if and only if $A\notto D$.
\end{dfn}

One should view the \Ds~$F$ as a characteristic obstacle, which prevents
a \Ds\ from being homomorphic to~$D$.

We introduce the following short notation for a symbolic description of duality pairs.

\begin{pgf}
Let
\begin{align*}
\Fnto &:= \{ A : F\notto A\},\\
\toD  &:= \{ A : A\to D \}.
\end{align*}
Then $(F,D)$ is a duality pair if and only if
\[\Fnto=\toD.\]
\end{pgf}

For proving that a certain pair of structures is indeed a duality pair,
the following characterisation is frequently convenient.

\begin{lemma}\label{lem:equivdual}
Let $F$ and $D$ be \Ds s. Then $(F,D)$ is a duality pair if and only if
$F\notto D$, and whenever $F\notto A$, then $A\to D$.
\end{lemma}

\begin{proof}
If $(F,D)$ is a duality pair, then $F\to F$ and hence $F\notto D$ by
the definition of a duality pair. This proves one implication.

Next is the proof of the opposite implication.  According to the
definition, it has to be shown that if $F\to A$, then $A\notto D$. Suppose
that $F\to A$ and $A\to D$. Then by composition of homomorphisms $F\to
D$. That is a contradiction.
\end{proof}

\begin{pgf}
\label{pgf:fdnot}
In our notation for duality pairs, the letter~$F$ stands for
\emph{forbidden} and the letter~$D$ for \emph{dual}.
\end{pgf}

As a warm-up, we state a few results, which are not difficult to prove
but they provide the first experience with homomorphism dualities.

The first proposition states the obvious: that replacing $F$ or $D$
in a duality pair with a homomorphically equivalent structure does not
change the duality.

\begin{prop}
Let $(F,D)$ be a duality pair. If $F'$ is homomorphically equivalent
to~$F$, and $D'$ is homomorphically equivalent to~$D$, then $(F',D')$
is also a duality pair.
\end{prop}

\begin{proof}
By Lemma~\ref{lem:equivdual} we have $F\notto D$. By the same lemma it
suffices to prove that $F'\notto D'$, and that whenever $F'\notto A$,
then $A\to D'$.

First we prove that $F'\notto D'$. Suppose on the contrary that $F'\to
D'$. Then $F\to D$ because $F\to F'\to D'\to D$ since homomorphisms
compose. That is a contradiction.

Now we show that whenever $F'\notto A$, then $A\to D'$. Suppose that
$F'\notto A$, then also $F\notto A$, so $A\to D$ because $(F,D)$ is
a duality pair. Since $D\to D'$, we have $A\to D'$ by composition of
homomorphisms.
\end{proof}

Next we show that in a duality pair $(F,D)$ the structure~$D$ is
uniquely determined by~$F$ up to homomorphic equivalence. Vice versa,
up to homomorphic equivalence the structure~$F$ is determined by~$D$.

\begin{prop}\label{prop:dp-are-eq}
Let $F$, $F'$, $D$, $D'$ be \Ds s. If $(F,D)$ and $(F,D')$ are duality
pairs, then $D$ and $D'$ are homomorphically equivalent. If $(F,D)$
and $(F',D)$ are duality pairs, then $F$ and $F'$ are homomorphically
equivalent.
\end{prop}

\begin{proof}
Suppose that both $(F,D)$ and $(F, D')$ are duality pairs.
Lemma~\ref{lem:equivdual} implies that $F\notto D$ since $(F,D)$~is a duality
pair, so $D\to D'$ because $(F,D')$~is a duality pair. Moreover $F\notto
D'$ because $(F,D')$~is a duality pair, so $D\to D'$ because $(F,D)$~is
a duality pair. Hence $D\heq D'$.

The proof of the second part is analogous.
\end{proof}

\begin{cor}
If $(F,D)$ is a duality pair, then there exists a unique core~$D'$
such that $(F,D')$~is a duality pair.
\qed
\end{cor}

This corollary motivates the following definition.

\begin{dfn}
Let $D$ be a core and let $(F,D)$ be a duality pair. Then the \Ds~$D$
is called \deph{the dual} of~$F$; it is denoted by~$D(F)$.
\end{dfn}

Now we present the first example of duality pairs.

\begin{exa}[\cite{NesPul:SubFac}]
\label{exa:dipath-dual}
Let us first consider digraphs, that is \Ds s with $\Delta=(2)$. Let $\vec P_k$ denote the
directed path with $k$~edges and $\vec T_k$ the transitive tournament on $k$~vertices. Then for
$k\ge1$ the pair~$(\vec P_k, \vec T_k)$ is a duality pair.
\end{exa}

\begin{proof}
We proceed by induction on~$k$. For $k=1$, clearly $\vec P_1\notto \vec T_1$, and if $\vec
P_1\notto A$, then $A$~has no edges, whence $A\to\vec T_1$. So $(\vec P_1, \vec T_1)$ is a
duality pair by Lemma~\ref{lem:equivdual}.

Now we prove the induction step for $k\ge 2$. Again, obviously $\vec
P_k\notto \vec T_k$. Moreover, suppose that $\vec P_k\notto A$. Then the
digraph~$A$~is acyclic (it contains no directed cycles). Therefore there
exists a vertex of~$A$ with out-degree zero, that is with no outward
edges going from it. Let $A'$~be the digraph created from~$A$ by
deleting all vertices with out-degree zero. It is easy to see that
$\vec P_{k-1}\notto A'$, so by induction we have that $f':A'\to\vec
T_{k-1}$. We can extend the homomorphism~$f'$ to a homomorphism from~$A$
to~$\vec T_k$: the vertices of~$A'$ are mapped in the same way as by~$f'$,
to the subgraph of~$\vec T_k$ consisting of non-zero out-degree vertices,
which is isomorphic to~$\vec T_{k-1}$; and the remaining vertices of~$A$
are mapped to the terminal vertex (sink) of~$\vec T_k$.
\end{proof}

The next example is due to Komárek~\cite{Kom:Somenew}. It is presented
here without proof.  However, let us point out that the example is the
essence of the mosquito construction (see Section~\ref{sec:3constr}) and
historically it was an important milestone on the way to the description
of all duality pairs.

\begin{exa}[\cite{Kom:Somenew}]
For two positive integers $m$, $n$, define the digraph $P_{m,n}=(V,E)$
to be the oriented path with vertices
\[V =\{a_0,a_1,\dotsc,a_m,b_0,b_1,\dotsc,b_n\}\]
and edges
\begin{multline*}
E = \bigl\{(a_j,a_{j+1}) : j=0,1,\dotsc,m-1\bigr\} \cup
	\bigl\{(b_j,b_{j+1}) : j=0,1,\dotsc,n-1\bigr\} \\
	\cup\bigl\{(b_0,a_m)\bigr\}
\end{multline*}
(see~Figure~\ref{fig:P52}).
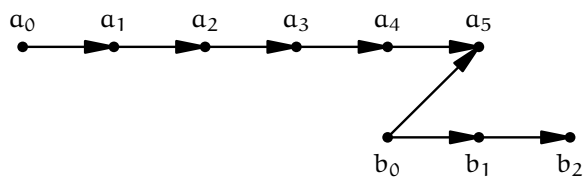
\picture{duality.1}{The path $P_{5,2}$}{P52}

Furthermore, let $D_{m,n}=(W,F)$ be the digraph with the vertex set defined by
\begin{align*}
W &=\bigl\{(i,j):i\ge0,\ j\ge0,\ 0\le i+j\le m+n-2\bigr\}\\
\intertext{and edges by}
F &= \Bigl\{\bigl((i,j),(i',j')\bigr) :
	i < i',\ j > j',\text{ and $i<m$ or $j<n$}\Bigr\}.
\end{align*}
An example is in Figure~\ref{fig:D23}.
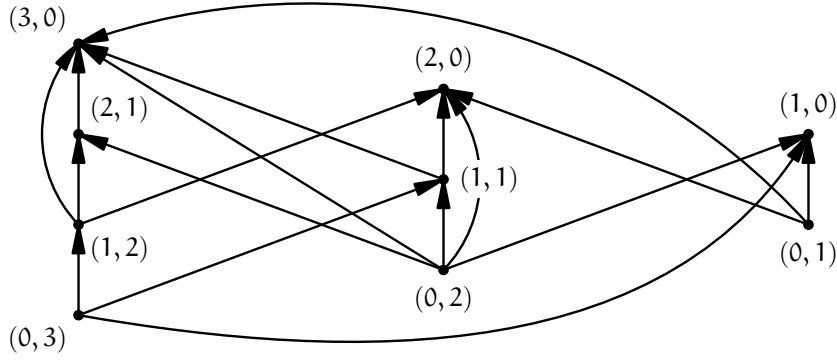
\picture{duality.2}{The digraph $D_{2,3}$}{D23}

Then $(P_{m,n},D_{m,n})$ is a duality pair.
\end{exa}

In the next example, we consider \Ds s with more than one relation. It
shows that if some relations of a \Ds~$F$ are empty, then in the dual
structure the corresponding relations contain all possible tuples of
vertices.

\begin{exa}
\label{exa:empty}
Let $\Delta=(\delta_i:i\in I)$, let $F=\bigl(\bs F, (R_i:i\in I)\bigr)$ be a \Ds\
and let $I'\subseteq I$ be a set of indices such that $R_i=\emptyset$
for all $i\in I\setminus I'$. Define $\Delta'=(\delta_i:i\in I')$ and let $F'$ be the
$\Delta'$-struc\-ture~$\bigl(\bs F,(R_i:i\in I')\bigr)$. Suppose that 
there exists a $\Delta'$-struc\-ture~$D'$ such that $(F',D')$~is a duality pair.
By adding complete relations to~$D'$ we get the \Ds~$D$ with
\begin{align*}
\bs D &= \bs{D'},\\
R_i(D) &=
\begin{cases}
R_i(D') & \text{if $i\in I'$},\\
(\bs D)^{\delta_i} & \text{if $i\in I\setminus I'$}.
\end{cases}
\end{align*}
Then the pair $(F,D)$ is a duality pair.
\end{exa}

\begin{proof}
We use Lemma~\ref{lem:equivdual}.
First, $F\notto D$ because $F'\notto D'$. Now let $A$ be a \Ds\
such that $F\notto A$. This means that $F'\notto A'=\bigl(\bs A,(R_i(A):{i\in I'})\bigr)$ and
therefore there exists a homomorphism $f': A'\to D'$.
Obviously $f: a\mapsto f'(a)$ is a homomorphism from~$A$ to~$D$.
\end{proof}

The following theorem characterises all duality pairs.

\begin{thm}[\cite{NesTar:Dual}]
\label{thm:dual-pairs}
If $(F,D)$ is a duality pair, then $F$~is homomorphically equivalent to
a \Dt. Conversely, if $F$~is a \Dt\ with more than one vertex,  then
there exists a unique (up to homomorphic equivalence) structure~$D$
such that $(F,D)$ is a duality pair.
\end{thm}

\begin{proof}

The proof is split into two parts. Here we prove that if $(F,D)$
is a duality pair and $F$~is a core, then $F$~is a tree.
The second part of the proof can be found in Section~\ref{sec:3constr},
which contains a construction of the dual structure for any tree~$F$.
Uniqueness follows from Proposition~\ref{prop:dp-are-eq}.

Our proof uses an idea of Komárek~\cite{Kom:Phd}, who proved the
characterisation of duality pairs for digraphs. We assume that in
a duality pair~$(F,D)$ the \Ds~$F$ is a core that is not a \Dt. The
idea of the proof is to construct an infinite sequence of \Ds s~$F_k$
such that $F$~is not homomorphic to any of them. By duality, all the
structures~$F_k$ are homomorphic to~$D$, and we will show that this
implies that $F\to D$. That is a contradiction with the definition of
a duality pair (see Lemma~\ref{lem:equivdual}).

Hence if $(F,D)$ is a duality pair and $F$~is a core, then $F$~is
a \Dt. If $(F',D)$ is a duality pair and $F'$~is arbitrary, then by
Proposition~\ref{prop:dp-are-eq} also the pair $(F,D)$ is a duality
pair, where $F$~is the core of~$F'$. So $F$~is a \Dt, and therefore
$F'$~is homomorphically equivalent to a \Dt.

This part of the proof is split into eight steps
(\ref{pgf:primal-connected}--\ref{pgf:dpchfin}).

\begin{pgf}\label{pgf:primal-connected}
\emph{$F$ is connected.} Because: If $F\to A+B$, then $A+B\notto D$, and so either $A$ or~$B$
is not homomorphic to~$D$, and so $F\to A$ or $F\to B$. This is one of the equivalent
descriptions of connectedness in Lemma~\ref{lem:connected}.
\end{pgf}

We proceed to define the structures~$F_k$ for all positive integers~$k$.

\begin{pgf}\label{pgf:dfn-Fk}
We suppose that $F$~is not a \Dt; thus $\Sh(F)$~contains a cycle. Let $\{u,v\}$ be an edge
of~$\Sh(F)$ that lies in a cycle. This edge was added to the shadow~$\Sh(F)$ because there is
an edge~$e\in R_j(F)$ that contains $u$ and~$v$ as consecutive vertices. Now, $F_1$~will be the
\Ds\ constructed from~$F$ by removing the edge~$e$. Furthermore, $F_2$~is constructed by taking
two copies of~$F_1$ and by joining them by two edges of the $j$th kind in the way depicted
in Figure~\ref{fig:F1}.

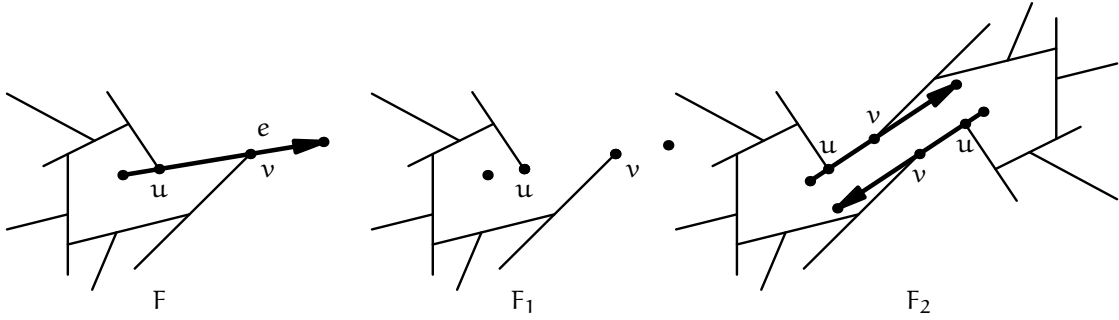
\picture{duality.4}{The construction of $F_1$ and $F_2$}{F1}

In general, for $k\ge 2$, take $k$~copies of~$F_1$ and join each with all
other copies by edges of the~$j$th~kind. Two edges are used to connect
each pair of copies of~$F_1$, so altogether $2\binom{k}{2}$~new edges
are introduced. The resulting structure is~$F_k$.

Formally, if $e=(e_1,\dotsc,u,v,\dotsc,e_{\delta_j})$, define $F_k$
in the following way:
\begin{align*}
\bs{F_k} &:= \{1,2,\dotsc,k\}\times\bs F,\\
R_j(F_k) &:= \Bigl\{\bigl((q,x_1),(q,x_2),\dotsc,(q,x_{\delta_j})\bigr) :
			e \ne (x_1,x_2,\dotsc,x_{\delta_j}) \in R_j(F),\
			1\le q\le k\Bigr\}\\
	&\qquad\cup\Bigl\{\bigl((q,e_1),(q,e_2),\dotsc,(q,u),
			(q',v),\dotsc,(q',e_{\delta_j})\bigr) :
			1\le q,q'\le k,\
			q\ne q'\Bigr\},\\
R_i(F_k) &:= \Bigl\{\bigl((q,x_1),(q,x_2),\dotsc,(q,x_{\delta_i})\bigr) :
			(x_1,x_2,\dotsc,x_{\delta_i}) \in R_i(F),\
			1\le q\le k\Bigr\}\\
	&\qquad\qquad\text{for $i\ne j$}.
\end{align*}
\end{pgf}

Next we prove several properties of the structures~$F_k$.

\begin{pgf}
\label{pgf:homh}
\emph{$F_k\to F$.} Because: The identity mapping is indeed a homomorphism
from~$F_1$ to~$F$. For $k\ge2$, mapping each vertex~$(q,x)$ to the corresponding vertex~$x$
of~$F$ provides a homomorphism (let us call this homomorphism~$h$; so $h:(q,x)\mapsto x$).
\end{pgf}

\begin{pgf}\label{pgf:F-not-to-F1}

\emph{$F\notto F_1$.} Because: The structures $F$ and $F_1$ have the same
number of vertices but $F_1$~has fewer edges. So a potential homomorphism
from~$F$ to~$F_1$ cannot be injective and so it would actually map~$F$
to a proper substructure of~$F$. This is impossible because $F$~is a core.
\end{pgf}

In the following four paragraphs, let $n:=|\bs F|$ and let the vertices
be enumerated in such a way that $\bs F=\{x_1,x_2,\dotsc,x_n\}$ with
$u=x_1$ and $v=x_2$.

\begin{pgf}\label{pgf:one-copy}
\emph{If $f:F\to F_k$ is a homomorphism, then the image~$f[F]$~consists of exactly one copy of
each vertex of~$F_1$ in~$F_k$. In particular, any homomorphism $f:F\to F_k$ is injective.}
Because: 
If $f:F\to F_k$ is a homomorphism, then the image~$f[F]$~contains vertices
of at least two distinct copies of~$F_1$ in~$F_k$ because of~\ref{pgf:F-not-to-F1}.
If some vertex appears in more than one copy, say $f(x_i)=(p,x_l)$
and $f(x_{i'})=(q,x_l)$ for some $p\ne q$, then the composed homomorphism $hf:F\to F$
is not an automorphism. This is a contradiction, since $F$~is a core.
Here $h$~is the homomorphism $(q,x)\mapsto x$ defined in~\ref{pgf:homh}.
\end{pgf}

\begin{pgf}\label{pgf:218}
\emph{If $F\to F_k$ for some~$k$, then $F\to F_2.$} 
Because:
Let $f:F\to F_k$. By~\ref{pgf:one-copy}, the image~$f[F]$~contains exactly one vertex $(p,x_1)$
and exactly one vertex~$(q,x_2)$. As $F$~is connected, the image contains vertices of
only at most two copies of~$F_1$ in~$F_k$, namely vertices in the form $(p,x_i)$ and $(q,x_i)$.
So the image~$f[F]$ lies entirely in a substructure of~$F_k$ that is isomorphic to~$F_2$.
Thus $F\to F_2$.
\end{pgf}

\begin{pgf}\label{pgf:fntfk}
\emph{$F\notto F_k$.}
Because:
Suppose $F\to F_k$. Then there is a homomorphism $f:F\to F_2$ because
of~\ref{pgf:218}. We may assume that $(1,x_1)$ and $(2,x_2)$ are images
under~$f$ of some vertices of~$F$.

Now, $f[\bs F]=V_1\cup V_2$, where $V_1$~contains vertices in the
first copy and $V_2$~vertices in the second copy of~$F_1$ in~$F_2$,
that is $V_l:=\bigl\{x:x\in f[\bs F]\text{ and } x=(l,x_i)\text{ for
some $i$}\bigr\}$ for $l=1,2$. Then $(1,x_1)\in V_1$ and $(2,x_2)\in V_2$.

For a \Ds~$A$, let $r(A)$ denote the number of all edges of~$A$, so
\[r(A):=\Bigl|\bigcup_{i\in I}R_i(A)\Bigr|.\]
Let $G$~be the substructure
induced in~$F_2$ by~${V_1\cup V_2}$. The homomorphism~$f$ is injective
by~\ref{pgf:one-copy}, and as a result it maps every edge of~$F$ to an
edge of~$G$. Thus $r(G)\ge r(F)$.

The edges of~$G$ are split into three groups: the edges
induced by~$V_1$ (let their number be~$m_1$), the edges
induced by~$V_2$ (let their number be~$m_2$) and the edge
\[\bigl((1,e_1),(1,e_2),\dotsc,(1,u),(2,v),\dotsc,(2,e_{\delta_j})\bigr).\]
Hence $r(G)=m_1+m_2+1$.

Let the base set $\bs F$ be partitioned into $W_1$ and $W_2$ by setting
$W_l:=\bigl\{x : (l,x)\in V_l\bigr\}$. So $W_1$~consists of those
vertices~$x$ of~$F$ whose copy~$(1,x)$ is in the image~$f[\bs{F}]$;
similarly for~$W_2$.

Because of~\ref{pgf:one-copy} we indeed
have that $\bs F=W_1\cup W_2$ and ${W_1\cap W_2=\emptyset}$. Observe
that the set~$W_1$ induces $m_1$~edges in~$F$ and the set~$W_2$
induces $m_2$ edges.  Recall that to get~$F_1$ we deleted the
edge~$e=(e_1,\dotsc,u,v,\dotsc,e_{\delta_j})$ of~$F$. The edge~$e$
contains vertices $e_1,\dotsc,u\in W_1$ and $v,\dotsc,e_{\delta_j}\in
W_2$, so it is not induced by either $W_1$ or~$W_2$. Moreover, $u$ and $v$
are vertices of a cycle in~$\Sh(F)$, so there has to be another edge in
the cycle with one end in~$W_1$ and the other end in~$W_2$. This edge
appears in yet another edge of~$F$. Therefore $r(F)\ge m_1+m_2+2$.

We conclude that $m_1+m_2+1=r(G)\ge r(F)\ge m_1+m_2+2$, a contradiction.
\end{pgf}

Now we can derive a contradiction, thus disproving the assumption that $F$~is \emph{not} a \Dt.

\begin{pgf}\label{pgf:dpchfin}
As a consequence of duality, $F_k\to D$ for all~$k$, because $F\notto
F_k$ by~\ref{pgf:fntfk}.  If
$e=(e_1,\dotsc,{e_t=u},{v=e_{t+1}},\dotsc,e_{\delta_j})$
is the deleted edge, let $k>|\bs D|^t$ and let $f:F_k\to D$
be a homomorphism. Then for some $p$,~$q$ such that $p\ne q$ we have
$f(p,e_l)=f(q,e_l)$ for all~$l$ satisfying that $1\le l\le t$.

Define $g:\bs F\to\bs D$ by setting $g(x):=f(p,x)$. If
$(x_1,x_2,\dotsc,x_{\delta_i})\in R_i(F)$ is an edge distinct from~$e$, then
$\bigl((p,x_1),(p,x_2),\dotsc,(p,x_{\delta_i})\bigr)\in R_i(F_k)$ and
therefore \[\bigl(g(x_1),g(x_2),\dotsc,g(x_k)\bigr)\in R_i(D).\] Besides,
\[\bigl(g(e_1),g(e_2),\dotsc,g(e_t)\bigr)\in R_j(D)\] because $f$~is a homomorphism,
\begin{multline*}
\bigl(g(e_1),g(e_2),\dotsc,g(e_t)\bigr)=\bigl(f(p,e_1),f(p,e_2),\dotsc,f(p,e_t)\bigr)\\
=\bigl(f(q,e_1),\dotsc,f(q,e_t),f(p,e_{t+1}),\dotsc,f(p,e_{\delta_j})\bigr)
\end{multline*}
and
\[\bigl((q,e_1),\dotsc,(q,e_t),(p,e_{t+1}),\dotsc,(p,e_{\delta_j})\bigr)\in
R_j(F_k).\] Hence $g$~is a homomorphism from~$F$ to~$D$, a contradiction with duality.
\end{pgf}
This proves that $F$~is a \Dt.
\end{proof}

\section{Three constructions}
\label{sec:3constr}

In this section, we first show two ways to construct the dual of a \Dt~$F$. The
first construction is by Komárek~\cite{Kom:Phd} for digraphs (we will call
it the \emph{mosquito construction}), the second by Nešetřil and
Tardif~\cite{NesTar:A-dualistic,NesTar:Short} for general \Ds s (called
here the \emph{bear construction}).

This section's main purpose is to present a more general construction, of
which both cited constructions are special cases.

\subsection*{Mosquito construction}

For the mosquito construction, we consider digraphs, so let $\Delta=(2)$. Let $F$ be a core tree.
The duals of directed paths are described in Example~\ref{exa:dipath-dual}, so here we suppose
that $F$~is a core digraph other than a directed path.

First we introduce some notation.

\begin{pgf}\label{pgf:mu-komar}
We define the function~$\mu:\bs F\to \N^2$ by $\mu(x):=(d,u)$, where
$d$~is the length of the longest directed path ending in the vertex~$x$
and $u$~is the length of the longest directed path starting from~$x$.

Further, let $r(F)$~be the length of the longest directed path
in~$F$, that is $r(f):=\max\{r:\vec P_r\to F\}$.

Note that $r(F)=
\max\bigl\{d+u: (d,u)=\mu(x),\ x\in\bs F\bigr\}$.
\end{pgf}

\begin{pgf}
Let $p(F)$ be the \deph{height} of~$F$ defined as the length of the shortest directed
path that $F$~is homomorphic to:
\[p(F):=\min\{r:F\to\vec P_r\}.\]

For any $p,q\ge0$ define
\[\Phi(p,q):=\bigl\{a\in\bs F: \mu(a)=(d,u),\ d\le p,\ u\le q\bigr\}.\]
\end{pgf}

Now follows the definition of the mosquito dual of the oriented tree~$F$.

\begin{dfn}
Let $F$ be an oriented tree.
The vertices of the \deph{mosquito dual}~$D_m(F)$ are triples $(p,q,\phi)$, where $0\le
p,q < p(F)$ and $p+q < p(F)$ and $\phi:\Phi(p,q)\to \bs F$ is a function such
that $\phi(a)$~is a neighbour (in-neighbour or out-neighbour) of~$a$ for all
$a\in \Phi(p,q)$.

The edges of the dual~$D_m(F)$ are pairs $\bigl((p,q,\phi),(p',q',\phi')\bigr)$ such that
\begin{itemize}
\item[(i)] $p<p'$ and $q>q'$, and
\item[(ii)] there is no edge $(a,b)$ of~$F$ such that $\phi(a)=b$ and $\phi'(b)=a$.
\end{itemize}
\end{dfn}

The correctness of the construction is asserted by the following theorem.

\begin{thm}[\cite{Kom:Phd}]
For an oriented tree~$F$ that is not a directed path, the pair $\bigl(F,D_m(F)\bigr)$ is a
duality pair.
\qed
\end{thm}

\subsection*{Bear construction}

Next we present a construction of the dual by Nešetřil and Tardif, which
works for \Ds s of an arbitrary type~$\Delta$.

\begin{dfn}[\cite{NesTar:Short}]
\label{pgf:bear-def}
Let $F$ be a core \Dt.  Remember the definition of $\Block(F)$ and $\Inc(F)$
from~\ref{dfn:incgraph}.  Define the \deph{bear dual}~$D_b(F)$ as the \Ds\ on the base set
\[
\bs{D_b(F)} := \Bigl\{f:\bs F\to\Block(F) :
	\bigl\{x,f(x)\bigr\}\in E\bigl(\Inc(F)\bigr) \text{ for all }x\in\bs F\Bigr\}\]
with relations
\begin{multline*}
R_i\bigl(D_b(F)\bigr) := \bigl\{ (f_1,f_2,\dotsc,f_{\delta_i}) :
	\text{for all $e=(x_1,x_2,\dotsc,x_{\delta_i})\in R_i(F)$}\\
	\text{there exists $j\in\{1,\dotsc,\delta_i\}$ such that $f_j(x_j)\ne(i,e)$} \bigr\}.
\end{multline*}
\end{dfn}

\begin{thm}[\cite{NesTar:Short}]\label{thm:bear-dual}
For any \Dt~$F$, the pair $\bigl(F,D_b(F)\bigr)$ is a duality pair.
\end{thm}

\subsection*{A generalisation: The animal construction}

In this subsection we provide a common framework for both
constructions. We present a construction with a parameter called a
\emph{positional-function family}. Depending on what family we take, we
get distinct dual structures. Later we will show what positional-function
families have to be considered to give the mosquito construction and
the bear construction.

We start with the definition of a positional-function family.

Recall the definition of an acyclic \Ds\ from~\ref{dfn:acyclic}.

\begin{dfn}\label{dfn:posifam}
Let $(Q,{\preceq})$ be a partially ordered set.
A \deph{positional-func\-tion family} is a family \[\{\mu_A:\text{$A$
is an acyclic structure}\}\] of functions  such that
\begin{itemize}
\item[(i)] $\mu_A:\bs A\to Q$ for all acyclic structures~$A$,
\item[(ii)] whenever $A$ and $B$ are acyclic structures and there exists
a homomorphism $f:A\to B$ such that $f(x)=y$, then $\mu_A(x)\preceq \mu_B(y)$,
\item[(iii)] for any non-empty finite downset~$S$ in~$Q$  (see~\ref{dfn:downset})
there exists an acyclic \deph{representing structure}~$\Theta=\Theta(S)$ with
\begin{enumerate}
\item a mapping $\Omega:\bs\Theta\to S$ such that for any homomorphism~$f$
from an acyclic structure~$A$ to~$\Theta$ and for any vertex $a\in\bs A$,
we have $\mu_A(a)\preceq\Omega\bigl(f(a)\bigr)$, and
\item a mapping $\Psi:S\to\bs\Theta$ such that for every acyclic
structure~$A$ with $\mu_A[A]\subseteq S$ the mapping given by
$a\mapsto\Psi\bigl(\mu_A(a)\bigr)$ is a homomorphism from~$A$ to~$\Theta$.
\end{enumerate}
\end{itemize}
\end{dfn}

In order to construct the dual of a tree~$F$, we need a
positional-function family satisfying a certain finiteness condition,
as we will see shortly. This condition is trivially satisfied if the poset~$Q$ is finite, as is
the case of the bear construction. The condition is as follows.

\begin{dfn}
Let $(Q,{\preceq})$ be a partially ordered set.
Let us have a po\-si\-tional-function family $M=\{\mu_A:\text{$A$ acyclic}\}$
with $\mu_A:\bs A\to Q$ and with representing structures~$\Theta(S)$ for all
non-empty finite downsets~$S$.

Let $F$ be a \Dt. Let \[T:=\bigcup_{\substack{F\notto A\\\text{$A$ acyclic}}} \mu_A[A]\]
and let $S(F):=\downs T$. So $S(F)$~is the smallest downset that contains~$T$ as a subset.
We say that the positional-function family~$M$ is \deph{suitable} for the \Dt~$F$ if the
following condition is satisfied:
\begin{equation}
\label{eq:sfinite}
\text{The downset $S(F)$ is finite.}
\end{equation}
\end{dfn}

\begin{pgf}
\label{pgf:obviousphi}
We define the mapping $\Phi:\bs\Theta\to2^{\bs F}$ by
\[\Phi(\theta):=\bigl\{y\in\bs F: \mu_F(y)\preceq\Omega(\theta)\bigr\}.\]
It is obvious from the definition that if $y\in\bs F$ is mapped
to~$\theta$ by some homomorphism from~$F$ to~$\Theta$, then $y\in\Phi(\theta)$.
\end{pgf}

Now we proceed to define the animal dual of a \Dt~$F$ with respect to
a positional-function family~$M=\{\mu_A:\text{$A$ acyclic}\}$ that is
suitable for~$F$.

\begin{dfn}\label{dfn:animal}
The base set~$\bs{D_a}$ of the \deph{animal dual}~$D_a$ consists
of pairs $(\theta,\phi)$ such that $\theta$~is a vertex of~$\Theta$
and $\phi:\Phi(\theta)\to\Block(F)$ is a mapping such that $y$ and
$\phi(y)$ are adjacent in~$\Inc(F)$, that is $y$~appears in the edge
given by~$\phi(y)$ in~$F$.

The $\delta_i$-tuple
$\bigl((\theta_1,\phi_1),(\theta_2,\phi_2),\dotsc,(\theta_{\delta_i},\phi_{\delta_i})\bigr)$
is an element of the relation~$R_i(D_a)$ if and only if
$(\theta_1,\theta_2,\dotsc,\theta_{\delta_i})\in R_i(\Theta)$
and there is no edge $e=(y_1,y_2,\dotsc,y_{\delta_i})\in R_i(F)$ such that
$\phi_1(y_1)=\phi_2(y_2)=\dotsb=\phi_{\delta_i}(y_{\delta_i})=(i,e)$.
(Some $\phi_j(y_j)$ may be undefined but that does not matter; in such a case
the edge~$e$ does not satisfy the equality.)
\end{dfn}

To prove that the animal construction is correct and the
structure~$D_a$ forms indeed a duality pair with~$F$, we will use
Lemma~\ref{lem:equivdual}. First we prove that $F$~is not homomorphic
to~$D_a$; this is shown by the following two lemmas. The other part of
the proof follows the statement of Theorem~\ref{thm:animal}.

\begin{lemma}
\label{lem:animal1}
Let $A$ be an arbitrary \Ds\ and let $f:A\to D_a$ be a homomorphism.
Let the mapping $g:\bs A\to\bs\Theta$ be defined by $g(a):=\theta$ such that
$f(a)=(\theta,\phi)$. Then $g$~is a homomorphism from~$A$ to~$\Theta$.
\end{lemma}

\begin{proof}
If $(a_1,a_2,\dotsc,a_{\delta_i})\in R_i(A)$, then
$\bigl(f(a_1),f(a_2),\dotsc,f(a_{\delta_i})\bigr)\in R_i(D_a)$ because $f$~is a homomorphism.
Hence $\bigl(g(a_1),g(a_2),\dotsc,g(a_{\delta_i})\bigr)\in R_i(\Theta)$ by the definitions of~$g$
and the dual structure.
\end{proof}

\begin{lemma}
\label{lem:animal2}
The tree $F$ is not homomorphic to~$D_a$.
\end{lemma}

\begin{proof}
Reductio ad absurdum: Suppose that there is a homomorphism $f:F\to
D_a$. Let $y$ be an arbitrary
element of~$\bs F$ and let $f(y)=(\theta,\phi)$. By Lemma~\ref{lem:animal1}, there is a
homomorphism~$g:F\to\Theta$ that maps $y$ to~$\theta$, whence $y\in\Phi(\theta)$
by~\ref{pgf:obviousphi}. So $\phi(y)$~is defined; let
\[\phi(y)=:\bigl(i,(y_1,y_2,\dotsc,y_{\delta_i})\bigr)\] and let
\[e:=(y_1,y_2,\dotsc,y_{\delta_i})\in R_i(F).\]

Since $f$~is a homomorphism, $\bigl(f(y_1),f(y_2),\dotsc,f(y_{\delta_i})\bigr)\in R_i(D_a)$.
By the definition of edges of~$D_a$, there is an index~$j$ such that $\phi_j(y_j)\ne (i,e)$ if
$f(y_j)=(\theta_j,\phi_j)$. Let $\phi_j(y_j)=(i',e')$. So we have a walk of length four
in~$\Inc(F)$, namely $y,e,y_j,e'$, with the property that $y\ne y_j$ and $e\ne e'$.

Repeating the procedure we can get an arbitrarily long walk
\[z_1,e_1,z_2,e_2,\dotsc,z_{n},e_n\] in~$\Inc(F)$ such that $z_j\ne z_{j+1}$
and $e_j\ne e_{j+1}$. This is a contradiction, because $\Inc(F)$~is a tree.
\end{proof}

Now we prove the correctness of the animal construction.

\begin{thm}
\label{thm:animal}
Let $(Q,\preceq)$ be a partially ordered set.  Let $\{\mu_A:\text{$A$
acyclic}\}$ be a po\-si\-tional-function family with representing
structure~$\Theta(S)$ for every non-empty finite downset~$S$ in~$Q$. If
$F$~is a \Dt\ such that the condition~\eqref{eq:sfinite} holds, and
$D_a$~is the \Ds\ defined in~\ref{dfn:animal}, then the pair $(F,D_a)$
is a duality pair.
\end{thm}

\begin{proof}
We have just seen that $F\notto D_a$, and so by Lemma~\ref{lem:equivdual}
it remains to show that whenever $F\notto X$, then $X\to D_a$. The proof
uses an idea of the proof for the bear dual~\cite{NesTar:Short}. This
idea is reworked so as to fit the animal construction.

Fix a labelling~$\ell$ of the vertices of~$\Inc(F)$ by positive integers
such that different vertices get different labels and the subgraph
of~$\Inc(F)$ induced by $\bigl\{u:\ell(u)\ge n\bigr\}$ is a connected
subtree for all positive integers~$n$. Such a labelling can be defined
by repeatedly labelling and deleting the leaves of~$\Inc(F)$.

For a vertex~$y\in\bs F$ and for its neighbour $b=(i,e)\in\Block(F)$
in~$\Inc(F)$, let $T_{y,b}$ be the maximal subtree of~$\Inc(F)$ that
contains $y$ and~$b$ but no other neighbour of~$y$. Let $F_{y,b}$ be the
\Dt\ such that $\Inc(F_{y,b})=T_{y,b}$ (so $F_{a,b}$~is a substructure
of~$F$). For a vertex~$y$ and for $b\ne b'$, the subtrees $F_{y,b}$ and
$F_{y,b'}$ intersect in exactly one vertex: the vertex~$y$.

Let $X$ be a \Ds\ such that $F\notto X$; we will define a mapping $f:\bs
X\to\bs{D_a}$ and prove that it is a homomorphism.

For every $x\in\bs X$ and $y\in\bs F$ we define
\begin{multline*}
K(x,y):=\bigl\{b\in\Block(F):
	\{b,y\}\in E(\Inc(F)),\\
	\text{and there is no homomorphism $g:F_{y,b}\to X$ s.t. $g(y)=x$}\bigr\}.
\end{multline*}
For every $x\in\bs X$ and $y\in\bs F$, the set $K(x,y)$ is
non-empty. Otherwise there would be a homomorphism $g_b:F_{y,b}\to X$
for all edges~$b$ incident with~$y$, satisfying $g_b(y)=x$, and so their
union would define a homomorphism from~$F$ to~$X$.

Let $x\in\bs X$ be an arbitrary vertex. Define $f(x):=(\theta_x,\phi_x)$ by setting
\[\theta_x:=\Psi\bigl(\mu_X(x)\bigr)\]
and let $\phi_x:\Phi(\theta_x)\to\Block(F)$ be defined by letting $\phi_x(y)$ be
the element~$b$ of~$K(x,y)$ with the smallest label~$\ell(b)$.
The element~$\theta_x$ is well-defined because $\mu_X[X]\subseteq S$ by the definition of~$S$.

Suppose that $(x_1,x_2,\dotsc,x_{\delta_i})\in R_i(X)$. We want to show that
\[\bigl(f(x_1),f(x_2),\dotsc,f(x_{\delta_i})\bigr)\in R_i(D_a).\]
Let $f(x_j)=(\theta_j,\phi_j)$.
Then \[(\theta_1,\theta_2,\dotsc,\theta_{\delta_i})=
\Bigl(\Psi\bigl(\mu_X(x_1)\bigr),\Psi\bigl(\mu_X(x_2)\bigr),\dotsc,
	\Psi\bigl(\mu_X(x_{\delta_i})\bigr)\Bigr)\in R_i(\Theta)\]
by the definition of a representing structure.

Thus it remains to prove that there is no edge $e=(y_1,y_2,\dotsc,y_{\delta_i})\in R_i(F)$ such
that $\phi_1(y_1)=\phi_2(y_2)=\dotsb=\phi_{\delta_i}(y_{\delta_i})=(i,e)$.
For the sake of contradiction, suppose that such an edge exists. Let
$N(y_j)$ denote the set of all neighbours of~$y_j$ in~$\Inc(F)$ different
from~$(i,e)$ and let $N:=\bigcup_{1\le j\le\delta_i}N(y_j)$. Among the
elements of~$N$, there
may be at most one with its $\ell$-label bigger than the label of~$(i,e)$ because of the way
$\ell$~was defined.

If there is no such element, then for any~$j$, no element of~$N(y_j)$ belongs to~$K(x_j,y_j)$,
because otherwise we would not have selected~$(i,e)$ as the value of~$\phi_j(y_j)$. Therefore
for all~$j$ and all~$b\in N(y_j)$ there is a homomorphism $g_{j,b}:F_{y_j,b}\to X$ such that
$g_{j,b}(y_j)=x_j$ and the union of all these homomorphisms defines a homomorphism from~$F$
to~$X$, a contradiction.

Thus there is a unique element $b'\in N$ such that
$\ell(b')>\ell(i,e)$. Hence $b'\in N(y_{j'})$ and for all $b\in N$ different
from~$b'$ we can find a homomorphism~$g_{j,b}$ as above.

Then the mapping~$g$ such that
\[g(y):=
\begin{cases}
	x_{j'}&		\text{if $y=y_{j'}$,}\\
	g_{j,b}(y)&	\text{if $y\in\bs{F_{y_j,b}}$ and $y\ne y_{j'}$,}
\end{cases}
\]
is a homomorphism from~$F_{y_j,\,(i,e)}$ to~$X$, proving that $(i,e)\notin K(x_{j'},y_{j'})$ and
so contradicting the value of~$\phi_{j'}$.
\end{proof}

\subsection*{Bear and mosquito are animals}

In the beginning of this section we promised a generalisation of both the bear and the mosquito
constructions. We have seen a \emph{metaconstruction}: it produces
different results depending on what positional-function family we plug
in. Here we show what to plug in to get the two previous constructions.
That proves that indeed the bear construction and the mosquito construction are
special cases of the animal construction.

\begin{exa}[bear]
Let $Q=\{\diamondsuit\}$ be a one-element poset. For an acyclic
structure~$A$ define $\mu_A$ to be the constant mapping that maps
all vertices of~$A$ to~$\diamondsuit$. The conditions (i) and~(ii) of
Definition~\ref{dfn:posifam} are trivially satisfied. Let $\Theta$ be the
\Ds\ defined by $\bs\Theta=Q$ and $R_i(\Theta)=\bs\Theta^{\delta_i}$
for all $i\in I$, and let $\Omega$ and $\Psi$ be the identity
mapping on~$Q=\bs\Theta$. The structure~$\Theta$ is the representing
structure~$\Theta(Q)$ for the only downset~$Q$, since the constant
mapping from any \Ds\ to~$\Theta$ is a homomorphism. Also the
condition~\eqref{eq:sfinite} is satisfied trivially.

It is easy to see that for any \Dt~$F$ the dual structures $D_a(F)$ and $D_b(F)$ are
isomorphic; the mapping defined by $(\theta,\phi)\mapsto \phi$ is an isomorphism.
\end{exa}

\begin{exa}[mosquito]

Now, let $(Q,\preceq)$ be the product ${(\N,\le)}\times{(\N,\le)}$, that
is $Q=\N\times\N$ and $(d,u)\preceq(d',u')$ if and only if $d\le d'$ and
$u\le u'$. For an acyclic structure~$A$, the positional function~$\mu_A$
is defined as in~\ref{pgf:mu-komar}: $\mu_A(a)=(d,u)$, where $d$~is
the length of the longest directed path ending in the vertex~$a$ and
$u$~is the length of the longest directed path starting in~$a$. As a
homomorphism to an acyclic structure maps a directed path bijectively,
the condition~(ii) of Definition~\ref{dfn:posifam} is satisfied.

For a downset~$S$, the representing digraph is $\Theta=(S,E)$, where
\[\bigl((p,q),(p',q')\bigr)\in E\text{ if and only if $p<p'$ and $q>q'$.}\]
Furthermore, we define $\Omega=\Psi=\id_S$. It can be checked easily that this correctly
defines a representing structure.

Because (in the context of digraphs) every tree is homomorphic to a directed path, for any
tree~$F$ the set~$S(F)$ contains pairs $(p,q)$ with $p+q$ bounded from above by the height
of~$F$. Therefore the condition~\eqref{eq:sfinite} is satisfied for any oriented tree~$F$.

Evidently, $D_a(F)$ and $D_m(F)$ are isomorphic.
\end{exa}

{\tolerance1000
\begin{problem}
Find other suitable positional-function families and representing structures to get essentially
new constructions of the dual structure for a general type~$\Delta$.
\end{problem}}

\section{Properties of the dual}

Two particular properties of dual structures are worthwhile to mention:
connectedness and irreducibility.

\subsection*{Dual is irreducible}

Recall that a \Ds~$D$ is called irreducible if $A\times B\to D$
implies that $A\to D$ or $B\to D$ for every two structures $A$, $B$ (see
Lemma~\ref{lem:irreducible}).

We show that the dual of a \Dt\ is irreducible. This statement is dual to the connectedness of
the left-hand side of a duality pair (see~\ref{pgf:primal-connected}).

\begin{prop}\label{prop:dual-irred}
If $(F,D)$ is a duality pair, then the \Ds~$D$ is irreducible.
\end{prop}

\begin{proof}
Let $A$, $B$ be structures such that $A\times B\to D$. By duality, $F\notto A\times B$.
Therefore $F\notto A$ or $F\notto B$; and using duality once more, it follows that $A\to D$ or
$B\to D$.
\end{proof}

\subsection*{Dual is connected}

We have seen that in a duality pair $(F,D)$ the core of~$F$ is connected
(\ref{pgf:primal-connected}). Now we prove that the dual~$D$
is connected too.

\begin{prop}
If $(F,D)$ is a duality pair and $D$~is a core, then the \Ds~$D$ is connected.
\end{prop}

\begin{proof}
Suppose that $F$ and $D$ are core \Ds s, the pair $(F,D)$ is a duality
pair and $D$~is not connected. Hence $F$~has edges of all kinds; otherwise the dual is
connected, see~\ref{exa:empty}.

Let $\Delta=(\delta_i:i\in I)$. Let $K$ be
the \Ds\ that consists of isolated edges, one edge of each kind; formally
\begin{align*}
\bs K &= \bigl\{(i,u) : i\in I \text{ and } 1\le u \le \delta_i\bigr\},\\
R_i(K) &= \Bigl\{\bigl((i,1), (i,2), \dotsc, (i,\delta_i)\bigr)\Bigr\}.
\end{align*}
Now, the structure~$J$ is obtained from~$K$ by gluing all edges at the first vertex, and
$J'$~is obtained by gluing them at the last vertex. In other words, we have two equivalence
relations on~$\bs K$, namely $\approx$ and ${\approx}'$, defined by
\begin{align*}
(i,u) \approx (j,v) &\text{ if } u=v=1 \text{ or } (i,u)=(j,v),\\ 
(i,u) \mathrel{{\approx}'} (j,v) &
	\text{ if } u=\delta_i \text{ and } v=\delta_j \text{ or } (i,u)=(j,v).
\end{align*}
Define $J$ and $J'$ as factor structures (see~\ref{dfn:factor}):
$J:=K/{\approx}$ and $J':=K/{\approx}'$. The construction is illustrated
in Figure~\ref{fig:jjdash}.

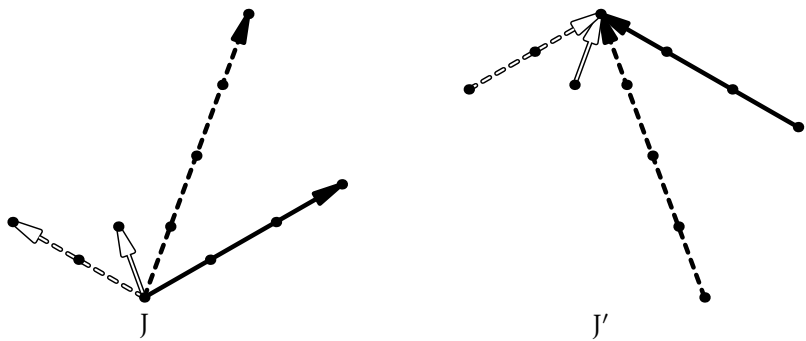
\picture{duality.5}{The construction of $J$ and $J'$}{jjdash}

Suppose $F$~is homomorphic to both $J$ and~$J'$. The structure~$F$ is
connected, so it has at least two incident edges~$e,e'$ of distinct
kinds. As $F\to J$, the vertex these two edges share is the starting
vertex of each of them, because all other vertices are distinct as they
are homomorphically mapped to distinct vertices of~$J$. But since $F\to
J'$, the starting vertices of $e$ and $e'$ are distinct; a contradiction.

Therefore either $J$ or $J'$ is homomorphic to~$D$ if at least two
relations have arity greater than two.

In any case, $D$~has a component that contains edges of all kinds.
This is satisfied trivially if there is only one kind of edges, that is
if $|I|=1$.

So we can connect all components of~$D$ with long zigzags, see
Figure~\ref{fig:zigzags-dual}. We get a connected \Ds~$D'$; if the zigzags
are long enough, any substructure~$E$ of~$D'$ induced by at most~$|\bs
F|$~vertices contains vertices of only one of the components of~$D$;
hence $E$~is homomorphic to~$D$. Because $F$~is \emph{not} homomorphic
to~$D$, it is not homomorphic to~$D'$ either.

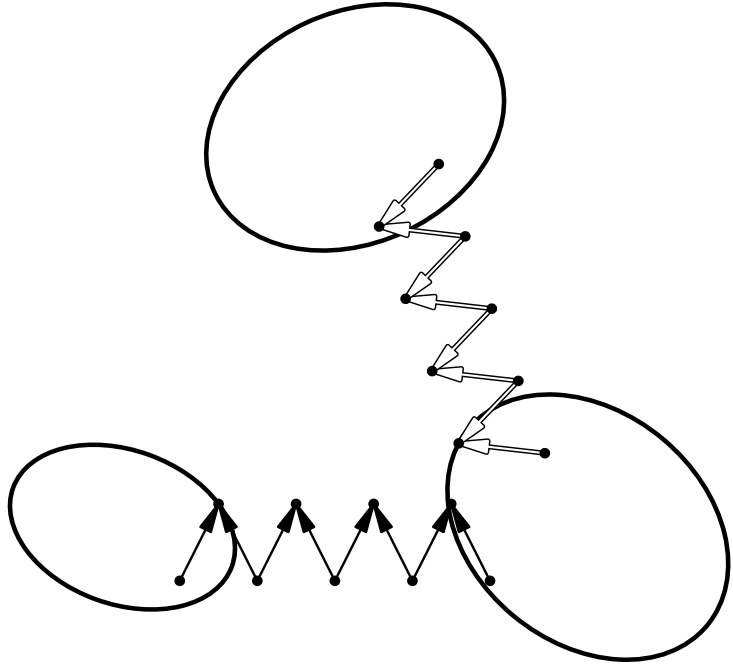
\picture{duality.6}{Connecting the components of a dual with zigzags}{zigzags-dual}

But then, by duality, the
connected structure~$D'$ is homomorphic to~$D$, so it is homomorphic to
a component of~$D$. Therefore $D$ is homomorphic to a proper substructure
of~$D$, a contradiction with $D$ being a core.
\end{proof}

Connectedness of duals was originally proved for digraphs by Nešetřil
and Švej\-da\-rová~\cite{NesIda:Diam} in a different way, by examining the
bear construction.

Dually, we would expect all trees to be irreducible. But this shows
another difference between the dual notions of connectedness and
irreducibility; not all trees, not even all paths are irreducible,
as the following example demonstrates.

\begin{exa}
\label{exa:tree-not-irreducible}
Let $\Delta=(2)$. Figure~\ref{fig:p-n-irr} shows two oriented paths $P_1$
and~$P_2$ and their product~$P_1\times P_2$. We can see that the core~$P$
of~$P_1\times P_2$ is an oriented path. Thus $P$~is homomorphically
equivalent to the product of two incomparable structures $P_1$ and $P_2$.
It follows from Lemma~\ref{lem:irreducible} that $P$~is not irreducible.
\picture{duality.3}{A non-irreducible path}{p-n-irr}
\end{exa}

\section{Finite dualities}
\label{sec:findual}

\subsection*{Introduction}

In this section, we generalise duality pairs in a natural way: instead of forbidding homomorphisms
from a single structure, we forbid homomorphisms from a finite set of structures; and on the
other side, we allow structures to map to some of a finite number of structures.

\begin{dfn}
Let $\F$ and~$\D$ be two finite sets of core \Ds s such that no
homomorphisms exist among the structures in~$\F$ and among the structures
in~$\D$. We say that $(\F,\D)$ is a \deph{finite homomorphism duality}
(often just a \deph{finite duality}) if for every \Ds~$A$
there exists $F\in\F$ such that $F \to A$
if and only if for all $D\in\D$ we have $A\notto D$.
\end{dfn}

\begin{pgf}
For a symbol like $\toDd$, two definitions are possible. It may denote
either the set of all structures that admit a homomorphism to \emph{some}
structure in~$\D$, or the set of structures that map to \emph{each}
structure in~$\D$.

It is convenient for us to use the positive symbols $\toDd$ and $\Ffto$
in the former sense, while the negative symbols $\Ffnto$ and $\ntoDd$
will be used in the latter sense. So we define
\begin{align*}
\Ffto &:= \{ A : F\to A \text{ for some } F\in\F\},\\
\Ffnto &:= \{ A : F\notto A \text{ for all } F\in\F\},\\
\noalign{\allowbreak}
\toDd  &:= \{ A : A\to D \text{ for some } D\in\D\},\\
\ntoDd  &:= \{ A : A\notto D \text{ for all } D\in\D\}.
\end{align*}

In this notation, $(\F,\D)$ is a finite homomorphism duality if and only if
\[\Ffnto=\toDd.\]

Since the classes $\Ffto$ and $\Ffnto$ are complementary and so are the
classes $\toDd$ and $\ntoDd$, the pair $(\F,\D)$ is a duality pair if
and only if
\[\Ffto=\ntoDd.\]
\end{pgf}

We begin exploring the world of finite dualities by considering
situations where the set~$\D$ has only one element. Here we remark that
forbidding homomorphisms from a finite set of structures is equivalent to
forbidding a finite number of substructures. This is again characteristic
of finite dualities, since all the classes~$\toD$ are characterised by
forbidden subgraphs (although in many cases by infinitely many of them).

\begin{prop}
Let $D$ be a core \Ds. Then the following are equivalent:
\begin{itemize}
\item[(1)] There exists a finite set~$\F$ of \Ds s such that the pair
$\bigl(\F,\{D\}\bigr)$ is a finite homomorphism duality.
\item[(2)] There exists a finite set~$\F'$ of \Ds s such that any \Ds~$A$ is homomorphic to~$D$
if and only if it contains no element of~$\F'$ as its substructure.
\end{itemize}
\end{prop}

\begin{proof}
Suppose (1) holds and set $\F'$ to be the set of all homomorphic images
of structures in~$\F$. Then (2) follows from the definition of duality.

Conversely, if (2) holds, let $\F$ be the set of all cores of the
structures in~$\F'$.  If $A\notto D$, then $A$~contains some element
$F\in\F'$ as a substructure, so the core of~$F$ is homomorphic to~$A$. And
if $A\to D$ and $F\in\F$, then $F\notto A$, for otherwise $F\to D$,
a contradiction with~(2).  Therefore $\bigl(\F,\{D\}\bigr)$ is a finite
duality.
\end{proof}

Next we characterise all finite dualities $(\F,\D)$ with a singleton
right-hand side, that is dualities such that $|\D|=1$.

\begin{thm}[\cite{NesTar:Dual}]\label{thm:finitary}
If $\bigl(\F,\{D\}\bigr)$ is a finite homomorphism duality, then all elements of~$\F$ are \Dt s
and 
\begin{equation}\label{eq:2.1}
D\heq\prod_{F\in\F} D(F).\end{equation}
Conversely, for any finite collection~$\F$ of \Dt s, if \eqref{eq:2.1} holds,
then the pair $\bigl(\F,\{D\}\bigr)$ is a finite homomorphism duality.
\end{thm}

\begin{proof}
For the proof we need some terminology and results of Section~\ref{sec:gaps}.

Let $\bigl(\F,\{D\}\bigr)$ be a finite duality.  Suppose some $F\in\F$
is disconnected, $F=F_1+F_2$. Then no element of~$\F$ is homomorphic to any of
$F_1$ and $F_2$, so $F_1\to D$ and $F_2\to D$, and also $F=F_1+F_2\to D$,
a contradiction with the fact that $\bigl(\F,\{D\}\bigr)$ is a duality.
Therefore all elements of~$\F$ are connected.

Let $F\in\F$. Suppose that $D\to A\to F+D$, but
that $A\notto D$. By duality, there exists $F'\in\F$ such that $F'\to A$. Since $F'$ is
connected, $F'\notto D$, and $A\to F+D$, we have that $F'\to F$. Hence
$F'=F$ because distinct elements of~$\F$ are incomparable. Therefore
$A\notto D$ and $D\to A\to F+D$ implies that $F+D\to A$. This proves that $(D,F+D)$ is a gap,
as $D<F+D$ because $F\notto D$ by duality.

By Proposition~\ref{prop:all-gaps}, there exists a duality pair~$(T',D')$ such
that $T'\to F+D\to T'+D'$ and $D\heq (F+D)\times D'$ (see the following diagram, in which
gaps are marked by double arrows).
\[\xy\xymatrix{
&&T'+D'\\
&F+D \ar[ru]\\
T' \ar[ru]&&D' \ar@{{}=>}[uu]\\
&D\heq (F+D)\times D' \ar@{{}=>}[uu] \ar[ur]\\
T'\times D' \ar@{{}=>}[uu] \ar[ur]
}\endxy\]
If $F\to D'$, then by duality $T'\notto F$, and because $T'$ is connected
(it is in fact a \Dt), $T'\to D\to D'$, a contradiction with $(T',D')$
being a duality. So $F\notto D'$.

Since $F$~is connected and $F\notto D'$, we get that $F\to T'$; moreover, from duality we know
that $T'\to F$. We conclude that $T'\heq F$ and $D'\heq D(F)$ is the dual of~$F$.
In this way, we have proved that $D\to D(F)$ for all $F\in\F$, and hence also
$D\to\prod_{F\in\F} D(F)$.

On the other hand, $F\notto\prod_{F\in\F} D(F)$ for any $F\in\F$,
and so the product $\prod_{F\in\F} D(F)$ is homomorphic to~$D$ because
$\bigl(\F,\{D\}\bigr)$ is a finite duality.
Therefore $D\heq\prod_{F\in\F} D(F)\to D$.

Conversely, if \eqref{eq:2.1}~holds then  $F\notto\prod_{F\in\F} D(F)$
for any $F\in\F$, and in addition
if $F\notto A$ for any $F\in\F$ then $A\to\prod_{F\in\F} D(F)$. Hence
$\Bigl(\F,\prod_{F\in\F} D(F)\Bigr)$ is a finite duality.
\end{proof}

\begin{cor}
If $\F$ is a finite set of \Dt s, then there exists a unique core~$D$ such that
$\bigl(\F,\{D\}\bigr)$ is a finite homomorphism duality.
\qed\end{cor}

\begin{pgf}
This uniquely determined dual core~$D$ is denoted by~$D(\F)$.
\end{pgf}

\subsection*{Transversal construction}

Having characterised all finite dualities with a singleton right-hand side, we carry on by
providing a construction of finite dualities. Later we will see that all finite dualities
result from this construction.

The construction, which we will call the \emph{transversal construction},
starts with a finite set of \Df s. The forests are decomposed into
components and we consider sets consisting of the components. Some of
these sets satisfy certain properties and are called \emph{transversals}.
Each transversal is a set of \Dt s. The dual side of the finite duality
is then constructed by taking the dual structures for each transversal
(structures from Theorem~\ref{thm:finitary}).

\begin{pgf}
\label{pgf:247}
Let $\F=\{F_1, F_2, \dotsc, F_m\}$ be an arbitrary fixed non-empty finite set of core
\Df s that are pairwise incomparable ($F_j\notto F_k$ for $j\ne k$).
Let $\F_c = \{C_1,\dotsc,C_n\}$ be the set of all distinct connected
components of the structures in~$\F$; each of these components is a
core \Dt.
\end{pgf}

First we define \emph{quasitransversals} to be certain sets of components
appearing in the structures in~$\F$.

\begin{dfn}
\label{dfn:quasi}
A subset $\M\subseteq \F_c$ is a \deph{quasitransversal} if it satisfies
\begin{itemize}
\item[(\textsc{t}1)] any two distinct elements of~$\M$ are incomparable,
and
\item[(\textsc{t}2)] $\M$ \deph{supports} $\F$, that is for every $F\in\F$ there exists
$C\in\M$ such that $C\to F$.
\end{itemize}
\end{dfn}

\begin{pgf}
For two quasitransversals $\M$, $\M'$ we define that $\M\preceq \M'$ if and only if
for every $C'\in \M'$ there exists $C\in \M$ such that $C\to C'$. Note that this order is
different from the homomorphism order of forests corresponding to the quasitransversals.
On the other hand, we have the following.
\end{pgf}

\begin{lemma}\label{lem:preceq}
Let $\M$, $\M'$ be two quasitransversals. Then the dual structures $D(\M)$ and $D(\M')$ exist, and
$D(\M) \to D(\M')$ if and only if $\M\preceq \M'$.
\end{lemma}

\begin{proof}
By Theorem~\ref{thm:finitary}, the dual structures $D(\M)$ and $D(\M')$ exist and
\[D(\M)=\prod_{C\in \M} D(C),\quad D(\M')=\prod_{C'\in \M'} D(C').\]

Let $\M\preceq \M'$; we want to show that $D(\M)\to D(\M')$.
By the properties of product, it suffices to show that $D(\M)\to
D(C')$ for any $C'\in \M'$. So, let $C'\in \M'$. Because $\M\preceq \M'$,
there exists $C\in \M$ such that $C\to C'$. By the definition of a
duality pair, $C\to C'$ implies that $C'\notto D(C)$ and this implies that
$D(C)\to D(C')$. We conclude that $D(\M)\to D(C) \to D(C')$.

For the converse implication, let $D(\M)\to D(\M')$. We want to show that
for any $C'$ in~$\M'$ there is $C$ in~$\M$ with $C\to C'$. Indeed, for
$C'\in \M'$ we have $D(\M)\to D(\M')\to D(C')$; using duality, $C'\notto
D(\M)$, and therefore $C'\notto D(C)$ for some $C\in \M$. By duality
$C\to C'$.
\end{proof}

\begin{lemma}\label{lem:quasitransorder}
The relation $\preceq$ is a partial order on the set of all quasitransversals.
\end{lemma}

\begin{proof}
Obviously, $\preceq$ is both reflexive and transitive (a preorder).

Suppose now that
$\M\preceq \M'$ and $\M'\preceq \M$, and let $C\in \M$. Then there exists $C'\in
\M'$ such that $C'\to C$ and there exists $C''\in \M$ such that $C''\to C'$.
Consequently $C''\to C$, hence by (\textsc{t}1) we have $C=C'=C''$, so $\M\subseteq
\M'$. Similarly we get that $\M'\subseteq \M$. Therefore $\M=\M'$ whenever $\M\preceq\M'$ and
$\M'\preceq\M$. So $\preceq$~is antisymmetric; it is a partial order.
\end{proof}

\begin{dfn}
\label{dfn:trans}
A quasitransversal~$\M$ is a \deph{transversal} if
\begin{itemize}
\item[(\textsc{t}3)] $\M$ is maximal with respect to the order~$\preceq$.
\end{itemize}
\end{dfn}

\begin{pgf}
\label{pgf:transcons}
Set $\D=\D(\F)=\bigl\{D(\M) : \text{$\M$ is a transversal}\bigr\}$.
\end{pgf}

\begin{lemma}[Transversal construction works]
\label{lem:constr}
The pair $(\F,\D)$ is a finite homomorphism duality.
\end{lemma}

Before presenting the proof, we illustrate the construction by three examples.

\begin{exa}
First, suppose that $\F=\{T_1, T_2,\dotsc, T_n\}$ is a set of pairwise incomparable
trees and let $D_1$, $D_2$,~\dots,~$D_n$ be their respective duals. By~(\textsc{t}2),
every quasitransversal contains all these trees. Therefore there exists
only one quasitransversal $\M=\{T_1, T_2, \dotsc, T_n\}$ and it is a
transversal. So $\D=\bigl\{D(\M)\bigr\}=\{D_1\times D_2\times\dotsb\times
D_n\}$. This corresponds to the situation in Theorem~\ref{thm:finitary}.
\end{exa}

\begin{exa}
Now, let $T_1$, $T_2$, $T_3$ and $T_4$ be pairwise incomparable trees with
duals $D_1$, $D_2$, $D_3$, $D_4$. Let $\F=\{T_1+T_2, T_1+T_3, T_4\}$.
The partial order~$\preceq$ of quasitransversals is depicted
in the following diagram:
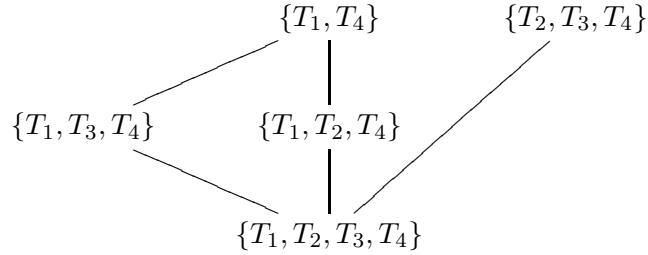
\[\xy\xymatrix{
&\{T_1,T_4\}&\{T_2,T_3,T_4\}\\
\{T_1,T_3,T_4\}\ar@{{}}[ur]&\{T_1,T_2,T_4\}\ar@{{}}[u]\\
&\{T_1,T_2,T_3,T_4\}\ar@{{}}[ul]\ar@{{}}[u]\ar@{{}}[uur]
}\endxy\]
We have two transversals $\{T_1, T_4\}$ and $\{T_2, T_3, T_4\}$; and
$\D=\{D_1\times D_4, D_2\times D_3\times D_4\}$.
\end{exa}

\begin{exa}
Finally, let $T_1\to T_3$ and $\F=\{T_1+T_2, T_3+T_4\}$.
This time, we get the following order of quasitransversals:
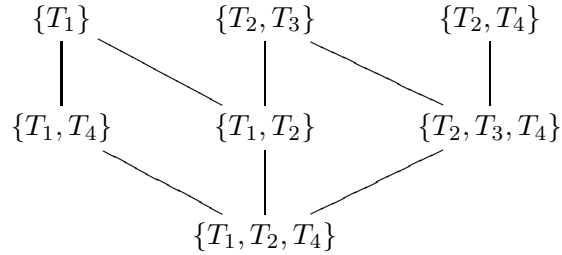
\[\xy\xymatrix{
\{T_1\}&	\{T_2,T_3\}&	\{T_2,T_4\}\\
\{T_1,T_4\}\ar@{{}}[u]&
	\{T_1,T_2\}\ar@{{}}[ul]\ar@{{}}[u]& \{T_2,T_3,T_4\}\ar@{{}}[ul]\ar@{{}}[u]\\
&\{T_1,T_2,T_4\}\ar@{{}}[ul]\ar@{{}}[u]\ar@{{}}[ur]
}\endxy\]
The transversals are $\{T_1\}$, $\{T_2, T_3\}$ and $\{T_2, T_4\}$. Hence
$\D=\{D_1,\mskip0mu plus 3mu D_2\times D_3,\penalty0\mskip0mu plus 3mu {D_2\times D_4}\}$.
\end{exa}

\begin{proof}[Proof of Lemma~\ref{lem:constr}]
By the definition of $\F$, any two distinct elements of~$\F$ are
incomparable.  Any two distinct elements of~$\D$ are incomparable too,
because any two transversals are incomparable with
respect to~$\preceq$ (they are all maximal in this order) and because
of Lemma~\ref{lem:preceq}.

Let $X$ be a \Ds\ such that $X\to D$ for some $D\in \D$. We want to prove
that $F_i\notto X$ for $i=1,\dotsc,m$. To obtain contradiction, assume that
$F_i\to X$ for some~$i$. Let $\M$ be the transversal for which $D(\M)=D$.
By~(\textsc{t}2), there exists $C\in \M$ such that $C\to F_i\to X$, therefore
$X\notto D(C)$.  This is a contradiction with the assumption that $X\to
D\to D(C)$ (here $D\to D(C)$ because $D$ is the product of the duals of the structures
in~$\M$, the component $C$ is an element of~$\M$, and the projection is a homomorphism).

Now, let $X$ be a \Ds\ such that $F_i\notto X$ for $i=1,\dotsc,m$. We want
to prove that there exists $D\in \D$ such that $X\to D$. Let $C_{j_i}$ be
a component of~$F_i$ such that $C_{j_i}\notto X$ for $i=1,\dotsc,m$. Let
$\M'=\min_{\to}\{C_{j_i}:i=1,\dots,m\}$, where by $\min_{\to} S$ we mean
the set of all elements of~$S$ that are minimal with respect to the
homomorphism order~$\to$. Because $\M'$ is a quasitransversal, there exists
a transversal~$\M$ such that $\M'\preceq \M$. We have that $C\notto X$ for
each $C\in \M$, and thus $X\to D(\M)\in\D$.
\end{proof}

\subsection*{Characterisation}

We will now prove that actually all finite homomorphism dualities are obtained
from the transversal construction.

\begin{pgf}
Let $(\F,\D)$ be a finite homomorphism duality.
Suppose
$\F=\{F_1,F_2,\dotsc,F_m\}$ and $\D=\{D_1,D_2,\dotsc,D_p\}$.
By definition, we assume that all the structures in~$\F$ and also all the structures
in~$\D$ are pairwise incomparable cores.  Consistently with the above
notation, let $\F_c=\{C_1,C_2,\dots,C_n\}$ be the set of all distinct
connected components of the structures in~$\F$.  Quasitransversals and
transversals are defined in the same way as above; notice that neither
for their definition nor for proving Lemma~\ref{lem:quasitransorder} we needed
the fact that the elements of~$\F_c$ are trees.
\end{pgf}

\begin{pgf}
For a quasitransversal $\M$, let $\overline{\M}=\{C'\in \F_c: C\in \M
\Rightarrow C\notto C'\}$ be the set of all components ``not supported''
by~$\M$.
\end{pgf}

\begin{lemma}\label{lem:f1}
If $\M\subseteq \F_c$ is a transversal, then there exists a unique \Ds\ ${D\in \D}$
that satisfies
\begin{itemize}
\item[(1)] $C\notto D$ for every $C\in\M$,
\item[(2)] $C'\to D$ for every $C'\in\overline\M$.
\end{itemize}
\end{lemma}

\begin{proof}
If $\overline\M=\emptyset$, let $D\in \D$ be arbitrary. Otherwise set
$S=\coprod_{C'\in\overline\M}C'$. Because $(\F,\D)$ is a finite homomorphism duality, either
there exists $F\in\F$ such that $F\to S$ or there exists $D\in\D$ such that $S\to
D$. If $F\to S$, by (\textsc{t}2) some $C\in\M$ satisfies $C\to F\to
S$, and since $C$ is connected, $C\to C'$ for some $C'\in\overline\M$, which is a
contradiction with the definition of~$\overline\M$. Therefore there exists $D\in
\D$ that satisfies $S\to D$.

Obviously, such $D$ satisfies (2).

Further, we will prove that $D$~satisfies~(1) as well.  For the sake of
contradiction, suppose that there is $C\in\M$ such that $C\to D$.

Consider $\M'=\M\setminus\{C\}$. The
set~$\M'$ is not a quasitransversal, because otherwise we would
have $\M\prec \M'$ and $\M$ would not satisfy~(\textsc{t}3). Hence $\M'$ fails to
satisfy~(\textsc{t}2), and we can find $F\in \F$ which is not supported by~$\M'$. It
follows that $C\to F$.

Consider $\Q'$, the set of all elements of~$\F$ that are not supported
by~$\M'$. We know that $\Q'$ is non-empty because $F\in \Q'$.

There exists $F'\in \Q'$ such that $C$ is a connected component of~$F'$:
otherwise let $\M^\ast$ be the set of all components $C^\ast$ of \Ds s
in~$\Q'$ such that $C\to C^\ast$, and let $\M'':=\min_\to(\M'\cup \M^\ast)$
be the set of all structures in the union of $\M'$ and $\M^\ast$ that
are minimal with respect to the homomorphism order. The set $\M''$ is a
quasitransversal but $\M\prec\M''$, contradicting the fact that $\M$~is
a transversal.

All the components of~$F'$ are elements of $\overline\M\cup\{C\}$. The assumption
that $C\to D$ leads, using (2), to the conclusion that $F'\to D$. That is a
contradiction with the definition of finite duality.

It remains to prove uniqueness. If $D,D'\in \D$ both satisfy (1) and
(2) and $D\ne D'$, that is $D\inc D'$, then $D+D'$ violates the definition
of finite homomorphism duality: $D+D'$~is homomorphic to no $\check D$ in~$\D$,
otherwise the elements of~$\D$ would not be incomparable, contradicting
the definition of finite duality; at the same time no $F$ in~$\F$
is homomorphic to~$D+D'$, because (by the definition of a transversal) for
every $F\in\F$ there is $C\in M$ such that $C\to F$, but $C\not\to D+D'$,
because $C$ is connected and by~(1) it is homomorphic to neither $D$ nor~$D'$.
\end{proof}

\begin{pgf}
For a transversal $\M$, the unique $D\in \D$ satisfying the conditions (1)
and~(2) above is denoted by~$d(\M)$.
\end{pgf}

\begin{lemma}\label{lem:f2}
$\D=\bigl\{d(\M): \text{$\M$ is a transversal}\bigr\}$.
\end{lemma}

\begin{proof}
Let $D\in \D$. We want to show that $D=d(\M)$ for some
transversal~$\M$. Let $\M'=\min_\to\{C'\in \F_c: C'\notto D\}$ be the
set of all components that are not homomorphic to~$D$, minimal in the
homomorphism order. The set $\M'$ is a quasitransversal: if some $F\in
\F$ is not supported by~$\M'$, then all its components are homomorphic
to~$D$, and so $F\to D$, a contradiction.

Let $\M$ be a transversal such that $\M'\preceq \M$. To prove that $D=d(\M)$,
it suffices to check conditions
(1) and~(2) of Lemma~\ref{lem:f1}.

If $C\in\M$, then there exists $C'\in\M'$ such that $C'\to C$. Therefore
$C\notto D$, so condition~(1)~is satisfied.

Now condition~(2): Suppose on the contrary that there exists $\check
C\in\overline\M$ such that $\check C\notto D$.
Consider the \Ds~$X=\check C+D$. If $F\to X$ for some $F\in \F$, then by the
property (\textsc{t}2) of~$\M$ there exists $C\in\M$ that is homomorphic to~$F$. But
since $\check C\in\overline\M$, we have that $C\notto \check C$, hence $C\to D$.
This is a contradiction with the condition~(1). It follows that $X\to\check D$
for some $\check D\in \D$, hence $D\to\check D$, so $D=\check D$. That is a
contradiction with $\check C\notto D$ and $\check C\to\check D$.
\end{proof}

\begin{lemma}\label{lem:f3}
For two distinct transversals $\M_1$, $\M_2$, we have
\begin{itemize}
\item[(a)] $\overline{\M_1}\cap{\M_2}\ne\emptyset$,
\item[(b)] $d(\M_1)\notto d(\M_2)$.
\end{itemize}
\end{lemma}

\begin{proof}
\mbox{}

(a) By (\textsc{t}3), $\M_1\npreceq \M_2$, and therefore there exists $C_2\in
\M_2$ such that $C_1\notto C_2$ for any $C_1\in \M_1$. Obviously $C_2\in
\overline{\M_1}\setminus\overline{\M_2}\subseteq\overline{\M_1}$. Since we
selected $C_2\in \M_2$, we have that $C_2\in\overline{\M_1}\cap \M_2$.

(b) Let $C_2\in\overline{\M_1}\cap \M_2$, as above.  Then $C_2\to d(\M_1)$
and $C_2 \notto d(\M_2)$. Consequently $d(\M_1)\not\to d(\M_2)$.
\end{proof}

\begin{lemma}\label{lem:f4}
If $\M$ is a transversal, then the pair $\bigl(\M,\{d(\M)\}\bigr)$ is a
finite homomorphism duality, and consequently $d(\M)=D(\M)$.
\end{lemma}

\begin{proof}
We want to prove that 
  \[{\M{\notto}} = {{\to}d(\M)}.\]

We claim that for a \Ds~$A$, the following statements are equivalent:

(1) $A\in {\M{\notto}} = \bigcap_{C\in\M}({C\notto})$

(2) $C\notto A$ for any $C\in\M$

(3) $C\notto A+\coprod_{\check C\in\overline\M}\check C$ for any $C\in\M$

(4) $A+\coprod_{\check C\in\overline\M}\check C \to d(\M)$

(5) $A\to d(\M)$

(6) $A\in{{\to}d(\M)}$

\noindent
Because: $(1)\Leftrightarrow(2)$ and $(5)\Leftrightarrow(6)$ by definition.
$(4)\Rightarrow(5)$ immediately.
$(5)\Rightarrow(2)$ by Lemma~\ref{lem:f1}(1).
$(2)\Rightarrow(3)$ follows from the definition of~$\overline\M$ and the fact that
$C$~is connected.

It remains to prove that $(3)\Rightarrow(4)$: Let $X=A+\coprod_{\check
C\in\overline\M}\check C$. If $F\to X$ for some $F\in \F$, then by (\textsc{t}2) there
exists $C\in\M$ such that $C\to F\to X$, a contradiction. Thus no element
of~$\F$ is homomorphic to~$X$, hence $X\to D$ for some $D\in \D$. By
Lemma~\ref{lem:f2}, $D=d(\M')$ for a transversal~$\M'$; by Lemma~\ref{lem:f1} and
Lemma~\ref{lem:f3}(a), $\M'=\M$.

The equivalence $\text{(1)}\Leftrightarrow\text{(6)}$ is precisely the
definition of finite duality.

By Theorem~\ref{thm:finitary}, the dual is uniquely determined if it is a core,
so $d(\M)=D(\M)$.
\end{proof}

Lemma~\ref{lem:f4} and Theorem~\ref{thm:dual-pairs} imply that any element of a
transversal is a \Dt, but we have not proved that every structure
in~$\F_c$ is an element of some transversal. However, we have the following lemma, for
whose proof we will once again use the characterisation of gaps in Section~\ref{sec:gaps}.

\begin{lemma}\label{lem:f5}
Each component $C\in \F_c$ is a \Dt.
\end{lemma}

\begin{proof}
Suppose that $C\in\F_c$ is not a tree. By Lemma~\ref{lem:f4},
Theorem~\ref{thm:dual-pairs}, and Theorem~\ref{thm:finitary}, $C$~is an
element of no transversal. Set
\[A=\coprod_{\substack{C'\in\F_c \\  C'< C}} C' +
\coprod_{\substack{C'\in\F_c \\ C'\inc C}} (C\times C').\]

Clearly, $A<C$ because all the summands are less than~$C$ and $C$~is
connected.  As $C$~is not a tree, it has no dual; because it is connected,
the pair $(A,C)$ is not a gap by Lemma~\ref{lem:gap-dual}. Let $X$
be a structure satisfying that $A<X<C$.

Then for any $C'\in\F_c$ such that $C\ne C'$, we have $C'\to X$ if and
only if $C'\to C$ and $X\to C'$ if and only if $C\to C'$. Indeed: if $C'\to
C$, then $C'\to A\to X$; if $C\to C'$, then $X\to C\to C'$. On the other hand, if
$C\inc C'$, then $X\to C'$ implies $X\to C\times C'\to A$ (because
$C\times C'$ is one of the summands in the above definition of~$A$),
a contradiction with $A<X$.  Moreover $C'\to X$ implies $C'\to C$.

Let $F\in\F$ be such that $C$ is a component of~$F$ and let $G$ be the
structure obtained from~$F$ by replacing $C$ with~$X$.

Suppose $F\to G$. Then $C\to G$. Because $C$ is connected, it
is homomorphic to a component of~$G$. Since $C\not\to X$, it is
homomorphic to some other component of~$F$, contradicting that $F$
is a core. Therefore $F\not\to G$.

In addition, $F'\notto G$ for any $F\ne F'\in\F$, because $F'\to G$
implies $F'\to F$. Therefore $G\to D$ for some $D\in\D$. Let $M$ be the
transversal such that $D=D(M)$. Recall that $C$ is an element of no
transversal, so $C\not\in M$. The
structure $D$ is a product of duals and hence $C'\notto G$ for any $C'\in M$;
therefore $C'\notto X$ and $C'\notto C$ for any $C'\in M$. Consequently
$C\to D$. We know that all components of~$G$ are homomorphic to~$D$, so all
components of~$F$ are homomorphic to~$D$ as well. We conclude that $F\to
D$, a contradiction.
\end{proof}

We finish this section by a theorem that characterises all finite dualities.

\begin{thm}[Characterisation of finite dualities]\label{thm:finite-character}
If $(\F,\D)$ is a finite homomorphism duality, then all elements of~$\F$
are \Df s and $\D=\D(\F)$ results from the transversal construction.
In particular, $\D$~is determined by~$\F$ uniquely up to homomorphic
equivalence.

Conversely, for any finite collection~$\F$ of core \Df s,
$\bigl(\F,\D(\F)\bigr)$ is a finite homomorphism duality.
\end{thm}

\begin{proof}
All elements of~$\F$ are forests because of of Lemma~\ref{lem:f5}.  The set
$\D$ is uniquely determined as a consequence of Lemma~\ref{lem:f2}
and because of Lemma~\ref{lem:f4} and Theorem~\ref{thm:finitary} it is
determined by the transversal construction.

The second part is Lemma~\ref{lem:constr}.
\end{proof}

Now we can view the notation $(\F,\D)$ from a different perspective:
the letter~$\F$ stands for \emph{forbidden}, as we mentioned in~\ref{pgf:fdnot},
but it may also be understood to stand for \emph{forests}.

\begin{bfdtf}
Our construction of duals from forests relied heavily on the fact that
every finite \Ds\ is a finite sum of components, structures that are
connected. Although we mentioned in~\ref{pgf:nondecomposable} that some
structures are not a finite product of irreducible structures. However, it
can be shown that the set~$\F$ in a duality pair~$(\F,\D)$ is determined
uniquely by~$\D$ too.

The characterisation of finite homomorphism dualities implies that dual
structures can be factored into a product of irreducible structures. A
construction dual to the transversal construction produces the forests
from the dual set. This is covered in more detail when we discuss a
complexity issue in Section~\ref{sec:decdual}.
\end{bfdtf}

\section{Extremal aspects of duality}
\label{sec:extrem}

Extremal theories are concerned with questions how large can an object be if it satisfies
certain conditions, or has certain properties.

\begin{pgf}
In the context of homomorphism dualities, we are interested in the following four questions.
\begin{itemize}
\item[(1)] Given a \Dt, how large can its dual be?
\item[(2)] Given a right-hand side of a duality pair, how large can the corresponding \Dt\ be?
\item[(3)] Given a finite set~$\F$ of \Df s, how large can the set~$\D(\F)$ be?
\item[(4)] Given a right-hand side of a finite duality, how large can the corresponding
left-hand side be?
\end{itemize}
Naturally, one has to define a suitable notion of size for this purpose.
\end{pgf}

Results have been published about questions (1) and~(2), partially also
about~(4). Finite dualities have been studied in full generality only
recently, so questions (3) and~(4) have not yet been thoroughly investigated.

\begin{pgf}
When examining extremal problems about homomorphism dualities, we consider
only cores. That is, we ask how large the \emph{core} of the dual of a \emph{core} \Dt\
is, and analogously for the other questions.

The reason for this is clear: in every class of homomorphic
equivalence there are arbitrarily large structures. The smallest
structure in such a class is the (unique) core in it.  By taking sums of
an arbitrary number of disjoint copies of the core, we can produce arbitrarily large
homomorphically equivalent structures.
\end{pgf}

Concerning question~(1), an upper bound on the size of a \Dt's dual
follows from the bear construction. The bound on the size of the base
set of the dual is exponential in terms of the size of the base set of
the \Dt.

\begin{thm}[\cite{NesTar:Short}]
\label{thm:dualsize}
Let $\Delta=(\delta_i:i\in I)$. Let $(F,D)$ be a duality pair
such that $D$~is a core and let $n:=|\bs F|$. Then \[|\bs D| \le
n^{n}.\]
\end{thm}

\begin{proof}
Let $F$ be a \Dt.  Consider the bear construction
from~\ref{pgf:bear-def}. A vertex of the dual~$D_b(F)$ is a function that assigns
each vertex~$x$ of~$F$ an edge~$e$ of~$F$ such that $x$~appears in~$e$.
Since $F$ is a tree, no two distinct edges contain more than one vertex in
common. Thus the number of edges containing a fixed vertex~$x$ is at
most~$n$.
Hence the number of vertices of~$D_b(F)$ is at most $n^n$, as we
were supposed to prove.
\end{proof}

Nešetřil and Tardif~\cite{NesTar:Short} also provide a construction
of paths whose duals indeed have exponential size.

\begin{thm}[{\cite[Theorem~8]{NesTar:Short}}]
For any sufficiently large positive integer~$N$ there exists a core
\Dt~$F$ such that $|\bs F|\ge N$ and if $D$~is the core for which
$(F,D)$~is a duality pair, then
\[|\bs D| \ge 2^{n/7\log_2 n},\]
where $n:=|\bs F|$.
\end{thm}

For question~(2), it is easier to measure the size of the forbidden \Dt\
in terms of its diameter.

We assume that the reader knows that the \emph{distance} of two vertices
in an undirected graph is the number of edges of the shortest path
connecting them, and that a graph's \emph{diameter} is the maximum
distance of a pair of its vertices. The notion of diameter of an
undirected graph is used to define the diameter for \Ds s.

\begin{dfn}
The \deph{diameter} of a \Ds~$A$ is half the diameter of its incidence graph~$\Inc(A)$.
\end{dfn}

Larose, Loten and Tardif~\cite{Tar:FOD} proved an upper bound on the
diameter of forbidden trees in a finite duality with a singleton right
hand side.

\begin{thm}
\label{thm:tree-bound}
If $\bigl(\F,\{D\}\bigr)$ is a finite duality, $F\in\F$ is a core,
and $n=|\bs D|$, then the diameter of~$F$ is at most $n^{n^2}$.
\qed
\end{thm}

This theorem implies a bound on the number of edges of such forbidden
trees, see Lemma~\ref{lem:edgebd}. However, this bound is very rough,
even though no examples are known that have an exponential number of
edges in terms of the number of vertices of the dual.

\setchapterpreamble[u]{%
	\dictum[Niccolo Machiavelli]{%
	There is nothing more~difficult to~take in~hand, more~perilous
	to~conduct or more~uncertain in~its~success than to~take
	the~lead in~the~introduction of~a~new~order of~things.}}

\chapter{Homomorphism order}
\label{chap:order}

\bigskip\bigskip

The relation of existence of a homomorphism on the class of all
\Ds s induces a partial order, called the \emph{homomorphism order}.

Properties of this partial order have been widely studied in algebraic,
category theory, random and combinatorial context. The homomorphism
order motivates several questions linked to problems of existence of
homomorphisms, see Section~\ref{sec:csp}.

Density- and universality-related issues have attracted special attention.
Universality of the homomorphism order of undirected graph was proved
already in 1969 by Hedrlín~\cite{Hed:Uni}; in particular, it was
shown that any countable partial order is an induced suborder of the
homomorphism order. Universality was also studied for special classes
of digraphs, and it has recently been proved that even the relatively
small class of all directed paths induces a universal countable partial
order~\cite{HubNes:Finite,HuNe:Paths}.

The examination of density has a long history too.  A complete description
of all non-dense parts, called \emph{gaps}, for the homomorphism order
of undirected graphs, was given in~1982 by Welzl~\cite{Wel:Dense}.

For directed graphs and general relational structures, density
has a non-obvious link to duality. It was shown by Nešetřil and
Tardif~\cite{NesTar:Dual} that all gaps correspond to duality pairs; we
survey their results in Section~\ref{sec:gaps}. This connection can be
extended from the homomorphism order to lattices satisfying some extra
axioms (Heyting algebras with finite connected decompositions). Such
an extension is presented in Section~\ref{sec:heyt}. Some
of the ideas are contained in a paper of Nešetřil, Pultr and
Tardif~\cite{NesPulTar:HeytDual}. We add the description of finite
dualities.

Our further interest concentrates on another issue. In
Section~\ref{sec:macs} we study finite maximal antichains in the
homomorphism order. In particular, we show that with a few characterised
exceptions finite maximal antichains have the splitting property. This
by itself provides a connection to finite homomorphism dualities but
in the case of relational structures with at most two relations we can
prove that even the exceptional antichains are formed from dualities. For
structures with more than two relations this question remains open.

\section{Homomorphism order}

\begin{pgf}
The relation $\to$ of being homomorphic is reflexive, as the identity
mapping is a homomorphism from a \Ds\ to itself, and it is transitive,
since the composition of two homomorphisms, if possible, is a homomorphism
too. Thus $\to$~is a preorder.

There are standard ways to transform a preorder into a partial order. It
may be done by identifying equivalent objects, or by choosing a particular
representative for each equivalence class. The resulting partial order
is identical in both cases.

For $\to$, a suitable representative for each equivalence class is
a core \Ds. We have already observed in~\ref{prop:allhavecore} that
there is a unique core in each class of homomorphic equivalence;
unique up to isomorphism.
\end{pgf}

\begin{prop}
\label{prop:homord}
Let $\Delta$ be a fixed type.  Then the relation $\to$ of being
homomorphic is a partial order on the set of all
core \Ds s (taken up to isomorphism).
\end{prop}

\begin{proof}
Follows from the discussion above.
\end{proof}

\begin{dfn}
The partial order~$\to$ from Proposition~\ref{prop:homord} is called the
\deph{homomorphism order} and denoted by~$\CD$.
\end{dfn}

Any treatise on the homomorphism order is substantially simplified
by talking about the order of \Ds s rather than cores or equivalence
classes. For instance, when we say that $A$~is less than~$B$ in the
homomorphism order, we mean that \emph{the core of~$A$} is less than
\emph{the core of~$B$} in the homomorphism order. Similarly, when we
(soon) say that $A\times B$ is the infimum of $A$ and $B$ in the
homomorphism order, we mean that \emph{the core of~$A\times B$} is
actually the infimum. This approach is fairly standard in algebra.

With all this in mind, we observe that the homomorphism order is a nice
partial order: it is a lattice, and moreover a Heyting algebra.

\begin{prop}
\label{prop:homheyt}
The homomorphism order~$\CD$ is a Heyting algebra. In particular, for $A,B\in\CD$
\begin{itemize}
\item[(1)] the product $A\times B$ is the infimum (meet) of $A$ and $B$,
\item[(2)] the sum $A+B$ is the supremum (join) of $A$ and $B$,
\item[(3)] one vertex with no edges
$\bigl(\{1\},(\emptyset,\emptyset,\dotsc\emptyset)\bigr)=:\bot$
is the least element, and one vertex with all loops
$\bigl(\{1\},(\{1\}^{\delta_i}:i\in I)\bigr)=:\top$ is the greatest
element in~$\CD$,
\item[(4)] the exponential structure $B^A$ is the Heyting operation $A\heyt B$.
\end{itemize}
\end{prop}

\begin{proof}
(1) By~\ref{cor:product}, $C\to A\times B$ if and only if $C\to A$
and $C\to B$. So $A\times B$ is the infimum of $A$ and~$B$.

(2) By~\ref{cor:sum}, $A+B\to C$ if and only if $A\to C$ and $B\to C$. So
$A+B$ is the supremum of $A$ and~$B$.

(3) Let $A$ be an arbitrary \Ds. By definition, $\bs A$ is non-empty
and clearly any function mapping~$v$ to an arbitrary element of~$\bs A$
is a homomorphism from~$\bot$ to~$A$. So $\bot$~is the least element.

On the other hand, let $f$ be the constant function from~$\bs A$
to~$\{v\}$ such that $f(a)=v$ for all $a\in\bs A$. Since $\top$~has
all loops, all edges of~$A$ are preserved by~$f$ and hence it is a
homomorphism from~$A$ to~$\top$. Therefore $\top$~is the greatest element.

(4) By~\ref{cor:expon}, $C\to B^A$ if and only if $A\times C\to
B$. Since products are infima, this is exactly the Heyting axiom (see
Definition~\ref{dfn:heytal}).
\end{proof}

To illustrate the homomorphism order's power, we give (without proof)
one more example of its properties. The homomorphism order is a universal
countable partial order. Several proofs of this can be found in the
literature~\cite{Hed:Uni,HubNes:Finite,HuNe:Paths,Nes:ColPos,PulTrn:Cat}.

\begin{thm}
Let $\Delta$ be a type with at least one relation of arity at least two.
Every countable partial order is an induced suborder of the homomorphism order of \Ds s.
\qed
\end{thm}

\section{Gaps and dualities}
\label{sec:gaps}

In this section we briefly survey the results of~\cite{NesTar:Dual}
about a connection between duality pairs (see Section~\ref{sec:pairs})
and gaps in the homomorphism order. The explicit description of gaps,
besides being of interest by itself, provides a different proof of the
characterisation of duality pairs (Theorem~\ref{thm:dual-pairs}). We
have also used it for proving Theorem~\ref{thm:finite-character}.

The first fact we state is that the top of a gap, if connected, is the
left-hand side of a duality pair.

\begin{lemma}[\cite{NesTar:Dual}]\label{lem:gap-dual}
Let $(A,B)$ be a gap pair and let $B$ be connected. Then $(B, A^B)$ is a duality pair.
\qed
\end{lemma}

Hence a connected top of a gap is (homomorphically equivalent) to a \Dt.

The following proposition characterises all gaps.

\begin{prop}[\cite{NesTar:Dual}]\label{prop:all-gaps}
Gaps are exactly all the pairs $(A,B)$ such that there exists a duality
pair $(F,D)$ with $F\to B\to F+D$ and $A\heq B\times D$. Moreover, $B\heq A+F$.
\qed
\end{prop}

The correspondence is depicted in the following two diagrams, in which
double arrows denote gaps. Here $(F,D)$~is a duality pair.
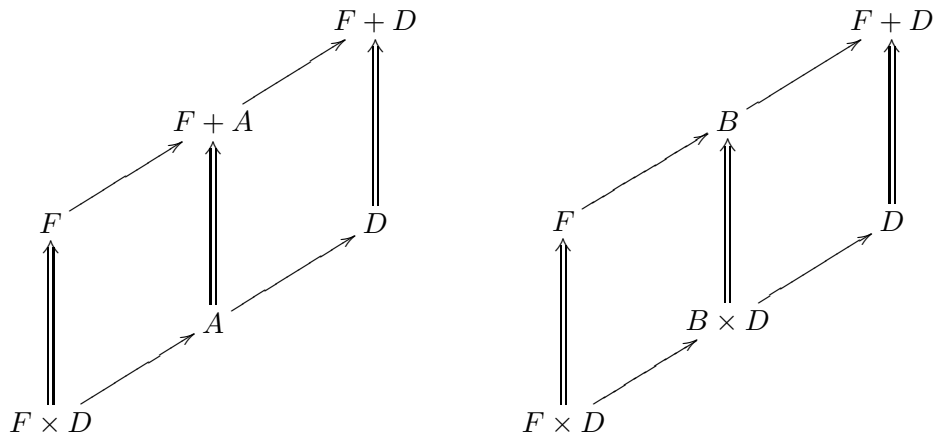
\begin{align*}
\xy\xymatrix{
&&F+D\\
&F+A \ar[ru]\\
F \ar[ru]&&D \ar@{{}=>}[uu]\\
&A \ar@{{}=>}[uu] \ar[ur]\\
F\times D \ar@{{}=>}[uu] \ar[ur]
}\endxy
&&
\xy\xymatrix{
&&F+D\\
&B \ar[ru]\\
F \ar[ru]&&D \ar@{{}=>}[uu]\\
&B\times D \ar@{{}=>}[uu] \ar[ur]\\
F\times D \ar@{{}=>}[uu] \ar[ur]
}\endxy
\end{align*}

\section{Dualities and gaps in Heyting algebras}
\label{sec:heyt}

The previous section presents a connection between finite dualities and
gaps in the homomorphism order. However, few properties typical
of the homomorphism order were used to prove them. Here we look at a more general
case. We provide conditions under which a theory of gaps and dualities
can be developed for partially ordered sets.

Gaps, duality pairs and combined dualities (which correspond to finite
dualities with a singleton right-hand side) in Heyting algebras have
been studied by Nešetřil, Pultr and Tardif~\cite{NesPulTar:HeytDual}.
We extend their results to a complete description of dualities.

At the same time, Proposition~\ref{prop:homheyt} implies that dualities
for relational structures are a special case of this general theory.

\begin{dfn}
\label{dfn:heytal}
A lattice~$P$ with an additional binary operation~$\heyt$ is a
\deph{Heyting algebra} if a least element and a greatest element exist
in~$P$ and for all $p,q,r\in P$,
\[p\preceq q\heyt r  \quad\text{if and only if}\quad  p\wedge q \preceq r.\]
\end{dfn}

Of course, not every lattice is distributive. However, it is well known
that every Heyting algebra is distributive.

\begin{lemma}
\label{lem:heyt-distr}
Every Heyting algebra is a distributive lattice.
\end{lemma}

\begin{proof}
First we show that $(a\wedge b) \vee (a\wedge c)
\preceq a\wedge(b\vee c)$. This is true in every lattice. Clearly $a\wedge
b\preceq a$ and $a\wedge b\preceq b\vee c$, so $a\wedge b\preceq a\wedge
(b\vee c)$. Similarly $a\wedge c\preceq a\wedge (b\vee c)$. Hence the
inequality holds.

Next, it suffices to prove that whenever $a\wedge b \preceq y$ and
$a\wedge c \preceq y$, then $a\wedge(b\vee c)\preceq y$. In connection
with the previous paragraph, it implies that $a\wedge(b\vee c)$ is the
supremum of $a\wedge b$ and $a\wedge c$.

So suppose that $a\wedge b\preceq y$ and $a\wedge c\preceq y$. Then
$b\preceq a\heyt y$ and $c\preceq a\heyt y$. Thus $b\vee c\preceq a\heyt
y$. Hence $a\wedge(b\vee c)\preceq y$.
\end{proof}

An important property for the development of duality theory for relational
structures was the existence of a decomposition of every relational
structure into connected components. We generalise connectedness in the
context of Heyting algebras.

\begin{dfn}
Let $L$ be a lattice. An element~$a$ of~$L$ is \deph{connected} if the equality
$a=b\vee c$ implies that $a=b$ or $a=c$.
\end{dfn}

Next we observe that connectedness in distributive lattices (and thus
in Heyting algebras) has the same equivalent descriptions (1)--(3)
as in Lem\-ma~\ref{lem:connected}.

\begin{lemma}\label{lem:heyt-conn}
Let $a$ be an element of a distributive lattice~$L$. Then the following
conditions are equivalent.

(1) If $a\preceq b\vee c$ for some elements $b$, $c$ of~$L$, then
$a\preceq b$ or $a\preceq c$.

(2) If $a=b\vee c$ for some elements $b$, $c$ of~$L$, then $b\preceq c$
or $c\preceq b$.

(3) The element $a$ is connected.
\end{lemma}

\begin{proof}\mbox{}

\mbox{(1)${}\Rightarrow{}$(2):}
If $a=b\vee c$, then $a\preceq b\vee c$, and using (1) we have $a\preceq
b$ or $a\preceq c$. In the first case $c\preceq b\vee c = a\preceq b$,
hence $c\preceq b$. In the latter case $b\preceq b\vee c = a\preceq c$,
and so $b\preceq c$.

\mbox{(2)${}\Rightarrow{}$(3):}
Suppose $a=b\vee c$. By (2) we have $b\preceq c$, and therefore $c=b\vee c=a$;
or we have $c\preceq b$, and then $b=b\vee c=a$.

\mbox{(3)${}\Rightarrow{}$(1):}
Let $a\preceq b\vee c$. Here we need distributivity: $(a\wedge
b)\vee(a\wedge c)= a\wedge(a\vee b)\wedge(a\vee c)=a$ since $a\preceq
a\vee b$ and $a\preceq a\vee c$.
\end{proof}

The existence of connected components is then generalised by the
following notion.

\begin{dfn}
We say that a lattice~$L$ \deph{has finite connected decompositions}
if each element~$x$ of~$L$ is a supremum of a finite set of connected
elements.
\end{dfn}

Analogously as for relational structures (Definition~\ref{dfn:dpair})
we define duality pairs for lattices.

\begin{dfn}
A pair $(f,d)$ of elements of a lattice~$L$ is a \deph{duality pair} if for any element
$x\in L$,
\[f\preceq x  \quad\text{if and only if}\quad  x\npreccurlyeq d.\]
\end{dfn}

\begin{dfn}
An element $f$ of a lattice~$L$ is called a \deph{primal} if there
exists $d\in L$ such that $(f,d)$~is a duality pair.  An element $d$
of a lattice~$L$ is called a \deph{dual} if there exists $f\in L$ such
that $(f,d)$~is a duality pair.
\end{dfn}

The next proposition is an analogue of~\ref{pgf:primal-connected}.

\begin{prop}[\cite{NesPulTar:HeytDual}]
In a distributive lattice, every primal is connected.
\end{prop}

\begin{proof}
We prove that if $(f,d)$ is a duality pair and $f\preceq b\vee c$
for some elements $b,c\in L$, then $f\preceq b$ or $f\preceq c$. By
Lemma~\ref{lem:heyt-conn} it follows that $f$~is connected.

So suppose that $f\preceq b\vee c$. By duality, $b\vee c\npreceq d$,
thus (by a property of join) $b\npreceq d$ or $c\npreceq d$. Using
duality once again we get that $f\preceq b$ or $f\preceq c$.
\end{proof}

Recall that a gap in a poset~$L$ is a pair $(p,q)$ of elements of~$L$
such that $p\prec q$ and no element~$r$ satisfies that $p\prec r\prec q$
(Definition~\ref{dfn:gap}). The connection between gaps and duality pairs
(Proposition~\ref{prop:all-gaps} for relational structures) is as follows.

\begin{thm}[\cite{NesPulTar:HeytDual}]
The gaps in a Heyting algebra~$L$ with finite connected decompositions
are exactly the pairs $(a,b)$ such that for some duality pair $(f,d)$
\[f\wedge d \preceq a \preceq d \quad \text{and} \quad b=a\vee f.\]
\qed
\end{thm}

\begin{dfn}
A pair $(F,D)$ of finite subsets of a lattice~$L$ is a \deph{finite duality} if
\begin{enumerate}
\item $f\npreccurlyeq f'$ if $f,f'\in F$ and $f\ne f'$,
\item $d\npreccurlyeq d'$ if $d,d'\in D$ and $d\ne d'$, and
\item for any $x\in L$ there exists $f\in F$ such that $f\preceq x$ if and only if
$x\npreccurlyeq d$ for any $d\in D$.
\end{enumerate}
\end{dfn}

The following theorem is an analogue of Theorem~\ref{thm:finitary}. It describes
finite dualities with a singleton right-hand side. The proof
is just a translation of the relational-structure proof of
Theorem~\ref{thm:finitary} into the language of Heyting algebras,
therefore we do not repeat it here.

\begin{thm}
\label{thm:heytfinitary}
Let $L$ be a Heyting algebra with finite connected decompositions. For
a finite subset~$F=\{f_1,f_2,\dotsc,f_n\}$ of~$L$ and for an element
$d\in L$, the pair $(F,\{d\})$ is a finite duality if and only if there
exist elements $d_1,d_2,\dotsc,d_n$ such that $(f_i,d_i)$ is a duality
pair for $i=1,2,\dotsc,n$ and $d=d_1\wedge d_2\wedge\dotsb\wedge d_n$.
\qed
\end{thm}

The \deph{transversal construction} of dualities in a lattice~$L$
with finite connected decompositions is defined analogously
to the definition in the context of \Ds s, contained in
paragraphs~\ref{pgf:247}--\ref{pgf:transcons}.

\begin{dfn}
Let $L$ be a lattice with finite connected decompositions.

Let $F=\{f_1, f_2, \dotsc, f_m\}$ be an arbitrary fixed non-empty finite
set of pairwise incomparable elements of~$L$. For each element~$f_i$
of~$F$ fix a finite connected decomposition
\[f_i=\bigvee_{j=1}^{k_i}c_{i,j}.\]
Now let $F_c$ be the set of all connected elements appearing in the decompositions, that is
\[F_c:=\bigcup_{i=1}^{m}\{c_{i,j}: 1\le j\le k_i\}.\]
Then $F_c$~is called a \deph{set of components} for~$F$.
\end{dfn}

Quasitransversals are defined analogously to quasitransversals for
relational structures (see~\ref{dfn:quasi}).

\begin{dfn}
A subset $M\subseteq F_c$ is a \deph{quasitransversal} if it satisfies
\begin{itemize}
\item[(\textsc{t}1)] any two distinct elements of~$M$ are incomparable,
and
\item[(\textsc{t}2)] $M$ \deph{supports} $F$, that is for every $f\in F$ there exists
$c\in M$ such that $c\preceq f$.
\end{itemize}
\end{dfn}

\begin{pgf}
For two quasitransversals $M$, $M'$ we define that $M\trianglelefteq M'$ if and only if
for every $c'\in M'$ there exists $c\in M$ such that $c\preceq c'$.
\end{pgf}

\begin{lemma}\label{lem:heytquasitransorder}
The relation $\trianglelefteq$ is a partial order on the set of all quasitransversals.
\end{lemma}

\begin{proof}
Obviously, $\trianglelefteq$ is both reflexive and transitive (a preorder).

Suppose now that $M\trianglelefteq M'$ and $M'\trianglelefteq M$, and let
$c\in M$. Then there exists $c'\in M'$ such that $c'\preceq c$ and there
exists $c''\in M$ such that $c''\preceq c'$.  As a result $c''\preceq
c$, thus by (\textsc{t}1) we have $c=c'=c''$, hence $M\subseteq M'$. Similarly
we get that $M'\subseteq M$. Hence $M=M'$ whenever $M\trianglelefteq M'$
and $M'\trianglelefteq M$. Therefore $\trianglelefteq$~is antisymmetric;
it is a partial order.
\end{proof}

\begin{dfn}
A quasitransversal~$M$ is a \deph{transversal} if
\begin{itemize}
\item[(\textsc{t}3)] $M$ is maximal with respect to the order~$\trianglelefteq$.
\end{itemize}
\end{dfn}

A characterisation similar to Theorem~\ref{thm:finite-character} follows.

\begin{thm}
Let $L$ be a Heyting algebra with finite connected decompositions.

Let $F\subseteq L$ be finite. If each element of~$F$ is a finite join
of primals, and $D$~is the result of the transversal construction,
then $(F,D)$~is a finite duality.

Conversely, if $(F,D)$~is a finite duality, then each element of~$F$
decomposes into a finite join of primals and $D$~is the result of the
transversal construction.
\qed
\end{thm}

In contrast to the homomorphism order, this does not in general mean that
the right-hand (dual) side of a finite duality is uniquely determined
by the left-hand (primal) side.  That is so because the decomposition
into connected components may not be
unique, so the transversal construction produces a different result for
different decompositions.

We do not give a detailed proof because the proof is a translation of the
proof of Theorem~\ref{thm:finite-character}. It suffices to check that
for proving the lemmas in Section~\ref{sec:findual} we did not use any
other properties of relational structures than the homomorphism order's
being a Heyting algebra with finite connected decompositions.
That is no longer true in the next section.

\section{Finite maximal antichains}
\label{sec:macs}

In this section we study finite maximal antichains in the homomorphism
order. In particular, we are interested in the \emph{splitting property}
of these antichains.

If $Q$~is an arbitrary maximal antichain in a poset~$P$, then every
element of~$P$ is comparable with some element~$q$ of~$Q$. In other
words, the poset~$P$ is the union of the downset generated by~$Q$
and the upset generated by~$Q$, that is $P=\ups{Q}\cup \downs{Q}$.
The splitting property of the antichain~$Q$ means that $Q$~can be split
into two subsets~$Q_1$,~$Q_2$ such that the poset~$P$ is the union of
the upset generated by~$Q_1$ and the downset generated by~$Q_2$. So
any element of~$P$ is either above some element in~$Q_1$ or below some
element in~$Q_2$ (see Figure~\ref{fig:split}).

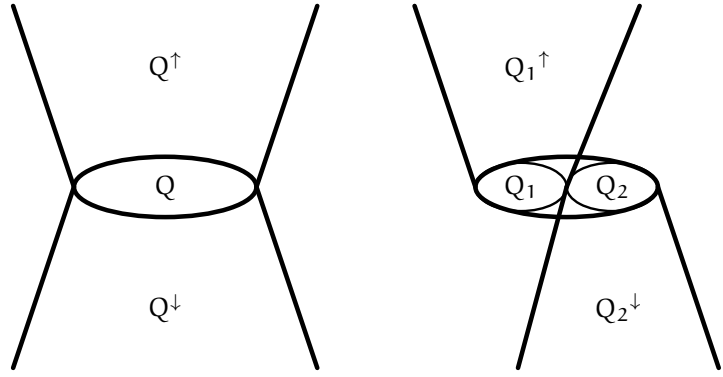
\picture{order.1}{The splitting of an antichain}{split}

A formal definition follows.

\begin{dfn}
We say that a maximal antichain $Q \subset P$
\deph{splits} if there exists a partition of~$S$ into disjoint
subsets $Q_1$ and~$Q_2$ such that $P=\ups{Q_1}\cup \downs{Q_2}$.
In such a case we say that $(Q_1,Q_2)$ is a \deph{splitting}
of the maximal antichain~$Q$.
\end{dfn}

The problem of splitting maximal antichains in the homomorphism
order took on significance when it was observed by Nešetřil and
Tardif~\cite{NesTar:MAC} that except for two small exceptions
finite maximal antichains of size two split in the order of digraphs.
That result is extended here to \Ds s and to finite maximal antichains
of any size; however, the description of exceptions is more involved.

Our approach is direct. In~\ref{pgf:macfd} we define a partition of
any finite maximal antichain and prove that~-- apart from exceptional
cases described later~-- this partition is a splitting of the antichain.

At the end of this section, we suggest an alternative approach that may
lead to a different proof of the splitting property of finite maximal
antichains. It is based on a general condition for the splitting
of antichains in arbitrary posets by Ahlswede, Erdős, Graham and
Soukup~\cite{AhlErdGra:A-splitting,ErdSou:How-to-split}.

There is also a link to homomorphism dualities because a splitting of a
finite maximal antichain is trivially a finite homomorphism duality. In
the case of relational structures with at most two relations, we 
show that the link is stronger: even those finite maximal antichains
that do not split correspond to homomorphism dualities.

For more than two relations this is unknown, but there is a significant
increase in the complexity of the homomorphism order. This suggests that
the property may not hold in this case.

\subsection*{Splitting finite antichains}

We would like to partition a finite maximal antichain~$\Q$ in the
homomorphism order~$\CD$ into disjoint
sets $\F$ and~$\D$ in such a way that $\ups\F\cup\downs\D=\CD$. A
partition is defined in the next paragraph. In the following we 
show that in many cases it satisfies the equality.

\begin{splitting}
\label{pgf:macfd}
Let $\Q=\{Q_1,Q_2,\dotsc,Q_n\}$ be a finite maximal antichain in~$\CD$. Recursively, define the
sets $\F_0$, $\F_1$,~\dots,~$\F_n$ in this way:
\begin{enumerate}
\item Let $\F_0=\emptyset$.
\item For $i=1,2,\dotsc,n$: check whether there exists a \Ds~$X$ satisfying
\begin{itemize}
\item[(i)] $Q_i < X$,
\item[(ii)] $F\notto X$ for any $F\in\F_{i-1}$, and
\item[(iii)] $Q_j\notto X$ for any $j>i$.
\end{itemize}
If such a structure~$X$ exists, let $\F_i=\F_{i-1}\cup\{Q_i\}$, otherwise let $\F_i=\F_{i-1}$.
\item Finally, let $\F=\F_n$ and $\D=\Q\setminus\F$.
\end{enumerate}
\end{splitting}

Because $\Q$ is a maximal antichain, $\ups\Q\cup\downs\Q=\CD$. In addition,
$\ups\F\subseteq\ups\Q$ and $\downs\D\subseteq\downs\Q$ since $\F\subseteq\Q$ and
$\D\subseteq\Q$. Therefore the equality
\[\ups\F\cup\downs\D=\CD,\]
which characterises the splitting of the antichain~$\Q$,
is equivalent to the pair of equalities
\begin{align*}
\ups\F &= \ups\Q,\\
\downs\D &= \downs\Q.
\end{align*}

The following lemma asserts that $\ups\F=\ups\Q$.

\begin{lemma}
\label{lem:qtofto}
Let $\Q$ be a finite maximal antichain and $\F$, $\D$ be defined in~\ref{pgf:macfd}. 
If $Q\in\Q$, $X$~is a \Ds, and $Q<X$, then there exists $F\in\F$ such that $F<X$.
\end{lemma}

\begin{proof}
Among the elements of~$\Q$ that are homomorphic to~$X$, let $Q_i$ be
the element of~$\Q$ with the greatest index~$i$. Then either $F\to X$
for some $F\in\F_{i-1}$, or all the conditions (i), (ii), (iii) are
satisfied and $Q_i\in\F$. So we have found $F\in\F$ such that $F\to X$.

If $F=Q$, then $X\notto F$ by the assumption that $Q<X$. If on the other
hand $F\ne Q$, then the existence of a homomorphism from~$X$ to~$F$ would
imply that $Q\to F$. This is a contradiction because $F$ and~$Q$ are
distinct elements of an antichain. Hence $X\notto F$ and therefore $F<X$.
\end{proof}

To prove that $(\F,\D)$ is a splitting of~$\Q$, it remains to show that
$\downs\D=\downs\Q$.  However, this is not true for all finite maximal
antichains. The following lemma provides a simple description of
antichains for which $(\F,\D)$ is not a splitting.

\begin{lemma}
\label{lem:condY}
Let $\Q$ be a finite maximal antichain and $\F$, $\D$ be defined in~\ref{pgf:macfd}. Then
exactly one of the following conditions holds:
\begin{itemize}
\item[(1)] The pair $(\F,\D)$ is a splitting of~$\Q$.
\item[(2)] There exists a structure~$Y$ such that $Q\notto Y$ for any $Q\in\Q$ and $Y\notto D$
for any $D\in\D$.
\end{itemize}
\end{lemma}

\begin{proof}
If $(\F,\D)$ is a splitting and $Y$~is an arbitrary structure such
that $Q\notto Y$ for any $Q\in\Q$, then $Y\to Q$ for some $Q\in\Q$ because
$\Q$~is a maximal antichain. Moreover, $Q\in\D$ because of splitting.

Conversely, suppose $(\F,\D)$ is not a splitting. So
there exists a structure~$Y$ that violates the definition of a splitting:
$\F\notto Y$ and $Y\notto\D$. Because of Lemma~\ref{lem:qtofto}
we have $\Q\notto Y$.
\end{proof}

If (1) holds, then the antichain~$\Q$ splits. Now we investigate
those maximal antichains that satisfy~(2). The structure~$Y$ has to be
comparable with some element of the maximal antichain~$\Q$, and because
of the condition~(2) there exists $F\in\F$ such that $Y<F$.
We show that all such $Y$'s are bounded from above
by a fairly simple structure.

However, first we need some preparation.
The recursive definition~\ref{pgf:macfd} assures that for every
element~$F$ of~$\F$ there is a witness~$X$ that forces $F$ to be added
to~$\F$. This is formally expressed in the following lemma. The witness
for~$F$ is denoted by~$\check F$.

\begin{lemma}
\label{lem:fhacek}
Let $\Q$ be a finite maximal antichain and $\F$, $\D$ be defined in~\ref{pgf:macfd}. 
For every $F$ in~$\F$ there exists a \Ds~$\check F$ such that $F<\check F$ and moreover $F$~is
the only element of~$\F$ that is homomorphic to~$\check F$.
\end{lemma}

\begin{proof}
The structure~$X$ satisfying the properties (i), (ii), (iii), which
caused $F=Q_i$ to be an element of~$\F_i$ has the required properties
for~$\check F$.
\end{proof}

We use a tool, which is a generalisation of a famous theorem of
Erdős~\cite{Erd:Gtp} (this was one of the first applications of the
then emerging probabilistic method).

\begin{thm}[\cite{ErdHaj:LGH,Lov:LGH,NesRod:LGH}]
\label{thm:girthchrom}
Let $\Delta$ be an arbitrary type, and let $g$ and $k$ be positive
integers. Then there exists a \Ds~$G=G(g,k)$ such that
\begin{itemize}
\item every substructure
of~$G$ induced by at most $g$ vertices is a \Df, and
\item whenever the vertices of~$G$ are coloured by fewer than~$k$ colours,
there exists a colour~$\kappa$ that induces an edge of each kind; that
is, for each kind $i\in I$ there is an edge $e\in R_i(G)$ such that all
vertices of~$e$ have colour~$\kappa$.
\qed
\end{itemize}
\end{thm}

Recall from Definition~\ref{dfn:balanced} that a balanced structure is
a structure that is homomorphic to a forest.

\begin{lemma}
Let $\Q$ be a finite maximal antichain and $\F$, $\D$ be defined in~\ref{pgf:macfd}. 
If $F\in\F$, then $F$~is balanced.
\end{lemma}

\begin{proof}
Let $F\in\F$ be arbitrary and let $\check F$ be the structure whose
existence is guaranteed by Lemma~\ref{lem:fhacek}. Furthermore, let $k:=\max\{|\bs
Q|^{|\bs{\check F}|} : Q\in\Q\}+1$ and let $G$ be
a \Ds\ such that any substructure of~$G$ on at most~$|\bs F|$~vertices is
a \Df\ and whenever the vertices of~$G$ are coloured by fewer than~$k$
colours, there exists a colour that induces an
edge of each kind (Theorem~\ref{thm:girthchrom}).

Consider the structure~$H=\check F\times G$. Suppose that $f:H\to Q$ for
some $Q\in\Q$. For every vertex~$u$ of~$G$, the mapping~$f_u:\bs{\check
F}\to\bs{Q}$ is defined by $f_u(x)=f(u,x)$. We consider this assignment of
mappings to vertices of~$G$ as a colouring of the vertices. Since there
are only $|\bs Q|^{|\bs{\check F}|}<k$ possible mappings from~$\bs{\check
F}$ to~$\bs{Q}$, there exists a colour that induces an edge of every kind; so there exists a
mapping~$g:\bs{\check F}\to\bs{Q}$ satisfying the condition that
\begin{multline}
\label{eq:3.1}
\text{for every $i\in I$ there exists an edge
$(u_1,u_2,\dotsc,u_{\delta_i})$ of~$G$}\\
\text{such that
$f_{u_1}=f_{u_2}=\dotsb=f_{u_{\delta_i}}=g$.}
\end{multline}

Then $g$~is a homomorphism from~$\check F$ to~$Q$: whenever $(x_1,x_2,\dotsc,x_{\delta_i})\in
R_i(\check F)$ for some $i\in I$, we have that
\begin{multline*}
\bigl(g(x_1),g(x_2),\dotsc,g(x_{\delta_i})\bigr)=
\bigl(f_{u_1}(x_1),f_{u_2}(x_2),\dotsc,f_{u_{\delta_i}}(x_{\delta_i})\bigr)\\
=\bigl(f(u_1,x_1),f(u_2,x_2),\dotsc,f(u_{\delta_i},x_{\delta_i})\bigr)\in R_i(Q),
\end{multline*}
because
$\bigl((u_1,x_1),(u_2,x_2),\dotsc,(u_{\delta_i},x_{\delta_i})\bigr)\in R_i(H)$
and $f$~is a homomorphism from~$H$ to~$Q$;
here $(u_1,u_2,\dotsc,u_{\delta_i})$ is the edge of~$G$ from~\eqref{eq:3.1}.
That is a contradiction, because $F<\check F$ and $\Q$~is an antichain containing~$F$.
We conclude that $H\notto Q$ for any $Q\in\Q$. 

By Lemma~\ref{lem:qtofto}
and because $\Q$~is a maximal antichain, there exists $F'\in\F$ that
is homomorphic to~$H=\check F\times G\to\check F$. But $F$~is the only
element of~$\F$ that is homomorphic to~$\check F$, so we have $F'=F$,
and consequently $F\to H\to G$. The image of a homomorphism from~$F$
to~$G$ has no more
than $|\bs F|$~vertices, whence it is a forest. This concludes the proof.
\end{proof}

The elements of~$\F$ are balanced for all finite maximal antichains,
even for those for which $(\F,\D)$ is a splitting.

In the following, we closely investigate the ``non-splitting''
antichains. We derive properties of the structures~$Y$ that violate the
splitting (as in condition~(2) of Lemma~\ref{lem:condY}). In particular,
we show that some paths~-- called \emph{forbidden paths}~-- are not
homomorphic to~$Y$.

\begin{dfn}
Every \Dp\ has a height labelling; we say that a core \Dp~$P$ is a
\deph{forbidden path} if it has two edges of the same kind whose vertices
are not labelled the same. (This property does not depend on what height
labelling we choose, see Proposition~\ref{prop:htlab}.)
\end{dfn}

\begin{lemma}
\label{lem:338}
Let $\Q$ be a finite maximal antichain
in~$\CD$ and let $\F$, $\D$ be defined in~\ref{pgf:macfd}.  If $Y$~is a
\Ds\ such that $Y\notto D$ for any $D\in\D$ and $Y<F$ for some $F\in\F$,
and $P$~is a forbidden path, then $P\notto Y$.
\end{lemma}

\begin{proof}
We may suppose that the two edges of the same kind~$i\in I$ that are not
labelled the same (which prove that $P$~is indeed a forbidden path) are the
end edges of~$P$. Otherwise we could take a subpath of~$P$ (the smallest
substructure of~$P$ that contains both these edges) and show that it is
not homomorphic to~$Y$; consequently $P$~is not homomorphic to~$Y$ either.

Let the two end edges of~$P$ be
$(x_1,x_2,\dotsc,x_{\delta_i}),(y_1,y_2,\dotsc,y_{\delta_i})\in
R_i(P)$. At most one of the vertices $x_1$, $x_{\delta_i}$ is contained
in another edge of~$P$, and so is at most one of~$y_1$,~$y_{\delta_i}$;
if a vertex is contained in only one edge, we call it \deph{free}.

Let $Z$ be a long zigzag: a path with $2m$ or $2m+1$ edges, depending
on the end edges of~$P$.  If both $x_1$ and $y_1$ are free or if both
$x_{\delta_i}$ and $y_{\delta_i}$ are free, we use an even number of
edges; otherwise we use an odd number of edges. All edges of the zigzag are of the same kind as
the end edges of~$P$. 

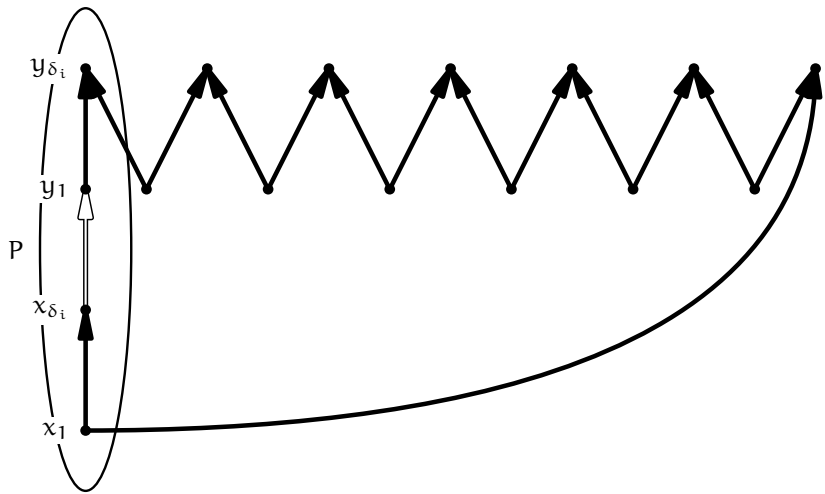
\picture{order.2}{Constructing $Z$ from a forbidden path~$P$}{zigzagZ}

Even though the definition of~$Z$ and~$W$ should be clear from
Figure~\ref{fig:zigzagZ}, we may also define them formally here. Suppose
$x_1$ and $y_1$ are free vertices in~$P$. Then
\begin{align*}
\bs Z &:= \bigl\{1,2,\dotsc,2m(\delta_i-1)+1\bigr\},\\
R_i(Z) &:= \bigl\{(k,k+1,\dotsc,k+\delta_i-1),(k+\delta_i-1,k+\delta_i,\dotsc,k+2\delta_i-2):\\
	&\qquad\quad k=1, 1+2(\delta_i-1), 1+4(\delta_i-1),\dotsc, 1+2(m-1)(\delta_i-1)\bigr\},
\end{align*}
$\approx$ is an equivalence relation on \bs{P+Z} with $x_1\approx 1$, $y_1\approx
2m(\delta_i-1)+1$, and $a\approx a$, and finally
\[W := (P+Z)/{\approx}.\]
If other vertices in the end edges of~$P$ are free, the description is analogous.

Clearly, any proper substructure of~$W$ that does not contain all vertices
of~$Z$ is homomorphic to~$P$. We choose the length of the zigzag (by a
suitable choice of~$m$) in such a way that the number of vertices of~$Z$
is bigger than the number of vertices of any structure in~$\F$.

Now observe that if there exists a height labelling of~$W$, the vertices
in the end edges of~$P$ have the same labels because they are joined by
the zigzag; at the same time, they have distinct labels, because they
are joined by the forbidden path~$P$. Therefore no height labelling
of~$W$ exists and by Proposition~\ref{prop:htlab} the structure~$W$
is not balanced.

Now consider the sum $W+Y$. It is comparable with some some element of
the maximal antichain~$\Q$. However, $W+Y\notto D$ for any $D\in\D$,
because $Y\notto D$ by our assumption on~$Y$; also $W+Y\notto F$ for any
$F\in\F$, because $W$~is not balanced and $F$~is (so $F$~is homomorphic
to a forest, but $W$~is not). Therefore $F\to W+Y$ for some $F\in\F$.

However, because the zigzag~$Z$ was very long, the image of a homomorphism
from~$F$ to~$W+Y$ does not contain all vertices of~$Z$. As we have
observed, therefore $F\to P+Y$. As $F\notto Y$ (by the definition of~$Y$),
necessarily $P\notto Y$.
\end{proof}

In the next lemma we prove that complex structures admit homomorphisms
from forbidden paths, and correspondingly structures that admit no
homomorphisms from forbidden paths are simple.

\begin{lemma}
\label{lem:339}
Let $C$ be a connected \Ds. If no forbidden path is homomorphic to~$C$, then $C$~is homomorphic
to a tree with at most one edge of each kind.
\end{lemma}

\begin{proof}
Suppose that no forbidden path is homomorphic to~$C$. If no height
labelling of~$C$ exists, then there exist vertices $u$ and~$v$
connected by two distinct paths in~$\Sh(C)$ such that counting forward
steps minus backward steps on these paths gives a different result
(see~\ref{pgf:htlabcnt}). Let $B$ be the minimal structure such that
its shadow~$\Sh(B)$ contains both of these paths ($\bs B\subseteq\bs
C$, but it need not be an \emph{induced} substructure; include only
those edges whose shadow edges lie in the two paths). Let us ``unfold''
this structure~$B$, which in a way resembles a cycle: choose an edge~$e$
of~$B$, it intersects two other edges. Then construct a path~$P$: its end
edges are two copies of~$e$, and the middle edges are the remaining
edges of~$B$.

Then $P\to C$ but $P$~is a forbidden path, because the copies of~$e$
get different labels. Therefore $C$~has a height labelling~$\gimel$.

Next observe that any two edges of the same kind are labelled in the
same way. If there were two edges with differently labelled vertices,
there would be a path in~$C$ with these two edges as end edges (because
$C$~is connected) and that would be a forbidden path.

Let $T$ be the structure with base set~$\gimel[\bs C]$, all labels used
on the vertices of~$C$. The edges of~$T$ are such that the identity mapping is a height
labelling of~$T$; in other terms
\begin{multline*}
R_i(T) = \bigl\{(x_1,x_2,\dotsc,x_{\delta_i}) : (x_{j+1})_{(i,j)} = (x_{j})_{(i,j)}+1,\\
	(x_{j+1})_{(i',j')} = (x_{j})_{(i',j')}\text{ for $(i',j')\ne(i,j)$}\bigr\}.
\end{multline*}

Because all the edges of the same kind in~$C$ have the same labelling,
$T$~has at most one edge of each kind. Moreover, any cycle in~$\Sh(T)$
would violate the height labelling of~$T$, so $T$~is a tree. And finally,
the height labelling~$\gimel$ of~$C$ is a homomorphism from~$C$ onto~$T$.
\end{proof}

\begin{cor}
\label{cor:ydstar}
Let $D^*$ be the sum of all \Dt s with at most one edge of each kind. If $Y$~satisfies the
condition~(2) of~\ref{lem:condY}, then $Y\to D^*$.
\end{cor}

\begin{proof}
The claim follows immediately from Lemmas~\ref{lem:338}~and~\ref{lem:339}.
\end{proof}

This shows that the cases when the antichain does not split are very
specific (and one would like to say they are rather rare).

\begin{thm}
\label{thm:aasplit}
Let $\Q$ be a finite maximal antichain in~$\CD$. Let $D^*$ be the sum
of all \Dt s with at most one edge of each kind. Suppose that every
element $Q\in\Q$ has the property that whenever $Y<Q$ and $Y\to D^*$
then there exists a \Ds~$X$ such that $Y<X<Q$ and $X\notto D^*$. Then
the antichain~$\Q$ splits; the pair~$(\F,\D)$ defined in~\ref{pgf:macfd}
is a splitting of~$\Q$.
\end{thm}

\begin{proof}
For the sake of contradiction, suppose that the pair~$(\F,\D)$ is
not a splitting of~$\Q$.  By Lemma~\ref{lem:condY}, there exists a
structure~$Y$ such that $Q\notto Y$ for any $Q\in\Q$ and $Y\notto D$
for any $D\in\D$. Since $\Q$~is a maximal antichain, $Y$~is comparable
with an element~$Q$ of~$\Q$; thus there exists $Q\in\Q$ with $Y<Q$.

Since $Y$~satisfies the condition~(2) of~\ref{lem:condY}, by
Corollary~\ref{cor:ydstar} we have $Y\to D^*$. Thus by assumption there
exists~$X$ such that $Y<X<Q$ and $X\notto D^*$. The structure~$X$ is not
homomorphic to any $D\in\D$, because otherwise we would have $Y\to D$
by composition. Hence $X$~satisfies the condition~(2) of~\ref{lem:condY}
as well, and by Corollary~\ref{cor:ydstar} it is homomorphic to~$D^*$,
a contradiction.
\end{proof}

The assumption on the elements of~$\Q$ posed in the previous theorem means
that no element of~$\Q$ is ``too small''. In particular, it is neither
homomorphic to~$D^*$ nor ``immediately'' above a structure homomorphic
to~$D^*$, that is, there is no gap $(Y,Q)$ such that $Y\to D^*$.

In fact, the assumption can be weakened (and the theorem strengthened)
by requiring that only elements of~$\F$ constructed from~$\Q$
by~\ref{pgf:macfd} have the property that whenever $Y<Q$ and $Y\to
D^*$ then there exists a \Ds~$X$ such that $Y<X<Q$ and $X\notto D^*$.
This is obvious from the proof, where we exploited the property only
for elements of~$\F$.

\subsection*{Connection to finite dualities}

Further examination reveals that in the case of structures with at most
two relations there are no infinite increasing chains below~$D^*$. From
that we can conclude that all elements of~$\F$ are \Df s and thus we
get the following theorem.

\begin{thm}
\label{thm:mac=dual}
Let $\Delta=(\delta_i:i\in I)$ be a type such that $|I|\le2$. Then all finite maximal
antichains in the homomorphism order~$\CD$ are exactly the sets
\begin{equation}
\label{eq:qfd}
\Q=\F\cup\{D\in\D: D\notto F\text{ for any } F\in\F\}
\end{equation}
where $(\F,\D)$~is a finite homomorphism duality.
\end{thm}

\begin{proof}
If $(\F,\D)$ is a finite duality, then for all $F\in\F$ and $D\in\D$
we have $F\notto D$, and so $\Q:=\F\cup\{D\in\D: D\notto F\text{ for
any } F\in\F\}$ is a finite antichain. Moreover, if $F\notto X$ for all
$F\in\F$, then $X$~is homomorphic to some element~$D$ of~$\D$. Either
$D\in\Q$ or $D\to F$ for some $F\in\F$; in any case, $X$~is comparable
with some element of~$\Q$. Hence $\Q$~is a finite maximal antichain.

Conversely, suppose that $\Q$~is a finite maximal antichain in~$\CD$. Let
$(\F,\D)$ be the partition of~$\Q$ defined in~\ref{pgf:macfd}. One of the
conditions of Lemma~\ref{lem:condY} is satisfied. If (1)~is satisfied,
then $(\F,\D)$ is a splitting of~$\F$, and so it is a finite duality
in which no element of~$\D$ is homomorphic to an element of~$\F$ and
$\Q=\F\cup\D$. It remains to examine the case that (2)~is satisfied.

In the case~(2), we first prove that all elements of~$\F$ are
forests. That implies that there exists a finite duality~$(\F,\D')$;
we then prove that this is the duality that satisfies~\eqref{eq:qfd}.

Suppose that condition (2) of Lemma~\ref{lem:condY} is satisfied,
$F\in\F$ and $C$~is a component of~$F$ that is not a tree. By
Lemma~\ref{lem:gap-dual} there is no gap below~$C$. Thus there exist
infinitely many structures $X_1,X_2,\dotsc$ such that $X_1<X_2<\dotsb<C$.
A simple case analysis reveals that the downset~$\downs{D^*}$ generated
by~$D^*$ contains no infinite increasing chain; this is only true for
structures with at most two relations. Hence $C\notto D^*$.

Consequently if $Y\to D^*$ and $Y\to F$ for some structure~$Y$, then
there exists~$X$ such that $Y<X<F$ and $X\notto D^*$ (just like in
Theorem~\ref{thm:aasplit}). It follows that no structure~$Y$ exists such
that $Y\to F$ but $Y\notto D$ for all $D\in\D$.

Now it is time to reuse the trick that served to prove
Lemma~\ref{lem:f5}. As there is no gap below~$C$, we can find a
structure~$B$ such that $B<C$, the structure $B$~is homomorphic to exactly
those components of structures in~$\Q$ as $C$~is, and exactly the same
components of structures in~$\Q$ are homomorphic to~$B$ as to~$C$. Let
$Y$~be the structure constructed from~$F$ by replacing its component~$C$
with~$B$. Clearly $Y\to F$, but $Y$~is homomorphic to no~$D\in\D$, because
$F$~is homomorphic to no~$D\in\D$. This is a contradiction. Therefore
all components of~$F$ are trees, and all elements of~$\F$ are forests.

Invoke the transversal construction on~$\F$ and get a finite
duality~$(\F,\D')$ (remember that $\D$~is defined by splitting the
antichain~$\Q$). We want to prove that $\D$~contains exactly the elements
of~$\D'$ that are not homomorphic to any element of~$\F$.

First, let $D'$ be an element of~$\D'$ such that $D'$~is homomorphic to
no $F\in\F$. As $\Q$~is a maximal antichain and $D'$~is incomparable to
every element of~$\F$, we have $D'\to D$ for some $D\in\D$. If $D\notto
D'$, then some element of~$\F$ is homomorphic to~$D$ by duality, a
contradiction with $\Q$~being an antichain. Hence $D\heq D'$, and so
$D=D'$ since both $D$ and~$D'$ are cores. Therefore $\D$~contains all
elements of~$\D'$ that are not homomorphic to anything in~$\F$.

Finally, we show that $\D$ contains no other elements. Suppose that
$D\in\D$. Because $\Q$~is an antichain, no $F\in\F$ is homomorphic
to~$D$. Thus by duality $D\to D'$ for some $D'\in\D'$.  However,
$D'$~is homomorphic to no~$F\in\F$ (otherwise $D$~would also be), and
so we know from the previous paragraph that $D'\in\D$ and consequently
$D'\in\Q$. Once again using the fact that $\Q$~is an antichain we conclude
that $D'=D$.

In this way, we have found a finite duality $(\F,\D')$ such that
\[
\Q=\F\cup\{D\in\D': D\notto F\text{ for any } F\in\F\}.
\qedhere
\]
\end{proof}

The case of three or more relations ($|I| \geq 3$)  is presently open.
There may be a ``quantum leap'' here as indicated by the following result,
which can be deduced from~\cite{HubNes:Finite}. It implies that for more
than two relations we cannot rely on the fact that the suborder induced
by preimages of~$D^*$ is simple.

\begin{prop}
Let $\Delta=(2,2,2)$. Then the suborder of~$\CD$ induced by all structures
homomorphic to~$D^*$ is a universal countable partial order; that is,
any countable partial order is an induced suborder of this order.
\qed
\end{prop}

\subsection*{Cutting points}

Finally, we show a connection of splitting antichains and cutting points.

\begin{dfn}
Let $P$ be a poset.
An element $y\in P$ is called a \deph{cutting point} if there
are $x, z \in P$ such that $x \prec y \prec z$ and $[x,z] = [x,y] \cup
[y,z]$. (The interval $[x,z] := \{y\in P : x \preceq y \preceq z\}$.)
\end{dfn}

The connection is that every finite maximal antichain
without a cutting point splits.  The following follows
from~\cite[Theorem~2.10]{ErdSou:How-to-split}.

\begin{thm}
If $S$ is a finite maximal antichain that does not contain a cutting
point, then $S$~splits.
\qed
\end{thm}

Thus a characterisation of all cutting points may potentially provide
another proof of the splitting property for finite maximal antichains.
Some cutting points are actually connected to dualities.

\begin{prop}
Let $T$ be a \Dt\ and let $D$ be its dual. Then the \Ds s $T+D$
and~$T\times D$ are cutting points in the homomorphism order~$\CD$.
\end{prop}

\begin{proof}
Consider the interval~$[\bot, T]$, which is equal to the downset generated
by~$T$. Suppose that $X$~is a \Ds\ such that $X<T$. Then $X\to D$,
because $T\notto X$. Thus $X\to T\times D$. Hence the interval~$[T\times
D, T]$ contains only its end-points, that is $[T\times D, T]=\{T\times
D, T\}$. Moreover, $[\bot, T\times D]\cup[T\times D, T]=[\bot, T]$,
so $T\times D$~is a cutting point.

Similarly, if $D<X$, then $T+D\to X$. Hence
$[D,T+D]\cup[T+D,\top]=[D,\top]$ and so $T+D$~is a cutting point.
\end{proof}

However, at present the general problem remains open.

\begin{problem}
Characterise all cutting points in the homomorphism order.
\end{problem}

\setchapterpreamble[u]{%
	\dictum[Winston Churchill]{%
	Out of intense complexities intense simplicities emerge.}}

\chapter{Complexity}
\label{chap:complex}

\bigskip\bigskip

This chapter contains remarks and results on complexity issues; most
of them are implied by the previous chapters. It leaves many questions
unanswered, though.

First, we introduce a generalisation of the \emph{Constraint Satisfaction
Problem (CSP)}, which is a decision problem whether a homomorphism exists
between two relational structures. We consider a parametrisation of CSP,
where the target structure is fixed and the input is the domain. Our
generalisation fixes a finite set~$\HH$ of structures as the parameter
and we ask whether there exists a homomorphism from the input structure
into some structure in the set~$\HH$. We observe that if the set~$\HH$
is the right-hand (dual) side of a finite homomorphism duality, the
problem is solvable by a \pt{} algorithm.

Next we examine the problem of deciding whether an input finite set of
relational structures is the right-hand side of a finite homomorphism
duality. The complexity of this problem with inputs restricted to sets
containing a single structure has recently been determined by Larose,
Loten and Tardif~\cite{Tar:FOD}. Using their result, we are able to
prove that this problem is decidable.

Finally we consider the decision problem whether an input finite set
of relational structures is a maximal antichain in the homomorphism
order. Our characterisation of finite maximal antichains from
Section~\ref{sec:macs} implies that this problem is decidable for
structures with at most two relations. Moreover, we show that the problem
is NP-hard. It is not known at present whether it belongs to the class~NP.

In this chapter, by a \deph{tractable problem} we mean a decision problem
that can be solved by a deterministic Turing machine using a polynomial
amount of computation time, that is a problem belonging to the class~$P$.

\section{Constraint satisfaction problem}
\label{sec:csp}

First, we define the constraint satisfaction problem.

\begin{dfn}
\label{dfn:csp}
Let $H$ be a fixed \Ds\ (called a \deph{template}).  The \deph{constraint
satisfaction problem}~$\CSP(H)$ is the problem to decide for
an input \Ds~$G$ whether there exists a homomorphism $G\to H$.
\end{dfn}

Several classical computational problems can be formulated as constraint
satisfaction, as the following three examples show.

\begin{exa}[$k$-colouring]
\label{dfn:kcol}
Recall from Example~\ref{exa:k-colouring} that a homomorphism to the
complete graph~$K_k$ is the same
as a $k$-colouring of~$G$. Therefore $\CSP(K_k)$ is nothing but $k$-colourability. This is
well-known to be tractable if $k\le 2$ and NP-complete if $k\ge3$.
\end{exa}

\begin{exa}[3-SAT]
The widely known \emph{3-SAT} or \emph{3-satisfiability problem} takes
as its input a propositional formula in conjunctive normal form such
that each clause contains three literals; the question is whether the
input formula is satisfiable,
that is whether logical values can be assigned to its variables in a way
that makes the formula true.

This problem is equivalent to $\CSP(H)$ for the following template~$H$:
let the type $\Delta=(3,3,3,3)$ and let $\bs H=\{0,1\}$. The \Ds~$H$
has four relations, namely $R_0$,~$R_1$, $R_2$, and~$R_3$. Let
\begin{align*}
R_0 &= \bs H^3 \setminus\bigl\{(0,0,0)\bigr\},\\
R_1 &= \bs H^3 \setminus\bigl\{(1,0,0)\bigr\},\\
R_2 &= \bs H^3 \setminus\bigl\{(1,1,0)\bigr\},\text{ and}\\
R_3 &= \bs H^3 \setminus\bigl\{(1,1,1)\bigr\}.
\end{align*}

For an input formula~$\phi$, we construct the \Ds~$G_\phi$ in such a way
that $\bs{G_\phi}$~will be the set for all variables appearing in~$\phi$
and for each clause of~$\phi$ we add a triple to one of the relations:
if there are exactly~$i$~negated literals in the clause, we add a
triple to~$R_i(G_\phi)$ consisting of the three variables appearing
in the clause, with the negated variables first.  For instance, for
the clause $x_1\land\neg x_2\land x_3$ add the triple $(x_2,x_1,x_3)$
to~$R_1(G_\phi)$.

It is straightforward that an edge is preserved by a mapping~$f$
from~$\bs{G_\phi}$ to~$\bs H$ if and only if the corresponding clause is
true in the assignment induced by~$f$ of logical values to variables.
Therefore this assignment makes $\phi$ true if and only if the
mapping~$f$ is a homomorphism from~$G_\phi$ to~$H$, and consequently
$\phi$~is satisfiable if and only if $G$~is homomorphic to~$H$.
\end{exa}

The next example is taken from Tsang's book~\cite{Tsa:Cons}. It is often
used for illustrating algorithms for solving~CSP.

\begin{exa}[$N$-queen problem]
Given any integer~$N$, the problem is to position $N$~queens on
$N$~distinct squares in an $N\times N$ chessboard, in such a way that
no two queens should threaten each other. The rule is such that a queen
can threaten any other pieces on the same row, column or diagonal.

This may be reformulated as the problem of assigning each of $N$~variables
(one for each row) a value from the set $\{1,2,\dotsc,N\}$, marking
the column in which the queen is positioned. It is possible to define
suitable template and relations that ensure that a homomorphism from
the input to the template exists if and only if the input encodes a
non-threatening position of queens on the chessboard.
\end{exa}

Some of the examples are to an extent artificial. However, many
common problems can be formulated as constraint satisfaction very
naturally. These problems appear in numerous areas such as scheduling,
planning, vehicle routing, networks, and bioinformatics. For more details
see, for example, the book~\cite{HCP}.

As the examples show, the complexity of $\CSP(H)$ depends on the
template~$H$. Considerable effort has recently gone into classifying
the complexity of all templates.  This complexity was determined for
undirected graphs by Hell and Nešetřil~\cite{HelNes:Dicho}. However,
already for directed graphs the problem is unsolved. Various results
have led to the following conjecture.

\begin{conj}[\cite{FedVar:SNP}]
\label{fedvarconj}
Let $H$ be a finite relational structure. Then
$\CSP(H)$ is either solvable in polynomial time or NP-complete.
\end{conj}

The following definition is motivated by finite dualities. 

\begin{dfn}
As an analogy to CSP, we define the \deph{generalised constraint
satisfaction problem}~$\GCSP(\HH)$ to be the following decision problem:
given a finite set~$\HH$ of \Ds s, decide for an input \Ds~$G$ whether
there exists $H\in\HH$ such that $G\to H$.
\end{dfn}

The existence of a
finite duality for a template set~$\HH$ ensures the existence of a
\pt\ algorithm for solving the particular generalised
constraint satisfaction problem. This is a classical observation for
the standard constraint satisfaction. We restate it here in view of the
fact that the description of finite dualities is a principal result of
the thesis.

\begin{thm}
If $(\F,\D)$ is a finite homomorphism duality, then $\GCSP(\D)$ is solvable by
a \pt\ algorithm.
\end{thm}

\begin{proof}
The key to the proof is to observe that the algorithm that checks for all possible mappings
from~$\bs F$ to~$\bs G$ whether they are a homomorphism or not runs in time~$O(|G|^{|F|})$,
polynomial in the size of~$G$.

Therefore for an input \Ds~$G$, it is possible to check for all $F\in\F$ whether $F\to G$; if
the response is negative for all~$F$, then $G$~is homomorphic to some~$H$ in~$\HH$, otherwise
it is not. Clearly this testing can be done in time polynomial in the size of the input
structure~$G$.
\end{proof}

As in Conjecture~\ref{fedvarconj}, one could ask whether there is a dichotomy for GCSP.
However, this problem is not very captivating, as the positive answer to the dichotomy
conjecture for CSP would imply a positive answer here as well.

\begin{thm}
Let $\HH$ be a finite nonempty set of pairwise incomparable \Ds s.
\begin{enumerate}
\item If $\CSP(H)$ is tractable for all $H\in\HH$, then $\GCSP(\HH)$ is tractable.
\item If $\CSP(H)$ is NP-complete for some $H\in\HH$, then $\GCSP(\HH)$ is NP-com\-plete.
\end{enumerate}
\end{thm}

\begin{proof}
The first claim is evident. For the second claim, there exists a polynomial reduction of
$\CSP(H)$ to $\GCSP(\HH)$. For an input $G$ of $\CSP(H)$, construct $G+H$ as an input for
$\GCSP(\HH)$. Using the pairwise incomparability of structures in~$\HH$, it is obvious that
$G\to H$ if and only if there exists $H'\in\HH$ such that $G+H\to H'$.
\end{proof}

\section{Deciding finite duality}
\label{sec:decdual}

We are interested in the following decision problem: For an input finite
set~$\D$ of \Ds, determine whether there exists a set~$\F$ of \Ds\
such that $(\F,\D)$~is a finite homomorphism duality.

The complexity of the problem was established by Larose, Loten and
Tar\-dif~\cite{Tar:FOD} in the special case where the input set~$\D$
is a singleton, that is $|\D|=1$.

\begin{thm}[{\cite[Theorem~5.1]{Tar:FOD}}]
\label{thm:fod51}
The problem of determining whether for a relational structure~$D$
there exists a finite set~$\F$ of relational structures such that
$\bigl(\F,\{D\}\bigr)$ is a finite duality, is NP-complete.
\qed
\end{thm}

This special case turns out to be essential. As a consequence we get that
the general problem is decidable, as we show in the rest of this section.

Theorem~\ref{thm:tree-bound} claims that in a duality pair
$\bigl(\F,\{D\}\bigr)$, the diameter of the elements of~$\F$, which
are cores by the definition of finite duality, is at most $n^{n^2}$,
where $n=|\bs D|$. We would like to generate all core trees
with small diameter. It has been established that their number is
finite~\cite[Lemma~2.3]{Tar:FOD}. By modifying the proof in the cited
paper we get a rough recursive estimate for the number of such trees
and the number of their edges.

\begin{lemma}
\label{lem:edgebd}
Let $\Delta=(\delta_i:i\in I)$ be a type and let $s:=|I|$ be the number
of relations and $r:=\max\{\delta_i:i\in I\}$ the maximum arity of
a relation.

Let $\T_d$ be the set of all core trees with a root such that the
distance of any vertex from the root is at most~$d$. Let $t_d:=|\T_d|$
be the number of such trees and let $m_d$~be the maximum number of edges
of a tree in~$\T_d$.

Then
\begin{align*}
t_0 &= 1, &	m_0 &= 0,\\
t_d &\le 2^{sr\cdot t_{d-1}^{r-1}}, &
		m_d &\le sr\cdot t_{d-1}^{r-1}\cdot
			\bigl(1+(r-1)\cdot t_{d-1}\cdot m_{d-1}\bigr).
\end{align*}
\end{lemma}

\begin{proof}
There is exactly one rooted tree with all vertices in distance at most~$0$
from the root: the tree~$\bot$ with one vertex and no edges. So $t_0=1$.
The tree~$\bot$ has no edges, hence $m_0=0$.

For a rooted tree with maximum distance at most~$d$ from the root~$v$,
we can encode every edge $(u_1,u_2,\dotsc,u_{\delta_i})$ that contains
the root by the name of the relation $i\in I$, the index~$k$ such that
$v=u_k$ and the trees rooted at~$u_j$ for $u_j\ne u_k$. Such trees belong
to~$\T_{d-1}$, thus the number of possible labels for edges is at most
$sr\cdot t_{d-1}^{r-1}$. In a core tree, all edges containing the root
have pairwise distinct labels. Hence the number of such trees is
$t_d\le 2^{sr\cdot t_{d-1}^{r-1}}$.

In a tree with maximum distance at most~$d$ from the root~$v$ there
are at most $sr\cdot t_{d-1}^{r-1}$ edges that contain the root. Hence
there are at most $(r-1)\cdot sr\cdot t_{d-1}^{r-1}$ vertices other
than the root~$v$. In each of these vertices, no more than~$t_{d-1}$
trees from~$\T_{d-1}$ are rooted. Therefore the number of edges in the
tree is at most
$sr\cdot t_{d-1}^{r-1} + (r-1)\cdot sr\cdot t_{d-1}^{r-1}\cdot t_{d-1}\cdot m_{d-1}=
sr\cdot t_{d-1}^{r-1}\cdot \bigl(1+(r-1)\cdot t_{d-1}\cdot m_{d-1}\bigr)$.
\end{proof}

As a consequence, in the duality pair~$\bigl(\F,\{D\}\bigr)$ the set~$\F$
is computable from~$D$.

\begin{lemma}
\label{lem:ffromd}
There exists an algorithm that computes a set~$\F$ of \Ds s from an input
\Ds~$D$ so that $\bigl(\F,\{D\}\bigr)$ is a finite homomorphism duality,
provided that such a set exists.
\end{lemma}

\begin{proof}
According to Theorem~\ref{thm:fod51} there is an algorithm that decides
wheth\-er such a set~$\F$ exists.

If it exists, by Theorem~\ref{thm:tree-bound} the diameter of all elements
of~$\F$ is bounded by~$d:=n^{n^2}$, where $n=|\bs D|$. Thus there is a
bound~$m:=m_d$ on the number of edges of the elements of~$\F$
that is computable from Lemma~\ref{lem:edgebd}. Let $\F'$ be the set
of all core \Dt s with at most $m$ edges that are not homomorphic
to~$D$. Theorem~\ref{thm:tree-bound} and Lemma~\ref{lem:edgebd} imply
that for any \Ds~$X$ that is not homomorphic to~$D$ there exists a \Dt~$F$
with at most $m$ edges such that $F\to X$ but $F\notto D$. Hence
\[ {\F'{\notto}} = {{\to}D}.\]
Therefore $\F$ is the set of all homomorphism-minimal elements of~$\F'$.

It is fairly straightforward to design an algorithm for constructing
all core trees with a bounded number of edges, as well as an algorithm
for determining the homomorphism-minimal elements of a finite set of
structures.
\end{proof}

We can conclude that the problem to determine whether an input set of
structures is a right-hand side of a finite duality is decidable.

\begin{thm}
There exists an algorithm that determines whether for an input finite
set~$\D$ of \Ds s there exists a finite set~$\F$ of \Ds s such that
$(\F,\D)$ is a finite homomorphism duality.
\end{thm}

\begin{proof}
The algorithm is as follows:

\begin{enumerate}

\item For each element~$D$ of~$\D$, determine whether it is a core. If
not, then no such duality can exist by definition.

\item For each element~$D$ of~$\D$, determine whether it
is a product of duals, that is whether there exists a finite
set~$\F$ such that $\bigl(\F,\{D\}\bigr)$ is a finite duality (see
Theorem~\ref{thm:fod51}). If an element of~$\D$ is not a product of
duals, then $\D$~is not a right-hand side of a finite duality because
of Theorem~\ref{thm:finite-character}.

\item For each element~$D$ of~$\D$, compute the set~$\F(D)$ such that
$\bigl(\F,\{D\}\bigr)$ is a finite duality (see Lemma~\ref{lem:ffromd}).

\item Check whether $\{\F_D:D\in\D\}$ is the set of all transversals
for some finite set~$\F$ of core \Df s. This can be done greedily by
considering all sets of core \Df s whose components appear in the
sets~$\F_D$, and by constructing the transversals (directly from
Definition~\ref{dfn:trans}).

\end{enumerate}

Theorem~\ref{thm:finite-character} implies that if the algorithm finds
a set~$\F$ in the last step, then $(\F,\D)$ is a finite homomorphism
duality; otherwise $\D$~is the right-hand side of no finite duality.
\end{proof}

\section{Deciding maximal antichains}
\label{sec:MAC-decide}

We consider the problem of deciding whether an input finite set of
relational structures forms a finite maximal antichain in the homomorphism
order. The problem is called the \emph{MAC problem}; the letters MAC
stand for ``\textsc{m}aximal \textsc{a}nti\textsc{c}hain''.

\begin{dfn}
The \deph{MAC problem} is to decide whether an input finite non-empty
set~$\Q$ of \Ds s is a maximal antichain in the homomorphism order~$\CD$.
\end{dfn}

The characterisation of finite maximal antichains for types with at most
two relations (Theorem~\ref{thm:mac=dual}) implies decidability of the
MAC problem.

\begin{thm}
\label{thm:mac-deci}
Let $\Delta=(\delta_i:i\in I)$ be a type such that $|I|\le2$. Then the
MAC problem is decidable.
\end{thm}

\begin{proof}
The algorithm is as follows:

\begin{enumerate}
\item
For each element of~$\Q$, check whether its core is a forest. The core
of a \Ds\ is computable (by checking the existence of a retraction to
every substructure). Deciding whether a \Ds\ is a forest is possible
even in polynomial time.

\item
Let $\F\subseteq\Q$ be the set of all such structures. Find
all transversals over~$\F$. This can be done directly from
Definition~\ref{dfn:trans}.

\item
For each transversal~$\M$, construct its dual~$D(\M)$. First use the
bear construction (\ref{pgf:bear-def}) to construct the dual of each
element of~$\M$ and then take the product of all these duals.

\item
Check whether $\Q\setminus\F$ is formed exactly by structures
homomorphically equivalent to the duals of transversals constructed in
the previous step.
\end{enumerate}

Theorem~\ref{thm:mac=dual} implies that the algorithm is correct.
\end{proof}

\begin{thm}
Let $\Delta=(\delta_i:i\in I)$ be a type such that $|I|\le2$. Then the
MAC problem is NP-hard.
\end{thm}

\begin{proof}
We will use the fact that for any type $\Delta$ there exists a \Dt~$T$
such that $\CSP(T)$ is NP-complete. We construct the following reduction
of $\CSP(T)$ to the MAC problem: For an input structure~$G$ of $\CSP(T)$,
let $\Q(G):=\bigl\{G+T, D(T)\bigr\}$. The set~$\Q(G)$ can be constructed
from~$G$ in polynomial time. By Theorem~\ref{thm:mac=dual}, $\Q(G)$~is
a finite maximal antichain if and only if $G\to T$.
\end{proof}

However, the algorithm given in the proof of Theorem~\ref{thm:mac-deci}
does not ensure that the MAC problem is in the class NP. This is not
known at present.

\begin{problem}
Is the MAC problem in NP?
\end{problem}

The hard part of the problem may actually consist in
finding the cores of the involved structures. Also, in our proof of
NP-hardness we actually reduce the decision whether the core of~${G+T}$
is~$T$ to the MAC problem. So it makes sense to ask whether the complexity
of the MAC problem changes when inputs are restricted to cores.

\begin{problem}
What is the complexity of the MAC problem if input is restricted to sets
of cores?
\end{problem}

\setbibpreamble{%
	\dictum[Mark Twain]{%
A classic is a book which people praise and don't read.}
\bigskip\bigskip}


\begin{thebibliography}{10}

\bibitem{AHS:AbsConCat}
J.~Ad{\'a}mek, H.~Herrlich, and G.~E. Strecker.
\newblock {\em Abstract and Concrete Categories. The Joy of Cats}.
\newblock John Wiley and Sons, 1990.

\bibitem{AhlErdGra:A-splitting}
R.~Ahlswede, P.~L. Erd{\H o}s, and N.~Graham.
\newblock A splitting property of maximal antichains.
\newblock {\em Combinatorica}, 15(4):475--480, 1995.

\bibitem{BarWel:Cat}
M.~Barr and C.~Wells.
\newblock {\em Category Theory for Computing Science}.
\newblock Les Publications CRM, Montr{\'e}al, 3rd edition, 1999.

\bibitem{DavPri:Introduction}
B.~A. Davey and H.~A. Priestley.
\newblock {\em Introduction to Lattices and Order}.
\newblock Cambridge University Press, 2nd edition, 2002.

\bibitem{Edm:PTF}
J.~Edmonds.
\newblock Paths, trees, and flowers.
\newblock {\em Canad. J. Math.}, 17:449--467, 1965.

\bibitem{EZaSau:Chro}
M.~El-Zahar and N.~Sauer.
\newblock The chromatic number of the product of two 4-chromatic graphs is 4.
\newblock {\em Combinatorica}, 5(2):121--126, 1985.

\bibitem{Erd:Gtp}
P.~Erd{\H o}s.
\newblock Graph theory and probability.
\newblock {\em Canad. J. Math.}, 11:34--38, 1959.

\bibitem{ErdHaj:LGH}
P.~Erd{\H o}s and A.~Hajnal.
\newblock On chromatic number of graphs and set-systems.
\newblock {\em Acta Math. Hungar.}, 17({1--2}):61--99, 1966.

\bibitem{ErdSou:How-to-split}
P.~L. Erd{\H o}s and L.~Soukup.
\newblock How to split antichains in infinite posets.
\newblock {\em Combinatorica}, 27(2):147--161, 2007.

\bibitem{FedVar:SNP}
T.~Feder and M.~Y. Vardi.
\newblock The computational structure of monotone monadic {SNP} and constraint
  satisfaction: {A} study through {D}atalog and group theory.
\newblock {\em SIAM J. Comput.}, 28(1):57--104, 1998.

\bibitem{FNT:WG06}
J.~Foniok, J.~Ne{\v s}et{\v r}il, and C.~Tardif.
\newblock Generalised dualities and finite maximal antichains.
\newblock In F.~V. Fomin, editor, {\em Graph-Theoretic Concepts in Computer
  Science (Proceedings of WG 2006)}, volume 4271 of {\em Lecture Notes in
  Comput. Sci.}, pages 27--36. Springer-Verlag, 2006.

\bibitem{FNT:Eurocomb7}
J.~Foniok, J.~Ne{\v s}et{\v r}il, and C.~Tardif.
\newblock On finite maximal antichains in the homomorphism order.
\newblock {\em Electron. Notes Discrete Math.}, 29:389--396, 2007.

\bibitem{FNT:GenDu}
J.~Foniok, J.~Ne{\v s}et{\v r}il, and C.~Tardif.
\newblock Generalised dualities and maximal finite antichains in the
  homomorphism order of relational structures.
\newblock {\em European J. Combin.}, to appear.

\bibitem{HHMN:Mult}
R.~H{\"a}ggkvist, P.~Hell, D.~J. Miller, and V.~Neumann~Lara.
\newblock On multiplicative graphs and the product conjecture.
\newblock {\em Combinatorica}, 8(1):63--74, 1988.

\bibitem{Hed:Con}
S.~T. Hedetniemi.
\newblock Homomorphisms of graphs and automata.
\newblock University of Michigan Technical Report 03105-44-T, University of
  Michigan, 1966.

\bibitem{Hed:Uni}
Z.~Hedrl{\'\i}n.
\newblock On universal partly ordered sets and classes.
\newblock {\em J. Algebra}, 11(4):503--509, 1969.

\bibitem{HelNes:Dicho}
P.~Hell and J.~Ne{\v s}et{\v r}il.
\newblock On the complexity of {$H$}-coloring.
\newblock {\em J. Combin. Theory Ser. B}, 48(1):92--119, 1992.

\bibitem{HelNes:GrH}
P.~Hell and J.~Ne{\v s}et{\v r}il.
\newblock {\em Graphs and Homomorphisms}, volume~28 of {\em Oxford Lecture
  Series in Mathematics and Its Applications}.
\newblock Oxford University Press, 2004.

\bibitem{HubNes:Finite}
J.~Hubi{\v c}ka and J.~Ne{\v s}et{\v r}il.
\newblock Finite paths are universal.
\newblock {\em Order}, 22(1):21--40, 2005.

\bibitem{HuNe:Paths}
J.~Hubi{\v c}ka and J.~Ne{\v s}et{\v r}il.
\newblock Universal partial order represented by means of oriented trees and
  other simple graphs.
\newblock {\em European J. Combin.}, 26(5):765--778, 2005.

\bibitem{Kom:Somenew}
P.~Kom{\'a}rek.
\newblock Some new good characterizations for directed graphs.
\newblock {\em \v Casopis P\v est. Mat.}, 109(4):348--354, 1984.

\bibitem{Kom:Phd}
P.~Kom{\'a}rek.
\newblock {\em Good characterisations in the class of oriented graphs}.
\newblock PhD thesis, Czechoslovak Academy of Sciences, Prague, 1987.
\newblock In Czech (Dobr{\'e} charakteristiky ve t{\v r}{\'\i}d{\v e}
  orientovan{\'y}ch graf{\r u}).

\bibitem{Tar:FOD}
B.~Larose, C.~Loten, and C.~Tardif.
\newblock A characterisation of first-order constraint satisfaction problems.
\newblock In {\em Proceedings of the 21st IEEE Symposium on Logic in Computer
  Science (LICS'06)}, pages 201--210. IEEE Computer Society, 2006.

\bibitem{Lov:LGH}
L.~Lov{\'a}sz.
\newblock On chromatic number of finite set-systems.
\newblock {\em Acta Math. Hungar.}, 19({1--2}):59--67, 1968.

\bibitem{N:TG}
J.~Ne{\v s}et{\v r}il.
\newblock {\em Theory of Graphs}.
\newblock SNTL, Prague, 1979.
\newblock In Czech (Teorie graf{\r u}).

\bibitem{Nes:ColPos}
J.~Ne{\v s}et{\v r}il.
\newblock The coloring poset and its on-line universality.
\newblock KAM-DIMATIA Series 2000-458, Charles University, Prague, 2000.

\bibitem{Nes:SS}
J.~Ne{\v s}et{\v r}il.
\newblock Combinatorics of mappings.
\newblock KAM-DIMATIA Series 2000-472, Charles University, Prague, 2000.

\bibitem{NesOss:ResDual}
J.~Ne{\v s}et{\v r}il and P.~Ossona~de Mendez.
\newblock Grad and classes with bounded expansion {III}. {R}estricted
  dualities.
\newblock KAM-DIMATIA Series 2005-741, Charles University, Prague, 2005.

\bibitem{NesPul:SubFac}
J.~Ne{\v s}et{\v r}il and A.~Pultr.
\newblock On classes of relations and graphs determined by subobjects and
  factorobjects.
\newblock {\em Discrete Math.}, 22(3):287--300, 1978.

\bibitem{NesPulTar:HeytDual}
J.~Ne{\v s}et{\v r}il, A.~Pultr, and C.~Tardif.
\newblock Gaps and dualities in {H}eyting categories.
\newblock {\em Comment. Math. Univ. Carolin.}, 48(1):9--23, 2007.

\bibitem{NesRod:LGH}
J.~Ne{\v s}et{\v r}il and V.~R{\"o}dl.
\newblock A short proof of the existence of highly chromatic hypergraphs
  without short cycles.
\newblock {\em J. Combin. Theory Ser. B}, 27(2):225--227, 1979.

\bibitem{NesIda:Diam}
J.~Ne{\v s}et{\v r}il and I.~{\v S}vejdarov{\'a}.
\newblock Diameters of duals are linear.
\newblock KAM-DIMATIA Series 2005-729, Charles University, Prague, 2005.

\bibitem{NesTar:Dual}
J.~Ne{\v s}et{\v r}il and C.~Tardif.
\newblock Duality theorems for finite structures (characterising gaps and good
  characterisations).
\newblock {\em J. Combin. Theory Ser. B}, 80(1):80--97, 2000.

\bibitem{NesTar:MAC}
J.~Ne{\v s}et{\v r}il and C.~Tardif.
\newblock On maximal finite antichains in the homomorphism order of directed
  graphs.
\newblock {\em Discuss. Math. Graph Theory}, 23(2):325--332, 2003.

\bibitem{NesTar:Short}
J.~Ne{\v s}et{\v r}il and C.~Tardif.
\newblock Short answers to exponentially long questions: Extremal aspects of
  homomorphism duality.
\newblock {\em SIAM J. Discrete Math.}, 19(4):914--920, 2005.

\bibitem{NesTar:A-dualistic}
J.~Ne{\v s}et{\v r}il and C.~Tardif.
\newblock A dualistic approach to bounding the chromatic number of a graph.
\newblock {\em European J. Combin.}, to appear.

\bibitem{NesZhu:Path}
J.~Ne{\v s}et{\v r}il and X.~Zhu.
\newblock Path homomorphisms.
\newblock {\em Math. Proc. Cambridge Philos. Soc.}, 120:207--220, 1996.

\bibitem{PulTrn:Cat}
A.~Pultr and V.~Trnkov{\'a}.
\newblock {\em Combinatorial, Algebraic and Topological Representations of
  Groups, Semigroups and Categories}, volume~22 of {\em North-Holland
  Mathematical Library}.
\newblock North-Holland, Amsterdam, 1980.

\bibitem{HCP}
F.~Rossi, P.~van Beek, and T.~Walsh, editors.
\newblock {\em Handbook of Constraint Programming}, volume~2 of {\em
  Foundations of Artificial Intelligence}.
\newblock Elsevier, 2006.

\bibitem{Tar:Mul}
C.~Tardif.
\newblock Multiplicative graphs and semi-lattice endomorphisms in the category
  of graphs.
\newblock {\em J. Combin. Theory Ser. B}, 95(2):338--345, 2005.

\bibitem{Tsa:Cons}
E.~Tsang.
\newblock {\em Foundations of Constraint Satisfaction}.
\newblock Academic Press, 1993.

\bibitem{Wel:Dense}
E.~Welzl.
\newblock Color families are dense.
\newblock {\em Theoret. Comput. Sci.}, 17(1):29--41, 1982.

\end{thebibliography}
\end{document}